[1994/12/01]
\documentclass[draft,12pt,oneside]{book}
\title{Consistency Strengths of Modified Maximality Principles}
\author{George Leibman}

\usepackage{amsmath,amsthm,amsfonts,amssymb,bm}
\usepackage{bbm}
\usepackage{newlfont}
\usepackage{ifthen}

\chardef\bslash=`\\ 

\newcommand{\beginThm}{\begin{thm}\setcounter{claim}{0}}
\newtheorem{thm}{Theorem}[chapter]
\newtheorem{cor}[thm]{Corollary}
\newtheorem{lem}[thm]{Lemma}

\newtheorem{claim}{Claim}
\newtheorem{con}{Conjecture}

\theoremstyle{definition}

\theoremstyle{remark}

\newcommand{\thmref}[1]{Theorem~\ref{#1}}

\newcommand{\lemref}[1]{Lemma~\ref{#1}}

\newcommand{\corref}[1]{Corollary~\ref{#1}}

\DeclareMathOperator{\cof}{cof} \DeclareMathOperator{\supp}{supp}
\DeclareMathOperator{\coll}{coll}

\newcommand{\eval}[2][\right]{\relax
  \ifx#1\right\relax \left.\fi#2#1\rvert}

\newcommand{\UnderTilde}[1]{{\setbox1=\hbox{$#1$}\baselineskip=0pt\vtop{\hbox{$#1$}\hbox to\wd1{\hfil$\sim$\hfil}}}{}}
\newcommand{\Undertilde}[1]{{\setbox1=\hbox{$#1$}\baselineskip=0pt\vtop{\hbox{$#1$}\hbox to\wd1{\hfil$\scriptstyle\sim$\hfil}}}{}}
\newcommand{\undertilde}[1]{{\setbox1=\hbox{$#1$}\baselineskip=0pt\vtop{\hbox{$#1$}\hbox to\wd1{\hfil$\scriptscriptstyle\sim$\hfil}}}{}}
\newcommand{\UnderdTilde}[1]{{\setbox1=\hbox{$#1$}\baselineskip=0pt\vtop{\hbox{$#1$}\hbox to\wd1{\hfil$\approx$\hfil}}}{}}
\newcommand{\Underdtilde}[1]{{\setbox1=\hbox{$#1$}\baselineskip=0pt\vtop{\hbox{$#1$}\hbox to\wd1{\hfil\scriptsize$\approx$\hfil}}}{}}

\newcommand{\CARD}{\hbox {\sc card}}
\newcommand{\CCC}{\hbox {\sc ccc}}
\newcommand{\PROPER}{\hbox {\sc proper}}
\newcommand{\SEMIPROPER}{\hbox {\sc semiproper}}
\newcommand{\COHEN}{\hbox{\sc cohen}}
\newcommand{\COLL}{\hbox {\sc coll}}
\newcommand{\MP}{\hbox {\sc mp}}
\newcommand{\ZFC}{\hbox {\sc zfc}}
\newcommand{\CH}{\hbox {\sc ch}}

\newcommand{\MPmod}[1]{\MP_{\hbox  {\scriptsize {\textrm #1}}}}
\newcommand{\bval}[1]{[{\kern -.15 em}[#1]{\kern -.15 em}]}

\newcounter{listcounter}

\newcommand\makeabstractpage[2]{
 \newpage
  \vskip 3.0em%
  \begin{center}%
    Abstract
  \end{center}
  \vskip 3.0em%

  #2
}

\begin{document}

\frontmatter

\maketitle



\makeabstractpage{Joel David Hamkins}{The Maximality Principle
\MP\ is a scheme which states that if a sentence of the language
of \ZFC\ is true in some forcing extension $V^\mathbb{P}$, and
remains true in any further forcing extension of $V^\mathbb{P}$,
then it is true in all forcing extensions of $V$. A modified
maximality principle $\MP_\Gamma$ arises when considering forcing
with a particular class $\Gamma$ of forcing notions.  A
parametrized form of such a principle, $\MP_\Gamma(X)$, considers
formulas taking parameters; to avoid inconsistency such parameters
must be restricted to a specific set $X$ which depends on the
forcing class $\Gamma$ being considered. A stronger
\textit{necessary} form of such a principle, $\Box\MP_\Gamma(X)$,
occurs when it continues to be true in all $\Gamma$ forcing
extensions.

This study uses iterated forcing, modal logic, and other
techniques to establish consistency strengths for various modified
maximality principles restricted to various forcing classes,
including \CCC, \COHEN, \COLL\ (the forcing notions that collapse
ordinals to $\omega$), $<\kappa$ directed closed forcing notions,
etc., both with and without parameter sets. Necessary forms of
these principles are also considered.}





\tableofcontents

\mainmatter

\chapter*{Introduction}
\addcontentsline{toc}{chapter}{\numberline{}Introduction}
\markboth{Introduction}{INTRODUCTION}

The \textbf{Maximality Principle} (\MP) states that if a sentence
of \ZFC\ is true in some forcing extension $V^{\mathbb{P}}$ of $V$
and remains true in any subsequent forcing extension of
$V^{\mathbb{P}}$, then it is true in $V$.  An equivalent form of
this principle says that the given sentence must then be true in
all forcing extensions of $V$. This principle (or any of its
variations) can play a role analogous to that of other forcing
axioms such as MA or PFA in deciding questions which are left
unanswered by \ZFC\ alone (for example, concerning the size of the
continuum).  It is a \textbf{dense-set free} forcing axiom---it
doesn't require that a filter meet some family of dense sets;
rather, it is expressed in terms of the forcing relation itself.
The maximality principle is discussed extensively in
\cite{ham:max}, and some variations are explored in
\cite{sta:vaa}.

The principle \MP\ can be cast in the language of modal logic by
regarding models of \ZFC\ as possible worlds and defining one such
world to be accessible from another if the first is a forcing
extension of the second.  (This relation is transitive and
reflexive, giving rise to the $S4$ axiom system of modal logic.)
This gives interpretations of the modal concepts of possibility
and necessity.  A statement $\phi$ is possible ($\Diamond\phi$ in
the notation of modal logic) if it is true in some forcing
extension and necessary (denoted by $\Box\phi$) if it is true in
every forcing extension.  A sentence is \textbf{possibly
necessary}, or \textbf{forceably necessary} if it is true in some
forcing extension and remains true in any subsequent forcing
extension. So \MP\ says that if a sentence is forceably necessary
then it is necessary. It is not a formula of \ZFC: it is a scheme,
the collection of all instances of the statement
``$\Diamond\Box\phi $ implies $\Box\phi$" (or its $S4$ equivalent,
``$\Diamond\Box\phi $ implies $\phi$") where $\phi$ ranges over
all sentences in the language of \ZFC.

The Maximality Principle is equiconsistent with \ZFC\
(\cite{ham:max}), but it has two basic variations of greater
consistency strength.  The first is a variation in which any
formula $\phi$ occurring in an instance of the scheme is allowed
to take parameters.  This variation was denoted
$\Undertilde{\mathrm{MP}}$ (``boldface" \MP) in \cite{ham:max}, in
which it was shown that the most general class of parameters which
does not render the principle patently false is $H(\omega_1)$, the
set of hereditarily countable sets.  It is easy to see that such a
principle with parameters allowed from a larger set would be
false.  Indeed, a hereditary uncountable set can be forced to have
countable transitive closure in any further extension, so if it
were allowed as a parameter in \MP, it would already be countable,
contradicting the assertion of its uncountability. Since any
hereditarily countable set can be coded by a real, we will denote
this variation of \MP\ by $\MP(\mathbb{R})$, which expresses the
idea that any statement in the language of set theory with
arbitrary real parameters that is forceably necessary is
necessary.  In exploring modifications to the maximality principle
the new notation, $\MP(X)$ for $\MP$ with parameter set $X$, is
more perspicuous when discussing how different modifications of
\MP\ lead to different natural classes of parameters and on the
other hand, how changing the class of parameters while keeping the
same class of forcing notions leads to a different maximality
principle, possibly of different consistency strength.

The third and final basic form, of greatest consistency strength,
is the Necessary Maximal Principle, or $\Box \MP(\mathbb{R})$,
which says that $\MP(\mathbb{R})$ is necessary --- it holds in $V$
and every forcing extension $V^{\mathbb{P}}$.  This principle has
greater strength (still unknown exactly, but it has been bounded
above and below by consistency strengths that put it well up in
the hierarchy, along with Woodin cardinals and Projective
Determinacy, see \cite{ham:max}).  This is because parameters of
the formula $\phi$ are not restricted to those with interpretation
in $V$; they may appear as a result of the forcing that produces
the model in which $\MP(\mathbb{R})$ is interpreted and remains
necessarily true.

\textbf{Modified maximality principles} arise if one considers
possible worlds (models of \ZFC) to be accessible only if they are
extensions obtained by forcing with forcing notions of specified
classes. If $\Gamma$ is such a restricted class of forcing
notions, such an accessibility relation can be expressed between
models $M_1$ and $M_2$ by saying that $M_2$ is a
\textbf{$\Gamma$-forcing extension} of $M_1$ (that is, that $M_2 =
M_1[G]$ where $G$ is $M_1$-generic over some $\mathbb{P}$ where
$\mathbb{P} \in \Gamma$).  A formula $\phi$ is
\textbf{$\Gamma$-necessary} (denoted by $\Box_\Gamma \phi$) when
it is true in all $\Gamma$-forcing extensions.  (A related
concept, used in \cite{sta:vaa}, defines $\phi$ to be
\textbf{$\Gamma$-persistent} if, whenever $\phi$ is true, it is
then true in all $\Gamma$-forcing extensions.  So if $\phi$ is
$\Gamma$-persistent it may not be true, in which case it need not
be $\Gamma$-necessary. A formula $\phi$ is
\textbf{$\Gamma$-forceable} (denoted by $\Diamond_\Gamma \phi$)
when it is true in some $\Gamma$-forcing extension.  (Note that,
as in other interpretations of possibility and necessity, that
$\Diamond$ and $\Box$ are dual to each other---$\Diamond\phi$ can
be defined as $\neg\Box\neg\phi$.)  A formula $\phi$ is
\textbf{$\Gamma$-possibly necessary} or \textbf{$\Gamma$-forceably
necessary} (denoted by $\Diamond_\Gamma \Box_\Gamma \phi$) if it
is $\Gamma$-forceable that $\phi$ is $\Gamma$-necessary.  (The
determination as to whether a forcing notion $\mathbb{P}$ in a
$\Gamma$-forcing extension is itself in $\Gamma$ is made
\textit{de dicto}.  That is, $\Gamma$ is a definable class in
\ZFC, as as such, its members are formally determined by
satisfying its defining formula.  So this formula is interpreted
in the model of \ZFC\ in which $\mathbb{P}$ will be forced with.)
The modified maximality principle that arises with these notions
is denoted $\MP_\Gamma$, which says that if a formula of \ZFC\ is
true in some $\Gamma$-forcing extension $V^{\mathbb{P}}$, (that
is, where $\mathbb{P} \in \Gamma$) of $V$ and remains true in any
subsequent $\Gamma$-forcing extension of $V^{\mathbb{P}}$, then it
is $\Gamma$-necessary in $V$ (hence true, if $\Gamma$ includes the
trivial forcing notion $\{\emptyset\}$, as it almost always will).
In terms of the symbols just introduced, the principle
$\MP_\Gamma$ can be expressed as the scheme
$\Diamond_\Gamma\Box_\Gamma\phi \Longrightarrow \Box_\Gamma\phi$
where $\phi$ can be any statement in the language of \ZFC.

If a modified maximality principle is allowed to take parameters,
then the notation indicated above will be used to indicate a
modified maximality principle with parameters from the class
specified in parentheses. Thus, a formula is
\textbf{$\CCC$-necessary} if it is true in all $\CCC$-forcing
extensions.  The principle $\MPmod{\CCC}(H(2^{\omega}))$ then says
that any formula with a parameter from $H(2^{\omega})$---the sets
hereditarily of cardinality less than $2^\omega$---that is
$\CCC$-forceably necessary is $\CCC$-necessary. It is a result of
\cite{ham:max} that $\MPmod{\CCC}$ with parameters is consistent
(relative to the L\'evy Scheme) only as long as the parameters are
taken from $H(2^{\omega})$, so this is a \textbf{natural class} of
parameters for this particular principle.  This means that two
conditions are met. First, that any instance of the principle with
any parameter in this class is consistent, and second, that that
any instance of the principle with any parameter not in this class
is inconsistent. In what follows, the first order of business upon
introducing a modified maximality principle will be to determine
the natural classes of parameters for which it is consistent.  The
class of $\CCC$ forcing notions itself will be simply denoted by
the abbreviation $\CCC$, and other classes of forcing notions will
be similarly abbreviated.

Modifications to the Necessary Maximality Principle will be
denoted using the same system.  For example $\Box
\MPmod{\CCC}(H(2^\omega))$ says that $\MPmod{\CCC}(H(2^\omega))$
is \CCC-necessary (the box operator is implicitly restricted to
\CCC\ in this notation). This way of forming new principles can be
applied equally to any other restricted type of forcing.

In this work I will try to answer questions concerning the
consistency strengths of these modified maximality principles. The
work is organized around the various classes of forcing notions
(forcing notions) that occur in various modifications of the
maximality principle.  After a chapter on general concepts,
Chapter 2 is devoted to classes of forcing notions that preserve
specified statements in the language of \ZFC, chapter 3 to the
class of $\CCC$ forcing notions, chapter 4 to Cohen forcing,
chapter 5 to forcing notions that collapse cardinals to $\omega$,
and chapter 6 to forcing notions parametrized by large cardinals.
A goal of this work has been to establish the consistency strength
of any given maximality principle with the most extensive class of
parameters that results in a consistent principle.

\chapter{Fundamental notions}
\section{Modal logic of forcing extensions}

\subsection{Set-theoretic semantics} The modal logic of forcing
extensions is a notational convenience for proving theorems in
\ZFC. (A good introduction to the formalism of modal logic used
here is \cite{fm:mod}.)  The notation is called  modal logic
because the symbols used, $\Box$ and $\Diamond$, obey certain
modal logic axioms when interpreted in the forcing context.  In
addition,regarding models of \ZFC\ as possible worlds corresponds
well with Kripke model semantics.  For any formula $\phi$,
$\Box\phi$ says that $\phi$ is true in all forcing extensions: for
all $\mathbb{P}$, $V^\mathbb{P}\models\phi$. $\Diamond\phi$ says
that $\phi$ is true in some forcing extensions.  As stated in the
introduction, this notation can be relativized to specific classes
of forcing notions.  This notation adds defined symbols to the
language of \ZFC. But saying $\phi$ is true in a forcing extension
is already an abbreviation of a sentence in the language of \ZFC\
that involves an abstract ``forcing language" using
$\mathbb{P}$-names based on some partial order $\mathbb{P}$. So in
adding modal symbols, we are just continuing a tradition of adding
another degree of separation from the primitive language
$\{\in\}$. In order to benefit from the ability to think
intuitively, using and mining familiar concepts founded in Kripke
models, we may regard forcing extensions as accessible possible
worlds. But while harvesting intuitions provided by the modal
symbols' interpretations as necessity and possibility, one should
remain centered in the realization that these are just statements
in the language of \ZFC\ being interpreted in a single model,
replete with names for sets in all forcing extensions.

For us, it is primarily the economy of expression of maximality
principles and the arguments they require that leads to the use of
these modal symbols.  Recall that their definitions given in the
introduction are in fact set-theoretic.  Thus $\Box$, $\Diamond$
and their relativized counterparts are defined symbols added to
the language of \ZFC; they can be safely eliminated from any
formal argument in which they are used.

An equivalent definition of $\Box$ and $\Diamond$ can be made in
terms of complete Boolean algebras, a commonly used foundation for
forcing arguments.  For a complete discussion of this approach see
\cite{bel:bvm}. If one is working in a Boolean-valued model
$V^\mathcal{B}$ of \ZFC, with Boolean values in the Boolean
algebra $\mathcal{B}$, the Boolean value of a formula $\phi$ is
denoted $\bval{\phi}^\mathcal{B}$.  If $\Gamma$ is a class of
complete Boolean algebras, we have dual definitions:
$\Diamond_\Gamma\phi$ if and only if there exists some
$\mathcal{B}$ in $\Gamma$ such that $\bval{\phi}^\mathcal{B}\neq
0$, and $\Box_\Gamma\phi$ if and only if for all $\mathcal{B}$ in
$\Gamma$, $\bval{\phi}^\mathcal{B}=1$.  All results herein will
follow from this definition, including results on iterated
forcing, by translating forcing notions to the corresponding
regular open algebras.  In fact, a Boolean-valued model approach
to forcing has a distinct advantage.  It allows one to formalize
the concept of ``forcing over $V$" by establishing a relationship
between $V$ and $V^\mathcal{B}$.  Since this occurs within \ZFC,
it applies over any model.  This frees one from having to base a
definition of forcing on countable transitive models.

\subsection{S4 forcing classes} Given this set-theoretic
definition of the modal operators, we can ask if the operators in
fact behave according to various axioms of modal logic.  In our
interpretation, an axiom of modal logic will be a scheme in \ZFC\
holding true for any formula of \ZFC\ when substituted into the
scheme.  If these schemes are true in \ZFC, we can conveniently
use them to prove assertions regarding maximality principles.  The
relevant axioms for our purposes will be the \textbf{K},
\textbf{S4} and \textbf{S5} axiom schemes. The \textbf{K} axiom
schemes are
\begin{list}
 {(\arabic{listcounter})}{\usecounter{listcounter}\setlength{\rightmargin}{\leftmargin}}
\item
$\Diamond_\Gamma\phi \leftrightarrow \neg\Box_\Gamma\neg\phi$
\item
$\Box(\phi\rightarrow\psi)\rightarrow(\Box\phi\rightarrow\Box\psi)$;
\end{list}
the \textbf{S4} axiom schemes include the \textbf{K} axiom schemes
as well as
\begin{list}
 {(\arabic{listcounter})}{\usecounter{listcounter}\setcounter{listcounter}{2}\setlength{\rightmargin}{\leftmargin}}
\item
$\Box\phi\rightarrow\phi$ and
\item
$\Box\phi\rightarrow\Box\Box\phi$;
\end{list}
and the \textbf{S5} axiom schemes include
\begin{list}
 {(\arabic{listcounter})}{\usecounter{listcounter}\setcounter{listcounter}{4}\setlength{\rightmargin}{\leftmargin}}
\item
$\Diamond\phi\rightarrow\Box\Diamond\phi$
\end{list}
together with $(1)$--$(4)$.

\textit{Remark.} By the duality between $\Box$ and $\Diamond$,
this last axiom is equivalent to \linebreak $(5')$
$\Diamond\Box\phi\rightarrow\Box\phi$, which is the maximality
principle MP if these operators are given the suggested set
theoretic meaning.

We now interpret these schemes in the context of the forcing
relation relativized to a fixed class $\Gamma$ of forcing notions.
If these schemes are interpreted as schemes in \ZFC, then the S4
axiom schemes are intended to hold for all formulas of the
language of \ZFC, with arbitrary parameters allowed.  On the other
hand, the S5 schemes are only consistent in \ZFC\ when the
parameters for the formulas are restricted to specific sets,
depending on the specific class of forcing notions used in
defining the modal operators.  This is due to the equivalence of
the S5 scheme with the modified maximality principle associated
with the forcing class.  This point is important, as the
maximality principles we work with will invariably impose
restrictions on the parameters that can be used in order to avoid
inconsistencies.

We first show that modal operators defined from any such class of
forcing notions obeys the K scheme.  From their definitions, they
clearly obey the duality of the classical modal operators.
Specifically,

\begin{lem}\label{lem_box_diamond_duality}
If $\Gamma$ is any class of forcing notions and $\phi$ is any
formula in the language of \ZFC, then $\Diamond_\Gamma\phi
\leftrightarrow \neg\Box_\Gamma\neg\phi$.
\end{lem}

\begin{proof}
If $\phi$ holds in some $\Gamma$-forcing extension $V[G]$, then,
for its negation to be $\Gamma$-necessary, $\neg\phi$ would have
to hold in $V[G]$ as well, a contradiction.  And, conversely, if
$\neg\phi$ is $\Gamma$-necessary  there can be no extension $V[G]$
in which $\phi$ is true.
\end{proof}

\begin{lem}\label{lem_box_distrib over implic}
Let $\Gamma$ be any class of forcing notions. If $\phi$ and $\psi$
are any formulas in the language of \ZFC, then
$\Box_\Gamma(\phi\rightarrow\psi)\Longrightarrow(\Box_\Gamma\phi\rightarrow\Box_\Gamma\psi)$.
\end{lem}

\begin{proof}
Suppose it is false that
$\Box_\Gamma\phi\rightarrow\Box_\Gamma\psi$, i.e.,
$\Box_\Gamma\phi$ is true but $\Box_\Gamma\psi$ is false.  Then
there is a $\Gamma$-forcing extension $V[G]$ where $\phi$ is true
and $\psi$ is false.  Then
$V[G]\models\neg(\phi\longrightarrow\psi)$.  So
$\Box_\Gamma(\phi\rightarrow\psi)$ must be false.
\end{proof}

The class $\Gamma$ is \textbf{closed under two-step iterations}
if, whenever $\mathbb{P}$ is in $\Gamma$ and $V^\mathbb{P}\models$
``$\mathbb{Q}$ is in $\Gamma$"\footnote{As noted earlier, deciding
whether a forcing notion is in a defined class will always be made
\textit{de dicto} in the forcing extension in which the forcing
notion is to be forced with. $\Gamma$ may be definable in a way
that depends on the model; for example, a forcing notion that is
not $\CCC$ in the ground model may become $\CCC$ in a forcing
extension. Consider the forcing notion $\mathbb{P}$ that collapses
$\omega_1$ to $\omega$, whose conditions are finite partial
injective functions $\omega_1\longrightarrow\omega$.  It has
antichains of size $\omega_1$ but $V^\mathbb{P}\models$
``$\mathbb{P}$ is $\CCC$".}, then $\mathbb{P}*\dot{\mathbb{Q}}$ is
in $\Gamma$.  We will say $\Gamma$ is \textbf{adequate} if
$\Gamma$ is $\Gamma$-necessarily closed under 2-step iterations
and $\Gamma$-necessarily contains trivial forcing, that is, the
forcing notion $\{\emptyset\}$ consisting of only one condition.
Examples of such classes include the class of all forcing notions,
\CCC\ (the class of all forcing notions which satisfy the
countable chain condition), \PROPER\ (the class of all proper
forcing notions) and \CARD\ (the class of all forcing notions that
preserve all cardinals).

\noindent\textit{Remark.} By the definition of adequacy, if
$\Gamma$ is adequate it must then be $\Gamma$-necessarily so
(i.e., this continues to hold in any $\Gamma$-forcing extension).
Again, these axioms are schemes of \ZFC, allowing arbitrary
parameters for the formula in any instance.

\begin{lem}\label{lem adequate ccc proper card}
The forcing classes \CCC, \PROPER, and \CARD\ are all adequate.
\end{lem}

\begin{proof}
They clearly all contain trivial forcing.  The classes \CCC\   and
\PROPER\ are well known to be necessarily closed under two-step
iterations. And, the two-step iteration of cardinal-preserving
forcing notions will itself preserve cardinals, hence it will be
in \CARD.
\end{proof}

Recalling our definitions of $\Gamma$-necessary and
$\Gamma$-forceable, define $\Gamma$ to be a \textbf{K or S4
forcing class} if the operators $\Box_\Gamma$ and
$\Diamond_\Gamma$ satisfy the K or S4 axioms of modal logic,
respectively, and if this is $\Gamma$-necessary (i.e., this
continues to hold in any $\Gamma$-forcing extension).

\begin{thm}\label{thm adequacy}
Every adequate class of forcing notions is an $S4$ forcing class.
\end{thm}

\begin{proof}
Let $\Gamma$ be adequate. We first show that the S4 axioms are
true, then that they are $\Gamma$-necessarily true.

$\Box_\Gamma\phi\rightarrow\phi$: Assume $\Box_\Gamma\phi$, that
$\phi$ is true in every $\Gamma$-forcing extension.  But $\Gamma$
includes trivial forcing, whose extension is the same as the
ground model, where $\phi$ is therefore true.

$\Box_\Gamma\phi\rightarrow\Box_\Gamma\Box_\Gamma\phi$:  Again
assuming $\Box_\Gamma\phi$, we now use the fact that $\Gamma$ is
$\Gamma$-necessarily closed under two-step iterations.  Let
$\mathbb{P}$ be any forcing notion in $\Gamma$.  Since
$\Box_\Gamma\phi$, $V^\mathbb{P}\models\phi$. Now let
$\dot{\mathbb{Q}}$ be any $\mathbb{P}$-name of a forcing notion in
$V^\mathbb{P}$ such that $V^\mathbb{P}\models$
``$\dot{\mathbb{Q}}$ is in $\Gamma$".  Since $\Gamma$ is closed
under two-step iteration, $\mathbb{P}*\dot{\mathbb{Q}}$ is in
$\Gamma$. And since $\Box_\Gamma\phi$,
$V^{\mathbb{P}*\dot{\mathbb{Q}}}\models\phi$. But
$\dot{\mathbb{Q}}$ was arbitrary in $\Gamma^{V^\mathbb{P}}$, so
$V^\mathbb{P}\models\Box_\Gamma\phi$. And since $\mathbb{P}$ was
any forcing notion in $\Gamma$, we have
$\Box_\Gamma\Box_\Gamma\phi$. So $\Gamma$ satisfies the S4 axioms.

Further, since $\Gamma$ is $\Gamma$-necessarily closed under
2-step iterations and trivial forcing, the foregoing argument
shows that, in any $\Gamma$-extension, $\Gamma$ satisfies the S4
axioms there as well. So $\Gamma$ is an S4 forcing class.
\end{proof}

\noindent\textit{Remark}. \thmref{thm adequacy} cannot be
strengthened to ``if and only if" for arbitrary models.  A
counterexample is letting $\Gamma$ be the class of all forcing,
with trivial forcing $\{\emptyset\}$ replaced by some other
forcing notion with the same effect, such as $\{1\}$.  Then
$\Gamma$ is S4 but not adequate. However, we have the following
theorem. We say two forcing notions are
\textbf{forcing-equivalent} if they both produce the same forcing
extensions.  And we define the operation of taking all conditions
below a given condition in a forcing notion to be
\textbf{restricting the forcing notion to a condition}.

\begin{thm}\label{thm_L[A]_adequate_equiv_S4}
If $V=L[A]$ for some set $A$, then every S4 class closed under
forcing-equivalence and the operation of restriction to a
condition is adequate.
\end{thm}

\begin{proof}
Let $\Gamma$ be an S4 class of forcing notions.  Let $\phi$ be the
first-order statement ``$V$ is a $\Gamma$-forcing extension of
$L[A]$." (This statement takes the set $A$ as parameter.)  Clearly
this is true in any $\Gamma$-forcing extension of $L[A]$.  So if
$V=L[A]$, then $\Box_\Gamma\phi$ is true.  So by the S4 axiom
scheme $\Box_\Gamma\phi\longrightarrow\phi$, $\phi$ is true in
$L[A]$.  But if $L[A]$ is a $\Gamma$-forcing extension of itself,
(that is, there is $\mathbb{P}$ in $\Gamma$ which is trivial below
some condition) then trivial forcing is forcing-equivalent to a
forcing notion in $\Gamma$, hence in $\Gamma$ . And if, in $L[A]$,
$\mathbb{P}$ is any forcing notion in $\Gamma$, and $\mathbb{Q}$
is any forcing notion in $\Gamma$ according to $L[A]^\mathbb{P}$,
then, by the S4 axiom scheme
$\Box_\Gamma\phi\longrightarrow\Box_\Gamma\Box_\Gamma\phi$, $\phi$
is true in $L[A]^{\mathbb{P}*\dot{\mathbb{Q}}}$.  That is,
$L[A]^{\mathbb{P}*\dot{\mathbb{Q}}}$ is a $\Gamma$-forcing
extension of $L[A]$.  So $\mathbb{P}*\dot{\mathbb{Q}}$ is forcing
equivalent to a forcing notion in $\Gamma$, hence in $\Gamma$.
This, together with trivial forcing being in $\Gamma$, can also be
shown to be $\Gamma$-necessary. For this, notice that $\Gamma$ is
$\Gamma$-necessarily S4, and that ``$V = L[A]$ for some set $A$"
is $\Gamma$-necessary.  Then repeat the above argument in any
$\Gamma$-extension. This proves that $\Gamma$ is an adequate
forcing class.
\end{proof}

The obstacle for a general proof of the equivalence of S4-ness
with adequacy seems to be the need to refer to the ground model in
the statement $\phi$.  This can be done in the case where $V=L[A]$
by using $A$  as a parameter in $\phi$.

If a class $\Gamma$ of forcing notions is S4 then the principle
$\mathrm{\MP}_\Gamma$ has two equivalent formulations.

\begin{lem}\label{lem_S4_MP}
If $\Gamma$ is an S4 class then the following schemes are
equivalent:
\begin{list}
 {(\arabic{listcounter})}{\usecounter{listcounter}\setlength{\rightmargin}{\leftmargin}}
\item
For any formula $\phi$,
$\Diamond_\Gamma\Box_\Gamma\phi\rightarrow\Box_\Gamma\phi$
\item
For any formula $\phi$,
$\Diamond_\Gamma\Box_\Gamma\phi\rightarrow\phi$.
\end{list}
\end{lem}

\begin{proof}
$(1)$ implies $(2)$:  By S4, $\Box_\Gamma\phi\rightarrow\phi$.
Combining this with
$\Diamond_\Gamma\Box_\Gamma\phi\rightarrow\Box_\Gamma\phi$ gives
$\Diamond_\Gamma\Box_\Gamma\phi\rightarrow\phi$.

$(2)$ implies $(1)$:  Substitute $\Box_\Gamma\phi$ for $\phi$ in
$\Diamond_\Gamma\Box_\Gamma\phi\rightarrow\phi$ to obtain
\begin{equation}\label{eqn boxbox box}
\Diamond_\Gamma\Box_\Gamma\Box_\Gamma\phi\rightarrow\Box_\Gamma\phi
\end{equation}
By the two S4 axioms,
$\Box_\Gamma\Box_\Gamma\phi\leftrightarrow\Box_\Gamma\phi$, so the
left side of this implication changes to give
$\Diamond_\Gamma\Box_\Gamma\phi\rightarrow\Box_\Gamma\phi$.
\end{proof}

Most of the forcing classes we will work with will be S4, so
$\Diamond_\Gamma\Box_\Gamma\phi\rightarrow \phi$ is the form of
$\mathrm{\MP}_\Gamma$ that will be used most frequently.

\begin{lem}\label{lem frc nec is nec}
Let $\Gamma$ be S4.
\begin{list}
 {(\arabic{listcounter})}{\usecounter{listcounter}\setlength{\rightmargin}{\leftmargin}}
\item
$\phi(x)$ is not $\Gamma$-forceably necessary if and only if this
fact is $\Gamma$-necessary.
\item
If members of $\Gamma$ are necessarily in $\Gamma$, then $\phi(x)$
is $\Gamma$-forceably necessary if and only if this fact is
$\Gamma$-necessary.
\end{list}
\end{lem}

\begin{proof}
The right-to-left implications of this lemma are obvious from the
S4 axiom $\Box_\Gamma\phi\longrightarrow\phi$.  So we proceed to
proving each of the left-to-right implications.

 \noindent$(1)$: If there is no
$\Gamma$-forcing extension $V[G]$ where $\phi(x)$ is
$\Gamma$-necessary, then in all forcing extensions it is not
$\Gamma$-forceably necessary.  So it is $\Gamma$-necessarily not
$\Gamma$-forceably necessary.  Formally,
\[\begin{array}{clcr} \neg\Diamond_\Gamma\Box_\Gamma\phi(x) &
\Longrightarrow\neg\Diamond_\Gamma\Diamond_\Gamma\Box_\Gamma\phi(x)
\\ & \Longrightarrow\Box_\Gamma\neg\Diamond_\Gamma\Box_\Gamma\phi(z).
\end{array}\]

This last modal derivation makes use of the property of S4 modal
logic that replacement of a subformula by an equivalent subformula
in a formula leads to an equivalent formula.

\noindent$(2)$: We use the additional hypothesis that members of
$\Gamma$ are necessarily in $\Gamma$. In this case, two forcing
notions $\mathbb{P}$ and $\mathbb{Q}$ of $\Gamma$ can be iterated
in $\Gamma$ since each remains in $\Gamma$ under forcing by the
other. Such iteration is simply product forcing. If $\phi(x)$ is
$\Gamma$-forceably necessary, then there is a $\Gamma$-forcing
extension $V^\mathbb{P}$ where $\phi(x)$ is $\Gamma$-necessary.
Let $\mathbb{Q}$ be in $\Gamma$. Then in $V^\mathbb{Q}$,
$\mathbb{P}$ still forces $\phi(x)$ to be $\Gamma$-necessary (in
$V^{\mathbb{Q}\times\mathbb{P}}$, since
$\mathbb{P}\times\mathbb{Q}=\mathbb{Q}\times\mathbb{P}$, and
$\phi(x)$ is $\Gamma$-necessary in
$V^{\mathbb{P}\times\mathbb{Q}}$, a $\Gamma$-forcing extension of
$V^\mathbb{P}$.
\end{proof}

\subsection{Possibilist quantification and rigid identifiers}
Various quantified modal logics exist, to handle whether or not
the existence of an object can transcend the confines of the world
in which some referring expression is being interpreted.
Quantification is \textbf{possibilist} if the scope of the
quantifier extends beyond the world in which interpretation is
occurring.  The alternative, \textbf{actualist} quantification, is
that which limits existence of an object to the world in which any
reference to it occurs, so quantifiers have their scope limited to
one world. Related to possibilist quantification is the concept of
a \textbf{rigid} designator, a name or referring expression in the
language of a theory whose meaning doesn't change from world to
world. This also occurs only when existence of an object is
allowed to extend to other possible worlds.  Since we are working
in \ZFC, we will appropriate this terminology to use when talking
about forcing extensions simply as models of set theory, not as
parts of a Kripke model, just as we have appropriated the modal
symbols for necessity and possibility without actually talking
about possible worlds or accessibility.  Rigidity can be
formalized using the technique of predicate abstraction.  We will
not do so here.

Whenever a model of \ZFC\ has another model of \ZFC\ as an end
extension (no new elements are added to any set), it is possible
to refer to elements in the first model as though they were in the
second. One opts for rigid designation in this discourse, in the
sense that the models are different but they have corresponding
objects that are structurally the same. This occurs in the set
theory literature, for example, when one writes of making a
supercompact cardinal, $\kappa$, indestructible by Laver's forcing
preparation. The name $\kappa$ is allowed to refer to a cardinal
which exists in two different models of \ZFC, but because the
ground model embeds into the extension, there is no harm in
regarding its elements as belonging to the extension. With forcing
extensions, new sets can be added that are not in the ground
model, but old ones cannot be removed. This property is known as
\textbf{montonicity}. In this situation, using possibilist
quantification, the Converse Barcan formula of quantified modal
logic is valid, which is expressed as either of the dual schemes

\begin{list}
 {(\arabic{listcounter})}{\usecounter{listcounter}\setlength{\rightmargin}{\leftmargin}}
\item
$\Box(\forall x)\phi(x)\rightarrow(\forall x)\Box\phi(x)$ or
\item
$(\exists x)\Diamond\phi(x)\rightarrow\Diamond(\exists x)\phi(x).$
\end{list}

We can easily verify the validity of the first scheme when
translated into pure \ZFC.  It simply says that if it is necessary
(true in all forcing extensions) that, for all $x$, $\phi(x)$ is
true, then, for all $x$, $\phi(x)$ is true in all forcing
extensions.  And this is clear: the ground model contains no
element that isn't in all forcing extensions, so the truth of
$\phi(x)$, for all x in the ground model, must hold in all forcing
extensions.

The concept which is converse to monotonicity is
\textbf{anti-monotonicity}.  It occurs in $\Gamma$-forcing
extensions when sets not in the ground model are never created.
The Barcan formula applies to such forcing classes:

\begin{list}
 {(\arabic{listcounter})}{\usecounter{listcounter}\setlength{\rightmargin}{\leftmargin}}
\item
$(\forall x)\Box\phi(x)\rightarrow\Box(\forall x)\phi(x)$ or
\item
$\Diamond(\exists x)\phi(x)\rightarrow(\exists x)\Diamond\phi(x).$
\end{list}

This will generally only occur in a class of forcing notions that
only contains trivial forcing. However, bounded quantification
over sets will produce an effect analogous to anti-monotonicity,
due to absoluteness of the set membership relation.

\begin{lem}[Restricted Barcan Formula]\label{lem_barcan}
Let $\Gamma$ be a class of forcing notions, and let $A$ be a set.
If for all $x$ in $A$, $\phi(x)$ is $\Gamma$-necessary, then it is
$\Gamma$-necessary that for all $x$ in $A$, $\phi(x)$.  Formally,
$(\forall x \in A)
\Box_\Gamma\phi(x)\longrightarrow\Box_\Gamma(\forall x\in
A)\phi(x)$.
\end{lem}

\begin{proof}
Suppose that for all $x$ in $A$ $\phi(x)$ is $\Gamma$-necessary.
We want to show that $\Box_\Gamma(\forall x\in A)\phi(x)$, which
can be rewritten as $\Box_\Gamma(\forall x) (x \in A
\longrightarrow\phi(x))$.  So let $\mathbb{P}$ be in $\Gamma$,
with $G$ a $V$-generic filter over $\mathbb{P}$, and let $x$ be in
$V[G]$, with $V[G]\models x\in A$. Then $V\models x \in A$ by
absoluteness. Then $\phi(x)$ is $\Gamma$-necessary by hypothesis.
So $V[G]\models\phi(x)$.  Thus $V[G]\models(x \in
A\longrightarrow\phi(x))$.  Since $\mathbb{P}$ was arbitrary in
$\Gamma$, and $x$ was arbitrary in $V[G]$, this gives
$\Box_\Gamma(\forall x)( x \in A\longrightarrow\phi(x))$.
\end{proof}

\subsection{Modal logic of forcing Kripke models} This section is
added as a bridge to the traditional interpretation of the modal
operators, since there is an intuitive interpretation of forcing
that allows one to think of forcing extensions as separate models
of \ZFC\ completely outside the universe of the ground model.
Although these are not the formal semantics we are using, they can
be intuitively thought of as a traditional \textbf{Kripke model},
a nonempty set of objects called \textbf{possible worlds} together
with a binary relation on those objects called an
\textbf{accessibility relation}.  All the possible worlds in a
Kripke model interpret the same fixed first-order language, which
is augmented with the symbols $\Box$ and $\Diamond$.  For any
formula $\phi$, $\Box\phi$ asserts that $\phi$ is
\textbf{necessary}, or true in every accessible world, while
$\Diamond\phi$ asserts that $\phi$ is \textbf{possible}, or true
in some accessible world.

Any Kripke model requires a bare-bones modal logic (called
\textbf{K}) that includes at least all axioms of classical
first-order logic (including the rule of inference \textit{modus
ponens}), together with the axiom scheme, for all formulas $\phi$
and $\psi$,
\[\Box(\phi\rightarrow\psi)\rightarrow(\Box\phi\rightarrow\Box\psi)\]
and the \textit{necessitation} rule of inference
\[\frac{\phi}{\Box\phi}\]

(The necessitation rule means that if the formula $\phi$ occurs in
a proof tableau or proof sequence, then the formula $\Box\phi$ is
a valid deduction, since $\phi$ has been proved according to
axioms valid in all worlds, hence $\phi$ is true in all accessible
worlds.  In short, if $\phi$ is provable, then $\Box\phi$.  It
does \textit{not} say $\phi\rightarrow\Box\phi$.)  Additional
axioms and rules of inference depend on the type of Kripke model
involved.  If the accessibility relation of a Kripke model obeys
reflexivity and transitivity, it can be shown to obey the
\textbf{S4} axiom schemes mentioned earlier.  If that relation is
also symmetric, it obeys the \textbf{S5} axiom scheme also given
earlier.  In fact the accessibility relation we will consider,
that of being a $\Gamma$-forcing extension, is not symmetric
unless $\Gamma$ consists only of trivial forcing.

We now interpret these ideas in the context of the forcing
relation.  Given any model $M$ of set theory, there are many
models of set theory derivable from $M$ through the mechanism of
forcing. We call the Kripke model that is the totality of these
models, together with the accessibility relation ``$M'$ is a
forcing extension of $M$" to be a \textbf{forcing Kripke model}.
If $\Gamma$ is a class of forcing notions we can relativize this
idea to a $\Gamma$-forcing Kripke model as well.

Again, in this work, the interpretation of modal logic symbols and
terminology will be within \ZFC, hence within the single universe
$V$ that models it, not in this forcing Kripke model of possible
worlds. But the correspondence of ideas will be utilized freely.
For example, without adding it to \ZFC\ as a new rule of
inference, the necessitation rule is used in the form that says
that if $\phi$ is provable it is then true in all models of \ZFC.

\section{Absolutely definable parameter sets} Consider
parameters to formulas in instances of $\MP_\Gamma$.  Notice that
if a set parameter has a $\Gamma$-absolute definition it needn't
really be a parameter at all, since it can be eliminated by
replacing its occurrence by its definition. So, for any set $X$
whose elements are definable, the principle $\MP_\Gamma$ implies
$\MP_\Gamma(X)$ for such $X$.  What is less clear, if $X$ itself
is definable but its elements are not required to be so, is
whether $\MP_\Gamma$ still implies $\MP_\Gamma(X)$. Situations of
this type raise questions such as the following.\vspace {12 pt}

\noindent\textbf{Question 1.} If $\Gamma\subseteq \CARD$, does
$\MP_\Gamma$ imply $\MP_\Gamma(\omega_1)$? \vspace {12 pt}

\noindent\textbf{Question 2.} If $\Gamma$ adds no reals, does
$\MP_\Gamma$ imply $\MP_\Gamma(\mathbb{R})$?\vspace {12 pt}

Note that these questions are about parameters that are contained
in $\omega_1$ or $\mathbb{R}$, not the sets themselves.  Thus the
principle  $\MP_\Gamma(\omega_1)$ asserts that
$\Diamond_\Gamma\Box_\Gamma\phi(\alpha)$ implies
$\Box_\Gamma\phi(\alpha)$ for any formula $\phi$ with one free
variable, whenever $\alpha$ is in $\omega_1$.

Formally, we will call a \textbf{definition} of a set $A$ a
formula $\psi(x)$ such that for all $x$, $\psi(x)$ if and only if
$x = A$.  Let $\Gamma$ be a class of forcing notions. A formula
$\phi(x)$ is \textbf{$\Gamma$-absolute} if for any $\mathbb{P}$ in
$\Gamma$, and for any $x$ in $V$, $V^\mathbb{P}\models\phi(x)$ if
and only if $\phi(x)$. A set defined by such a formula $\phi(x)$
in all $\Gamma$-forcing extensions as well as in $V$ is said to be
\textbf{$\Gamma$-absolutely definable}.

Throughout this chapter $\Gamma$ will be an S4 forcing class. Let
$S$ be a $\Gamma$-absolutely definable set (so the definition of
$S$, when interpreted in any $\Gamma$-forcing extension of $V$,
results in the same set---no new elements are added).  We give
some examples to illustrate the idea:
\begin{list}
 {(\arabic{listcounter})}{\usecounter{listcounter}\setlength{\rightmargin}{\leftmargin}}
\item
$\omega_1^L$ is absolutely definable over all forcing notions.
\item
Let $\Gamma_{\omega_1}$ be the class of forcing notions that
preserve $\omega_1$.  Then the cardinal $\omega_1$ is
$\Gamma_{\omega_1}$-absolutely definable.
\item
The set of reals, $\mathbb{R}$, is not absolutely definable over
the class of all forcing notions.
\end{list}

Let $\MP_{\Gamma}(S)$ be the form of the principle $\MP_{\Gamma}$
in which statements take parameters from the set $S$.  The fact
proven in this chapter is that $\MP_{\Gamma}\Longrightarrow
\MP_{\Gamma}(S)$. The strategy in this proof is, for every formula
$\phi(\alpha)$ which takes a parameter $\alpha$ in $S$, to
consider the statement $\forall \alpha \in S (\Diamond \Box
\phi(\alpha) \longrightarrow \phi(\alpha))$.  Since this formula
has no free variables (after replacing $S$ with its definition),
it is subject to the principle $\MP_{\Gamma}$ if it holds. If we
can show that this statement is $\Gamma$-forceably necessary,
then, by $\MP_{\Gamma}$, it will be true.  But taking this over
all formulas $\phi$ we have the scheme $\MP_{\Gamma}(S)$.

We use a lemma whose proof is simply an exercise in modal logic,
using the S4 axioms which characterize the Kripke model of any
class of forcing extensions in which accessibility satisfies
reflexivity and transitivity.  Since $\Gamma$, the forcing notion
class under discussion, satisfies these, we will let the modal
operators $\Box$ and $\Diamond$ be shorthand for $\Box_{\Gamma}$
and $\Diamond_{\Gamma}$ respectively throughout the rest of this
chapter.

\begin{lem}
Let $\Gamma$ be an S4 forcing class.  If $\phi(z)$ is a formula in
the language of \ZFC\ with parameter $z$ then the statement
``$\Diamond\Box\phi(z)$ implies $\phi(z)$" is forceably necessary.
\end{lem}

\begin{proof}

Case I: $\phi(z)$ is forceably necessary ($\Diamond\Box\phi(z)$):
Then there is a forcing extension $V[G]$ after which all forcing
extensions satisfy $\phi(z)$.  Therefore they satisfy that
$\phi(z)$ is forceably necessary implies $\phi(z)$.

(Formally, in S4,
\[\begin{array}{clcr}\Diamond\Box\phi(z) &
\Longrightarrow\Diamond\Box(\Diamond\Box\phi(z) \rightarrow
\phi(z)).)
\end{array}\]

(This derivation is an instance of the S4 tautology
\[\Diamond\Box A\longrightarrow\Diamond\Box (B\longrightarrow A).)\]

Case II: $\phi(z)$ is not forceably necessary
($\neg\Diamond\Box\phi(z)$): Then in all forcing extensions
$\phi(z)$ cannot be forceably necessary, whence it is forceably
necessary that $\phi(z)$ is not forceably necessary.  But in any
world where the forceable necessity of $\phi(z)$ is false, the
forceable necessity of $\phi(z)$ implies everything, including
$\phi(z)$. (Formally, in S4,
\[\begin{array}{clcr} \neg\Diamond\Box\phi(z) &
\Longrightarrow\Box\neg\Box\phi(z) \\ &
\Longrightarrow\Box\Box\neg\Box\phi(z) \\ &
\Longrightarrow\Box(\neg\Diamond\Box\phi(z)) \\ &
\Longrightarrow\Diamond\Box(\neg\Diamond\Box\phi(z)) \\ &
\Longrightarrow\Diamond\Box(\Diamond\Box\phi(z)\rightarrow
\phi(z)).)
\end{array}\]

\end{proof}

We will say that a class $\Gamma$ is closed under iterations of
order type $\gamma$ with \textbf{appropriate support} if there is
an ideal on $\gamma$ (or other definable choice of support at
limit stages) which when used as support for any such iteration
produces a forcing notion which is still in $\Gamma$.  Moreover,
we require that if $\mathbb{P}_\Gamma$ is a $\gamma$-iteration of
forcing in $\Gamma$ and $\alpha<\gamma$, then $\mathbb{P}_\gamma =
\mathbb{P}_\alpha*\mathbb{P}_{\{\alpha,\gamma\}}$ such that
$\mathbb{P}_{\{\alpha,\gamma\}}$ is an iteration of forcing
notions of $\Gamma$ with support in the ideal, as defined in
$V^{\mathbb{P}_\alpha}$.   An example of such support would be
where either direct or inverse limits are always taken; this is
Lemma 21.8 in \cite{jec:sth}.

\begin{lem}\label{lem_sigma_Gamma_forc_nec}
Let $\Gamma$ be an S4 forcing class, $\Gamma$-necessarily closed
under iterations of arbitrary length with appropriate support. Let
$S$ be any set.  Let $\phi(x)$ be a formula in the language of set
theory with one free variable. Let $\sigma$ be the formula
``$\forall s\in S (\Diamond\Box\phi(s)$ implies $\phi(s))$". Then
$\sigma$ is $\Gamma$-forceably necessary.
\end{lem}

\begin{proof}
Let ordinal $\gamma$ be the order type of a fixed wellordering of
$S=\langle s_\alpha\mid \alpha<\gamma\rangle$.  Let $\sigma(s)$
stand for the formula ``$\Diamond\Box\phi(s)$ implies $\phi(s)$",
where $s$ is a parameter of $\phi$.  By the previous lemma,
$\sigma(s_0)$ is forceably necessary, and, more to the point, this
fact is a theorem of \ZFC, hence, by necessitation, absolute to
all $\Gamma$-forcing extensions.  We begin a finite support
$\gamma$-iteration of $\Gamma$-forcing notions by taking the
extension $V^{\mathbb{P}_1}$ in which $\sigma(s_0)$ is necessary.
Since $\sigma(s_\alpha)$ is still forceably necessary in
$V^{\mathbb{P}_1}$ for any $\alpha$ in $\gamma$, we can repeat
this construction for another stage, but this time using $1$ for
$\alpha$. This gives a model, $V^{\mathbb{P}_2}$, in which both
$\sigma(s_0)$ and $\sigma(s_1)$ are necessary, whence
$\sigma(s_0)\wedge\sigma(s_1)$ is true and remains true in any
extension. Continue in this manner for $\gamma$ stages: at stage
$\alpha$, let $\mathbb{Q}_\alpha$ be that $\Gamma$ forcing notion
which forces the $\Gamma$-necessity of $\sigma(s_\alpha)$.  Taking
a $\gamma$-iteration with appropriate support $\mathbb{P} =
\mathbb{P}_{\gamma}$, we obtain a $\Gamma$-forcing extension
$V^{\mathbb{P}_{\gamma}}$. Let $G$ be $\mathbb{P}$-generic over
$V$.  We claim first that, for every $\alpha$ in $\gamma$, $V[G]
\models \sigma(s_\alpha)$. To see this, choose any $\alpha$ in
$\gamma$.  At stage $\alpha$, we can factor $\mathbb{P} =
\mathbb{P}_\alpha*\dot{\mathbb{P}}_{TAIL}$, and take $G_\alpha$ as
$V$-generic over $\mathbb{P}_\alpha$.  But at stage $\alpha$,
$V[G_\alpha]\models$ ``$\mathbb{Q}_\alpha$ forces that
$\sigma(s_\alpha)$ is $\Gamma$-necessary".  So refactoring,
$V[G_{\alpha+1}]\models$ ``$\sigma(s_\alpha)$ is
$\Gamma$-necessary". But $V[G]$ is a $\Gamma$-forcing extension
over $V[G_{\alpha+1}]$ (since $\Gamma$ is $\Gamma$-necessarily
closed), so $V[G]\models$ ``$\sigma(s_\alpha)$ is
$\Gamma$-necessary". Formally, $V[G]\models (\forall \alpha \in
\gamma) \Box \sigma(s_\alpha)$. By \lemref{lem_barcan}, $(\forall
\alpha\in\gamma)\Box\sigma(s_\alpha)\rightarrow\Box(\forall \alpha
\in\gamma)\sigma(s_\alpha)$.  Modus ponens then gives the
statement $\Box(\forall \alpha\in\gamma)\sigma(s_\alpha)$, that
is, $\sigma$ is $\Gamma$-forceably necessary.
\end{proof}

\begin{thm}\label{thm_MP_Gamma_S}
Let $\Gamma$ be an S4 forcing class, $\Gamma$-necessarily closed
under iterations of arbitrary order type with appropriate support.
Let $S$ be a set which is $\Gamma$-absolutely definable. Then
$\MP_{\Gamma}$ if and only if $\MP_{\Gamma}(S)$.
\end{thm}
\begin{proof}
For any formula $\phi$ in the language of \ZFC, and for any $s$ in
$S$, the formula
\[\forall s \in S
(\Diamond\Box\phi(s)\longrightarrow \phi(s))\] is
$\Gamma$-forceably necessary by the last lemma.  Let $\psi$ be the
$\Gamma$-absolute definition of $S$.  Then the above formula is
equivalent to the formula \[\forall s (\psi(s)\longrightarrow
(\Diamond\Box\phi(s)\longrightarrow \phi(s)))\] which is free of
parameters.  So by $\MP_{\Gamma}$ it is true.
\end{proof}

Another way to express this result is to say that $\MP_\Gamma$
alone allows free use of parameters in $S$.  And that is the point
of this section---that as long as the parameter set is absolutely
definable, the actual parameters need not be definable themselves.
There is a way to generalize still further. If $\MP_\Gamma(X)$ is
assumed a priori, the argument above allows the adjoining to $X$
of a new set of parameters $P$ which is $\Gamma$-absolutely
definable.

\begin{cor}\label{cor MP_Gamma_S_X}
Let $\Gamma$ be a class of forcing notions obeying the S4 axioms
of modal logic and closed under arbitrary iterations with
appropriate support, and let $S$ be a set $\Gamma$-absolutely
definable from parameters in $X$.  Then $\MP_\Gamma(X)$ if and
only if $\MP_\Gamma(S\cup X)$.
\end{cor}

\begin{proof}
In the presence of $\MP_\Gamma(X)$, allow the formula $\phi$ to
take parameters from $X\cup S$.  Suppose $\phi(\alpha)$ is
forceably necessary.  If $\alpha$ is in $X$, then $\phi(\alpha)$
is true by virtue of $\MP_\Gamma(X)$. Otherwise, $\alpha$ is in
$S$, which is a set whose definition is $\Gamma$-absolutely
definable from parameters in $X$.  So the sentence
$\forall\alpha\in S\cup X (\Diamond\Box\phi(\alpha)\longrightarrow
\phi(\alpha))$ is forceably necessary by the above argument.  But
this sentence can be modified to take only parameters from $X$ by
using the definition of $S$, in which form it will be true by
applying $\MP_\Gamma(X)$. The collection of all such sentences
then holds, and this collection constitutes the scheme
$\MP_\Gamma(S\cup X))$.
\end{proof}

Finally we consider the union of all $\Gamma$-absolutely definable
parameter sets which consists of the least initial segment which
includes all $\Gamma$-absolutely definable sets.  We can call this
the \textbf{$\Gamma$-absolutely definable cut} of the universe,
$V_{\Gamma}$.  This is analogous to the definable cut of the
universe, the least initial segment which includes all definable
sets. (In contrast with the latter class, however, $V_{\Gamma}$
might not be a model of \ZFC.)

\begin{cor}
$\MP_\Gamma$ if and only if $\MP_\Gamma(V_{\Gamma})$.
\end{cor}

\begin{proof}
Assume $\MP_\Gamma$.  Any instance of the scheme
$\MP_\Gamma(V_{\Gamma})$ will be an instance of the scheme
$\MP_\Gamma(S)$ for some $\Gamma$-absolutely definable parameter
set $S$, which is true by \thmref{thm_MP_Gamma_S}
\end{proof}

Again, one can say that $\MP_\Gamma$ alone allows free use of
parameters in $V_{\Gamma}$.  Returning to the questions asked at
the beginning of this section, we can now provide answers.
Responding to question 1,  let \CARD\ be the class of
cardinal-preserving notions of forcing.

\begin{cor}
If the ordinal $\alpha$ is $\CARD$-absolutely definable, then
$\MPmod{\CARD}$ implies $\MPmod{\CARD}(\omega_\alpha)$.  In
particular one can freely use $\omega_1$, $\omega_3$,
$\omega_{\omega_{19}}$, and beyond, as parameter set $X$ in the
$\MPmod{\CARD}$ scheme.
\end{cor}

\begin{proof}
By \thmref{thm_MP_Gamma_S}.
\end{proof}

Responding to question 2:

\begin{cor}
Let $\Gamma$ be the class of all forcing notions that add no new
reals.  If $\MP_\Gamma$ holds, then $\MP_\Gamma(\mathbb{R})$
holds.
\end{cor}

\begin{proof}
Clearly $\Gamma$ is adequate, hence S4.  Also, $\mathbb{R}$ is
$\Gamma$-absolutely definable.  So apply \thmref{thm_MP_Gamma_S}.
\end{proof}

\section{Elementary submodels of $V$} Several theorems in this
work require construction of a forcing iteration where each
successor stage forces a particular instance of the maximality
principle whose model is being sought; the principle is then true
in the iterated extension because each instance has been handled
at some stage.  But in order to define a forcing notion that
forces even one instance of a maximality principle requires
expressing that a sentence is forceably necessary.  The relation
$p\Vdash\phi$ is defined in $V$ with a separate definition for
each $\phi$, constructed by induction on the complexity of $\phi$.
So the general relation $p\Vdash\phi$ cannot be expressed as a
formula with arguments for $p$ and $\phi$.  To express this
relation requires the scheme consisting of each particular case.
This changes for a set model $M$ within $V$.  In this case one can
express truth (and forcing) of $\phi$ in $M$ as a predicate
interpreted in $V$.

Since the forcing predicate is only definable over set models of
\ZFC, we employ the strategy of \cite{ham:max} of using an initial
segment of the universe as an elementary submodel of it.  One
should be warned that the models of \ZFC\ in which this occurs
might possibly be nonstandard (nonwellfounded) models\footnote{An
alternative approach, not taken here, would be to avoid these
nonstandard models by constructing a model for a maximality
principle that applies to any finite fragment of \ZFC\ and then
finally applying the compactness theorem to find a model of the
full maximality principle.  Our approach will be to apply the
compactness theorem earlier on, and allow iterations to construct
models of the full maximality principles.}. In particular, we
generalize Lemma 2.5 of \cite{ham:max}. We first add a constant
symbol $\delta$ to the language of \ZFC, intended to stand for
some ordinal.  The following arguments will take place in this
expanded language, together with an expanded model of \ZFC\ to
interpret this constant.

Let ``$V_\delta \prec V$" stand for the scheme that asserts, of
any formula $\phi$ with a parameter $x$, that

\begin{center}
 for every $x \in V_\delta$, $V_\delta \models
\phi[x]$ if and only if $\phi(x)$.
\end{center}

\begin{lem}\label{lem_Vdelta_submod_V}
Let $T$ be any theory containing \ZFC\ as a subtheory. Then

$Con(T)$ if and only if $Con(T+V_\delta\prec V)$
\end{lem}

\begin{proof}

The implication to the left is trivial.  To obtain the implication
to the right, let $M$ be a model for $T$.  I will show that, with
a suitable interpretation of $\delta$, $M$ is a model for any
finite collection of formulae in $T+V_\delta\prec V$, which is
therefore consistent. The conclusion then follows by the
compactness theorem.

Let $\Psi^*$ be any finite collection of instances of the scheme
$V_\delta\prec V$. Let $\Psi$ be the collection of formulas
$\psi(x)$, in the language of $T$, for which there is an instance
in $\Psi^*$ of the form $``\forall x \in V_\delta \ V_\delta
\models \psi[x] \ \textup{if and only if}\ \psi(x)"$. We can write
$\Psi^* = \{``\forall x \in V_\delta \ V_\delta \models \psi[x] \
\textup{if and only if}\ \psi(x)"|\psi \in \Psi\}$.

$\Psi$ is finite, so by L\'evy reflection there is an initial
segment, $M_\gamma$, of $M$ such that for all $\psi \in \Psi$, and
for all $x \in M_\gamma$, if $M_\gamma \models \psi[x]$ then
$\psi(x)$. So, interpreting $\delta$ as $\gamma$, $M \models
\Psi^*$ (so $\Psi^*$ is consistent).  As $\Psi^*$ was an arbitrary
finite collection of formulas from $V_\delta\prec V$, and since
$M$ also satisfies any finite fragment of $T$ (being a model of
$T$), the entire theory $T+V_\delta\prec V$ is therefore
consistent.
\end{proof}

The proof of \lemref{lem_Vdelta_submod_V} relies on a compactness
argument to establish a model $M'$ for $\ZFC+V_\delta\prec V$
given the existence of a model $M$ for \ZFC\ alone.   In fact, a
sharper observation can be made which relates these models,
namely, it can be arranged that $M\prec M'$.  This can be shown by
``reproving" the compactness theorem in building $M'$ as an
ultraproduct.  An alternate proof is also given which follows a
Henkin style proof.

We first introduce a notation for an idea that will recur through
the chapter.  Let $T$ be any theory containing \ZFC\ as a
subtheory. Let $I=\{\Phi \subseteq T|\Phi$ is finite $\}$ be the
set of finite collections of formulas in the theory $T$.  For each
$\Phi$ in $I$, let $V_\delta\prec_\Phi V$ denote the collection of
statements $\{``\forall x \in V_\delta (V_\delta \models \phi[x]$
if and only if $\phi(x))"|\phi \in \Phi\}$, that is, those
instances of $V_\delta\prec V$ only mentioning those formulas
$\phi$ in $\Phi$. Notice that any finite subcollection of the
scheme $V_\delta\prec V$ can be so represented.

\begin{lem}\label{lem_Vdelta_submod_V_in_elem_ext}
Let $T$ be any theory containing \ZFC\ as a subtheory.  Suppose
$M\models T$.  Then there is $M'\models T+V_\delta\prec V$ such
that $M\prec M'$.
\end{lem}

\begin{proof}[First proof]
Let $M$ be a model for theory $T$.  By the L\'evy Reflection
Theorem, there is a $\delta = \delta_\Phi$ for which
$V_\delta\prec_\Phi V$.  So the expansion of the model $M$ to
$\langle M,\delta_\Phi\rangle$ is a model of $T+V_\delta\prec_\Phi
V$. Denote this expanded model by $M_\Phi$. (Notice that if $\Phi
\subseteq \Psi$, then $M_\Psi \models V_\delta\prec_\Phi V$.)  Now
construct an ultrafilter on $I$ as follows: For each $\Phi \in I$,
define $d_\Phi = \{\Psi\in I|\Phi \subseteq \Psi\}$, the set of
finite collections of formulas containing $\Phi$ as a
subcollection. Then $D_I = \{d_\Phi|\Phi \in I\}$ is easily seen
to be a filter on $I$.  And by Zorn's Lemma, there is an
ultrafilter $U_I \supseteq D_I$. We can now take the ultraproduct
$\langle M',\delta \rangle=\prod M_\Phi/U_I$. Then $\langle M',
\delta\rangle$ is a model of $V_\delta\prec_\Phi V$, for any
$\Phi$ in $I$ (just apply \L os's Theorem to the set $\{\Psi \in
I|M_\Psi \models V_\delta\prec_\Phi V\} \supseteq d_\Phi \in
U_I$).

Since $\langle M', \delta\rangle$ is a model of any finite
subcollection of $V_\delta\prec V$ (where $\delta$ is the element
of the ultraproduct which represents the equivalence class of the
mapping $I$ to $M$ via $\Phi \mapsto \delta_\Phi$), it must
satisfy the entire scheme $V_\delta\prec V$.  Finally, the reduct
$M'$ of $\langle M',\delta\rangle$ is simply the ultrapower of the
model $M$ over the ultrafilter $U_I$, so $M' \models T$ and $M
\prec M'$.
\end{proof}

\begin{proof}[Second proof]
Let $M$ be a model for theory $T$.  Let $T'$ be the elementary
diagram of M together with the scheme $V_\delta \prec V$.  By a
L\'evy Reflection argument, $T'$ is finitely consistent: any
finite set of statements of $T'$ is consistent.  So by the
compactness theorem, there is a model $M'$ of $T'$.  But any such
model is an elementary extension of $M$, since each element of $M$
interprets its own constant in the elementary diagram which
subsequently has an interpretation in $M'$.
\end{proof}

The next lemma says that $V_\delta\prec V$ persists over forcing
extensions when forcing with forcing notions contained in
$V_\delta$.  We will refer to such forcing extensions as ``small".
(The condition $\mathbb{P}\in V_\delta$ precludes, say, collapsing
$\delta$ to $\omega$, which would destroy the scheme
$V_\delta\prec V$.)  In the following results, $V_\delta[G]\prec
V[G]$ will mean the obvious thing, namely that $V[G]\models
V_\delta\prec V$, where $V_\delta$ is interpreted as $V_\delta[G]$
in $V[G]$.  The expression $V_\delta[G]$ is unambiguous, since in
our usage, $G$ is always generic over small forcing, in which case
$(V_\delta)[G]=(V[G])_\delta$.

\begin{lem}\label{lem_Vdelta_submod_V_forcing_ext}
Let $V_\delta\prec V$, let $\mathbb{P}\in V_\delta$ be a notion of
forcing and let $G$ be $V$-generic over $\mathbb{P}$. Then
$V_\delta[G]\prec V[G]$.
\end{lem}

\begin{proof}
Suppose $x\in V_\delta[G]$ such that $V_\delta[G]\models\phi(x)$.
Then there is a condition $p \in G\subseteq \mathbb{P}$ such that
$V_\delta\models p\Vdash\phi(\dot{x})$, so by elementarity
$V\models p\Vdash \phi(\dot{x})$, hence $V[G]\models \phi(x)$.
\end{proof}

This is the way an initial segment $V_\delta$, an elementary
submodel of the universe, is used in a forcing iteration to obtain
a model in which a desired maximality principle holds. Once a
forcing notion has been found at each stage to force the necessity
of a particular forceably necessary formula in the elementary
submodel, a generic $G$ can then be taken to produce an actual
extension which is also an elementary submodel in which the next
iteration is definable.

Another direction in which \lemref{lem_Vdelta_submod_V} can be
generalized is to show there is an unbounded class of cardinals
$\delta$ for which $V_\delta\prec V$. This uses the following
strong form of the reflection theorem.

\begin{lem}[L\'evy]\label{lem_Levy_club_class}
For any finite list of formulas $\Phi$ in the language of \ZFC\
the following is a theorem of \ZFC:

There is a club class of cardinals such that for each $\delta$ in
$C$, $V_\delta\prec_\Phi V$.
\end{lem}
\begin{proof}
Without loss of generality, suppose the list $\Phi = \{
\phi_1\ldots\phi_n\}$ is closed under subformulas.  Define a class
function $f:ORD\rightarrow ORD$ as follows: For $\alpha$ in $ORD$,
let $f(\alpha)$ be the least ordinal $\gamma$ such that, for all
$x$ in $V_\alpha$, and for all $i=1\ldots n$ there exists $y$ such
that $\phi_i(x,y)$ and $y$ is in $V_\gamma$.  Let $D = \{\beta \in
ORD |f``\beta \subseteq \beta \} =$ the closure points of $f$. It
is easy to see that $D$ is a club class.  For any $\delta$ in $D$,
absoluteness for $\phi_1\ldots\phi_n$ between $V_\delta$ and $V$
can be proven by induction on the complexity of each $\phi_i$:
Absoluteness will be preserved under boolean connectives, and the
same is true under existential quantification (the Tarski-Vaught
criterion is satisfied since $\delta$ is a closure point of the
function $f$).  Since the cardinals also form a club class, the
intersection of them with $D$ will be the desired club $C$.
\end{proof}

This leads to the next theorem, another variation of
\lemref{lem_Vdelta_submod_V}.  As in
\lemref{lem_Vdelta_submod_V_in_elem_ext}, the ``ground" model is
expanded to interpret a new predicate symbol in the language,
which in this case is a name for a club class.

\begin{thm}\label{thm_Vdelta_submod_V_club_class}
Let $T_0$ be a theory containing \ZFC.  Then the following are
equivalent:
\begin{list}
 {(\arabic{listcounter})}{\usecounter{listcounter}\setlength{\rightmargin}{\leftmargin}}
\item
$Con(T_0)$
\item
$Con(T_0+T)$
\end{list}
where $T$ is the theory in the language $\{\in, C\}$ which
contains all instances of the Replacement and Comprehension axiom
schemes augmented with the relation symbol C added to the
language, and which asserts
\begin{list}
 {(\roman{listcounter})}{\usecounter{listcounter}\setlength{\rightmargin}{\leftmargin}}
\item
$C = \langle\delta_\alpha|\alpha\in ORD\rangle$ is an unbounded
class of cardinals.
\item
For any formula $\phi$, the following is true: for all $\delta$ in
$C$, $V_\delta \prec V$.
\item
for all $\alpha$ in $ORD$, $\delta_\alpha <
\cof(\delta_{\alpha+1})$.
\end{list}
\end{thm}

\begin{proof}
(2) implies (1): Trivial.

(1) implies (2):  We first give a model where theory $T$ only
includes assertions (i) and (ii).  Let $M$ be a model of theory
$T_0$. It will suffice to show that every finite subtheory of
$T_0+T$ is consistent.  So let $F$ be such a finite subtheory.
That is, $F\subseteq \ZFC+\{\sigma_1,\ldots \sigma_n\}\cup\{$``$C$
is club"$\}$, where $\sigma_i = ``\forall \delta\in C \forall x
\in V_\delta (V_\delta\models \phi_i[x]$ if and only if
$\phi_i(x))$". We show that $M$ is a model for $F$: Fix
$\{\phi_1,\ldots,\phi_n\}$. By \lemref{lem_Levy_club_class}, there
is a definable club class of cardinals $C$ such that for all $i$
in $1,\ldots, n$, $M \models$ ``$C$ is a club
class"+``$\forall\delta\in C$, $V_\delta \models \phi_i[x]$ if and
only if $\phi_i(x)$". But this is $\sigma_i$, so $M\models
\sigma_i$. And since $M$ models $T_0$, it is a model for $F$.  And
by compactness, $M$ all of $T_0$ and $T$, where $T$ includes
assertions (i) and (ii).  But from $M$ we can find a model which
has a club $C$ which satisfies assertion (iii) as well.  Simply
define a new class club to be a continuous subsequence of the
original one by inductively defining, for each ordinal $\alpha$,
$\delta_{\alpha+1} = $ the ${\delta_\alpha^+}^{th}$ element that
follows $\delta_\alpha$ in $C$.  Take suprema (which are, in fact,
unions) at limits. Interpreting the symbol $C$ to be this thinned
club ensures that $\delta_\alpha <\cof(\delta_{\alpha+1})$.
\end{proof}

One should stress that in all applications of
\thmref{thm_Vdelta_submod_V_club_class} one is using the language
of \ZFC\ expanded to have the relation symbol $C$, and that the
theory includes theory $T$ as described in the theorem as well as
\ZFC\ in this expanded language.

\lemref{lem_Vdelta_submod_V_in_elem_ext} makes it conceptually
easier to include in the theory $T$ statements about some element
$\kappa$ of $V$ referring to it by a name added to the language of
\ZFC.  This is done by expanding the model $M$ to interpret the
name of $\kappa$. Equiconsistency of such statements together with
$V_\delta \prec V$ follows from
\lemref{lem_Vdelta_submod_V_in_elem_ext} since the name for
$\kappa$ is ``rigid"---the model for $V_\delta \prec V$ can be
taken to be an elementary extension, so the same $\kappa$ can be
found in both models. The next lemma, illustrating this, is an
enhancement of \lemref{lem_Vdelta_submod_V} that provides a
condition on $\delta$.

\begin{lem}\label{lem_Vdelta_submod_V_big_cof}
Let $T$ be any theory containing \ZFC\ as a subtheory.  Let
$\kappa$ be any ordinal in a model $M$ of $T$ expanded to include
$\kappa$. Then there is an elementary extension $M'$ of $M$ which
is a model of $T+V_\delta\prec V+\cof(\delta)>\kappa$
\end{lem}

\begin{proof}
We proceed exactly as in \lemref{lem_Vdelta_submod_V}, performing
additional work to address the requirement on $\delta$:

By L\'evy Reflection \lemref{lem_Levy_club_class}, for any fixed
finite $\Phi\subseteq T$, the class $\{\alpha \in
ORD|V_\alpha\prec_\Phi V\}$ is closed and unbounded.  So it has a
${\kappa^+}^{th}$ member.  We interpret $\delta$ as this member,
giving $\cof(\delta)=\kappa^+$ and $V_\alpha\prec_\Phi V$.  The
rest of the proof of  \lemref{lem_Vdelta_submod_V} now gives
$M'\models T+V_\delta\prec V $+$ \cof(\delta) > \kappa$.
\end{proof}

A typical application of this lemma is to ensure that the
cofinality of $\delta$ is greater than $\omega$.

We include one equiconsistency result to be used in the next
chapter that uses the results of this section.

\begin{lem}\label{lem_MP_S4_X_one_step}
Let $X$ be any set, and let $\Gamma$ be an S4 class of forcing
notions, closed under iterations of length $\kappa=|X|$ with
appropriate support. Suppose $V_\delta \prec V$ and $\cof(\delta)
> \kappa$. Then there is a forcing notion $\mathbb{P}$ in $\Gamma$
such that $V^{\mathbb{P}} \models \MP_{\Gamma}(X^V)$ and
$\mathbb{P} \in V_\delta$.
\end{lem}

\noindent\textit{Remark}.  The notation ``$X^V$" for the parameter
set is to emphasize the \textit{de re} interpretation of the
symbol $X$--- that $X$ represents the same set in $V^\mathbb{P}$
as when interpreted in $V$.  This is a different situation from
\lemref{lem_MP_ccc_R_one_step}, where the parameter set
$\mathbb{R}^V$ is different from $\mathbb{R}$ in $V^\mathbb{P}$.
This is because $\mathbb{R}$ is a definition of a set, and is
interpreted \textit{de dicto}, hence differently in $V$ and
$V^\mathbb{P}$.

\begin{proof}
Let $\kappa = |X|$. Let $\pi:\omega\times
\kappa\longrightarrow\kappa$ be a bijective pairing function.
Enumerate all formulas with one parameter in the language of set
theory as $\langle \phi_n|n \in \omega \rangle$ and all elements
$x$ of $X$ as $\langle x_\mu|\mu \in \kappa \rangle$. Define a
$\kappa$-iteration $\mathbb{P} = \mathbb{P}{_\kappa}$ of $\Gamma$
forcing notions, with appropriate support, as follows. At
successor stages, let $\mathbb{P}_{\alpha+1} =
\mathbb{P}_\alpha*\dot{\mathbb{Q}}_\alpha$, where, if $\alpha =
\pi( n,\mu)$ and $V_\delta^{\mathbb{P}_\alpha} \models$
``$\phi_n(x_\mu)$ is $\Gamma$-forceably necessary", then
$V_\delta^{\mathbb{P}_\alpha} \models$ ``$\mathbb{Q}_\alpha$ is a
forcing notion in $\Gamma$ that forces `$\phi_n(x_\mu)$ is
$\Gamma$-necessary' "; otherwise, $\mathbb{Q}_\alpha$ is
$\{\emptyset\}$, the trivial forcing.  Use appropriate support at
limit stages.  Note that, since $\cof(\delta) > \kappa$,
$\mathbb{P}_\alpha$ is in $V_\delta$ for all $\alpha < \kappa$.

Let $G$ be $V$-generic over $\mathbb{P}$.  The claim is that $V[G]
\models \MP_{\Gamma}(X)$.  To prove the claim, suppose $x \in X$
and $V[G] \models$ ``$\phi(x)$ is $\Gamma$-forceably necessary".
It will suffice to show that $V[G] \models$ `` $\phi(x)$ is
$\Gamma$-necessary".  Let $\phi = \phi_n$ and $x = x_\mu$ for some
$\alpha = \pi(n,\mu)$. By factoring, $V[G]=V[G_\alpha][G_{TAIL}]
\models$  ``$\phi(x)$ is $\Gamma$-forceably necessary", so
$V[G_\alpha] \models$
 ``$\phi(x)$ is $\Gamma$-forceably necessary" as well.  By elementarity,
$V_\delta[G_\alpha] \models$  ``$\phi(x)$ is $\Gamma$-forceably
necessary".  But at stage $\alpha$, the forcing notion
$\mathbb{Q}$ in $V_\delta$ has been defined to force
$\Box_{\Gamma}\phi_n(x_\mu)$. So $V_\delta[G_{\alpha+1}] \models$
``$\phi(x)$ is $\Gamma$-necessary". Again by elementarity,
$V[G_{\alpha+1}] \models$ ``$\phi(x)$ is $\Gamma$-necessary". And
since $V[G]$ is a $\Gamma$-forcing extension of $V[G_{\alpha+1}]$,
$V[G] \models$ ``$\phi(x)$ is $\Gamma$-necessary". This proves the
claim. Finally, since the iteration of $\mathbb{P}$ has
appropriate support, $\mathbb{P}$ is in $\Gamma$.  And since
$\cof(\delta) > \kappa$, $\mathbb{P}$ is in $V_\delta$.
\end{proof}

\noindent\textit{Remark}. A clear similarity can be seen between
\lemref{lem_MP_S4_X_one_step} and
\lemref{lem_sigma_Gamma_forc_nec}.  The key difference is that, in
\lemref{lem_sigma_Gamma_forc_nec}, we are establishing a
proto-maximality principle for a single formula, for which a
definition of the forcing relation is available.  Here, in
\lemref{lem_MP_S4_X_one_step}, we need a uniform definition of the
forcing relation for all formulas, hence the need for the set
model $V_\delta$.

\chapter{$\mathrm{MP_{\Psi}}$ and variations}
In this chapter we explore the properties of maximality principles
of the form $\MP_\Psi = \MP_{\Gamma_\Psi}$, where $\Gamma_\Psi$ is
the class of forcing notions preserving the truth of some
particular sentence $\Psi$. To introduce classes of this type, we
first discuss one sentence, $\CH$, the continuum hypothesis.  Let
$\Gamma_{\CH}$ be the class of all forcing notions $\mathbb{P}$
that preserve $\CH$. Thus, if it holds, then $\mathbb{P}$ is in
$\Gamma_{\CH}$ if and only if $V^\mathbb{P}\models \CH$.  And, if
$\CH$ does not hold, $\Gamma_{\CH}$ is the class of all forcing.
Let $\MPmod{\CH}$ denote the modified maximality principle
restricted to $\Gamma_{\CH}$.  We will denote the modal operator
$\Gamma_{\CH}$-necessary by $\Box_{\CH}$ and denote
$\Gamma_{\CH}$-forceable by $\Diamond_{\CH}$. As a warmup to more
general results, we will first show that the principle
$\MPmod{\CH}$ is actually logically equivalent to \MP\ itself
together with $\CH$. We show the two directions of implication in
two separate theorems, since they each generalize to apply to more
general classes of forcing notions.

\begin{lem}\label{lem gamma CH is adequate}
$\Gamma_{\CH}$ is adequate, and therefore an S4 forcing class.
\end{lem}

\begin{proof}
Clearly trivial forcing preserves \CH.  And if $\mathbb{P}$
preserves \CH\ and $V^\mathbb{P}\models$ ``$\mathbb{Q}$ preserves
\CH", then $\mathbb{P}*\dot{\mathbb{Q}}$ preserves \CH\ as well.
And these facts are $\Gamma_{\CH}$-necessary, since they are
provable in \ZFC.
\end{proof}

\textit{Remark.}  Since the modal logic $S4$ applies to the
necessity and possibility operators $\Box_{\CH}$ and
$\Diamond_{\CH}$, the principle $\MPmod{\CH}$, or
``$\Diamond_{\CH}\Box_{\CH}\phi$ implies $\phi$", is equivalent to
``$\Diamond_{\CH}\Box_{\CH}\phi$ implies $\Box_{\CH}\phi$". So for
the rest of this section the form of the maximality principle we
will use in proofs will be $$\Diamond_{\CH}\Box_{\CH}\phi
\textnormal{ implies } \phi,$$ that is, a formula which is
$\Gamma_{\CH}$-forceably necessary must be true.

\begin{thm}\label{thm_MPCH}
The principle $\MPmod{\CH}$ is equivalent to $\MP+\CH$.
\end{thm}

\begin{proof} To prove the forward implication, assume $\MPmod{\CH}$ ($\Diamond_{\CH}\Box_{\CH}\phi$
implies $\phi$, for all statements $\phi$).  We first claim that,
in fact, $\CH$ must be true: $\CH$ is certainly
$\Gamma_{\CH}$-forceable ($\Diamond_{\CH}\CH$).  And in the
extension where $\CH$ becomes true, $\CH$ is
$\Gamma_{\CH}$-necessary.  So $\CH$ is $\Gamma_{\CH}$-forceably
necessary ($\Diamond_{\CH}\Box_{\CH}\CH$), so by $\MPmod{\CH}$,
$\CH$ is true.

To establish \MP, suppose that $\phi$ is forceably necessary.  It
will suffice to infer that $\phi$ is true.  There exists a forcing
notion $\mathbb{P}$ such that $V^\mathbb{P}$ satisfies that $\phi$
is necessary.  Let $\mathbb{Q}$ be any forcing notion that forces
$\CH$ to be true in $V^{\mathbb{P}*\dot{\mathbb{Q}}}$.  Notice
that $\phi$ is necessary in $V^{\mathbb{P}*\dot{\mathbb{Q}}}$ as
well, since $V^\mathbb{P} \models \Box\phi$.  So
$V^{\mathbb{P}*\dot{\mathbb{Q}}}$ satisfies both $\Box\phi$ and
$\CH$.  Since $\mathbb{P}*\dot{\mathbb{Q}}$ is in $\Gamma_{\CH}$,
we have that $\Box\phi$ is $\Gamma_{\CH}$-forceable. But
$\Box\phi$ implies $\Box_{\CH}\phi$, so
$\Diamond_{\CH}\Box_{\CH}\phi$. Finally, by $\MPmod{\CH}$, $\phi$
is true.

To prove the reverse implication, suppose $\phi$ is
$\Gamma_{\CH}$-forceably necessary, that is, $\Box_{\CH}\phi$ is
$\Gamma_{\CH}$-forceable. $\Box_{\CH}\phi$ is therefore forceable.
Said another way, there is a forcing extension in which any
further forcing satisfies $\phi$, unless it is a forcing that is
not in $\Gamma_{\CH}$ and does not preserve $\CH$.  So $\phi \vee
\neg \CH$ is forceably necessary, and by \MP, it is true.  But
$\CH$ is true by hypothesis, so $\phi$ is true, showing that
$\MPmod{\CH}$ holds.
\end{proof}

\begin{thm}\label{thm_MPCH X}
Let $X$ be any set. The principle $\MPmod{\CH}(X)$ is equivalent
to $\MP(X)+\CH$.
\end{thm}

\begin{proof}
The same as the previous proof, since whether the formula
instances of $\MPmod{\CH}$ used in that proof had parameters or
not was irrelevant.
\end{proof}

\begin{cor}
The following are equivalent:
\begin{list}
 {(\arabic{listcounter})}{\usecounter{listcounter}\setlength{\rightmargin}{\leftmargin}}
\item
$Con(\ZFC+\MPmod{\CH})$
\item
$Con(\ZFC+\MP)$
\item
$Con(\ZFC)$.
\end{list}
\end{cor}

\noindent\textit{Remark}. This answers a question posed by Hamkins
in \cite{ham:max}.

\begin{proof}
$(2)\Longleftrightarrow(3)$ by \cite{ham:max}.  Since neither
$\CH$ nor its negation is forceably necessary by general forcing,
$\ZFC+\MP+\CH$ is equiconsistent with $\ZFC+\MP$.  Then
$(1)\Longleftrightarrow (2)$ by \thmref{thm_MPCH}.
\end{proof}

The preceding result can be generalized by abstracting the
essential properties of the classes of forcing notions just
discussed. If $\Gamma_1$ and $\Gamma_2$ are each classes of
forcing notions, when can we say that $\MP_{\Gamma_1} $ implies $
\MP_{\Gamma_2}$? The preceding result suggests possible answers.
Suppose $\Gamma_1$ and $\Gamma_2$ are classes of forcing notions.
Let us say that $\Gamma_2$ \textbf{absorbs}\footnote{The name is
vaguely inspired by the absorption law of propositional calculus.
Another way to express the relationship between $\Gamma_2$ and
$\Gamma_1$ is to say that $\Gamma_2$ is cofinal in $\Gamma_1$
where all forcing notions are ordered so that $\mathbb{P} \leq
\mathbb{P}*\dot{\mathbb{Q}}$ for any $\mathbb{Q}$.} $\Gamma_1$ if
it is true and $\Gamma_1$-necessary that $\Gamma_2$ is contained
in $\Gamma_1$ and for any $\mathbb{P}$ in $\Gamma_1$ there exists
$\mathbb{Q}$ such that $V^\mathbb{P}\models$ ``$\dot{\mathbb{Q}}$
is in $\Gamma_2$" and ``$\mathbb{P}*\dot{\mathbb{Q}}$ is in
$\Gamma_2$". One example of this just seen is the following.

\begin{thm}\label{thm Gamma CH absorbs all}
$\Gamma_{\CH}$ absorbs the class of all forcing notions.
\end{thm}

\begin{proof}
Let $\Gamma$ = all forcing notions.  Clearly ``$\Gamma_{\CH}
\subseteq \Gamma$" is $\Gamma$-necessary. Let $\mathbb{P}$ be in
$\Gamma$.  One can find $\mathbb{Q}$ in $V^\mathbb{P}$ to force
$\CH$, so $\mathbb{Q}$ is in $\Gamma_{\CH}$.  And since
$\mathbb{P}*\dot{\mathbb{Q}}$ forces $\CH$,
$\mathbb{P}*\dot{\mathbb{Q}}$ is in $\Gamma_{\CH}$.
\end{proof}

\noindent\textit{Remark}. Another example of a class of forcing
notions that behaves as $\Gamma_{\CH}$ does in \thmref{thm Gamma
CH absorbs all} is \COLL.  The forcing notion
$\coll(\omega,\theta)$ is the forcing notion that collapses
$\theta$ to $\omega$.  The forcing class \COLL\ consists of
$\coll(\omega,\theta)$, for all $\theta> \omega$.
\lemref{lem_COLL_absorbs_all_forcing} will show that \COLL\
absorbs the class of all forcing notions.

\begin{lem}[Absorption Lemma]\label{absorption lemma} Suppose $\Gamma_1$
and $\Gamma_2$ are classes of forcing notions such that $\Gamma_2$
absorbs $\Gamma_1$, and $\phi$ is any formula, with arbitrary set
parameters for $\phi$. If $\phi$ is $\Gamma_1$-forceably necessary
then it is $\Gamma_2$-forceably necessary as well.
\end{lem}

\noindent\textit{Remark}.  We allow arbitrary set parameters for
$\phi$, since this does not affect the argument.

\begin{proof}
Assume that $\phi$ is $\Gamma_1$-forceably necessary. Then there
is $\mathbb{P}$ in $\Gamma_1$ such that $V^{\mathbb{P}}\models$
``$\phi$ is $\Gamma_1$-necessary".  By absorption, there exists
$\dot{\mathbb{Q}}$ such that $V^{\mathbb{P}}\models
``\dot{\mathbb{Q}}$ is in $\Gamma_2$" and
$\mathbb{P}*\dot{\mathbb{Q}}$ is in $\Gamma_2$.  Since
$V^{\mathbb{P}}\models \Gamma_2\subseteq\Gamma_1$,
$V^{\mathbb{P}}\models$ ``$\dot{\mathbb{Q}}$ is in $\Gamma_1$ and
$\phi$ is $\Gamma_1$-necessary".   So
$V^{\mathbb{P}*\dot{\mathbb{Q}}}\models$ ``$\phi$ is
$\Gamma_1$-necessary".  And since
$V^{\mathbb{P}*\dot{\mathbb{Q}}}\models$ ``$\Gamma_2$ is contained
in $\Gamma_1$", $V^{\mathbb{P}*\dot{\mathbb{Q}}}\models$
``$\Box_{\Gamma_1}\phi$ implies $\Box_{\Gamma_2}\phi$".  So
$V^{\mathbb{P}*\dot{\mathbb{Q}}}\models$ ``$\phi$ is
$\Gamma_2$-necessary".  This establishes that $\phi$ is
$\Gamma_2$-forceably necessary.
\end{proof}

\begin{cor}\label{cor_absorbs_MP_Gamma}
Let $\Gamma_1$ and $\Gamma_2$ be S4 classes of forcing notions. If
$\Gamma_2$ absorbs $\Gamma_1$ then $\MP_{\Gamma_2}$ implies $
\MP_{\Gamma_1}$.
\end{cor}
\begin{proof}
Assume $\MP_{\Gamma_2}$ (in the S4 form
$\Diamond_{\Gamma_2}\Box_{\Gamma_2}\phi$ implies $\phi$). Next
suppose $\phi$ is $\Gamma_1$-forceably necessary, with the goal of
proving the truth of $\phi$. By absorption and \lemref{absorption
lemma}, $\phi$ is $\Gamma_2$-forceably necessary. This gives the
truth of $\phi$ by $\MP_{\Gamma_2}$.
\end{proof}

\begin{cor}\label{cor_absorbs_MP_Gamma_X}
Let $\Gamma_1$ and $\Gamma_2$ be S4 classes of forcing notions.
Let $X$ be any set.  If $\Gamma_2$ absorbs $\Gamma_1$ then
$\MP_{\Gamma_2}(X)$ implies $ \MP_{\Gamma_1}(X)$.
\end{cor}
\begin{proof}
Since \lemref{absorption lemma} allows formulas with arbitrary
parameters, just replace $\phi$ in the proof of
\corref{cor_absorbs_MP_Gamma} with $\phi(x)$, where $x$ is in $X$.
\end{proof}

\corref{cor_absorbs_MP_Gamma} corresponds exactly to how we proved
that $\MPmod{\CH}$ implies $\MP$ in \thmref{thm_MPCH}, where
$\Gamma_{\CH}$ absorbed the class of all forcing notions.

To generalize the reverse direction of \thmref{thm_MPCH}, one
generalizes the $\CH$-preserving property. Let $\Psi$ be a
sentence of the language of \ZFC.  A forcing notion $\mathbb{P}$
is $\Psi$-\textbf{preserving} if $\Psi$ implies that it is true in
$V^{\mathbb{P}}$. I will denote the class of $\Psi$-preserving
forcing notions by $\Gamma_\Psi$.  The associated maximality
principle is $\MP_\Psi$.

\begin{lem}\label{lem gamma Psi is adequate}
Let $\Psi$ be a sentence of the language of \ZFC. $\Gamma_{\Psi}$
is S4.
\end{lem}

\begin{proof}
This proof emulates the proof of \lemref{lem gamma CH is
adequate}.

Clearly trivial forcing preserves $\Psi$.  And if $\mathbb{P}$
preserves $\Psi$ and $V^\mathbb{P}\models$ ``$\mathbb{Q}$
preserves $\Psi$", then $\mathbb{P}*\dot{\mathbb{Q}}$ preserves
$\Psi$ as well. And these facts are $\Gamma_{\Psi}$-necessary,
since they are provable in \ZFC.
\end{proof}

\begin{lem}\label{lem MP Psi implies S4}
Let $\Psi$ be a forceable sentence of the language of \ZFC. The
principle $\MP_{\Psi}$ implies $\Psi$.
\end{lem}

\begin{proof}
If $\neg \Psi$ then $\Gamma_{\Psi}$ is all forcing.  Since $\Psi$
is forceable, $\Diamond_{\Psi}\Psi$. And in the extension where
$\Psi$ becomes true, $\Psi$ is $\Gamma_{\Psi}$-necessary. So
$\Diamond_{\Psi}\Box_{\Psi}\Psi$, so by $\MP_{\Psi}$, $\Psi$ is
true.
\end{proof}

\begin{thm}\label{thm_MP_Psi_preserve}
Let $\Psi$ be a forceable sentence of the language of \ZFC.  Let
$\Gamma$ be an S4 class of forcing notions. If $\Gamma_\Psi
\subseteq \Gamma$ is $\Gamma$-necessary then $\MP_{\Gamma}$+$\Psi$
implies $ \MP_\Psi$.
\end{thm}
\begin{proof}
Assume $\MP_{\Gamma}$.  By S4 we can write it in the form that
says $\Diamond_{\Gamma}\Box_{\Gamma}\phi$ implies $\phi$. Next
suppose $\phi$ is $\Gamma_\Psi$-forceably necessary, with the goal
of proving $\phi$ from this.  Since $\Gamma_\Psi \subseteq
\Gamma$, we have $\phi$ is $\Gamma$-forceably
$\Gamma_\Psi$-necessary. But $\Box_{\Gamma_\Psi}\phi$ implies
$\Box_{\Gamma}(\phi \vee \neg \Psi)$, since any forcing notion in
$\Gamma$ that is not in $\Gamma_\Psi$ must produce a model in
which $\Psi$ is false. Combining the last two sentences gives
$\Diamond_{\Gamma}\Box_{\Gamma}(\phi \vee \neg \Psi)$. Now we
apply $\MP_{\Gamma}$ to get $\phi \vee \neg \Psi$. But $\Psi$ is
true by hypothesis, hence $\phi$ is true.
\end{proof}

The concept that unites both directions of implication is
resurrectibility.  A sentence $\Psi$ of the language of \ZFC\ is
$\Gamma$-\textbf{resurrectible} if there is $\Gamma$-necessarily a
forcing notion $\mathbb{P} \in \Gamma$ such that $V^{\mathbb{P}}
\models \Psi$, in other words, if
$\Box_{\Gamma}\Diamond_{\Gamma}\Psi$ is true.

\noindent\textit{Remark}.  \CH\ is resurrectible---it can be
forced to be true over any model of \ZFC.

\begin{thm}\label{lem_resurrect_MP_psi}
Let $\Gamma$ be an S4 class of forcing notions. If $\Psi$ is a
$\Gamma$-resurrectible sentence of the language of \ZFC\ and
$\Gamma_\Psi\subseteq \Gamma$ is $\Gamma$-necessary, then
$\MP_\Psi+\Psi \Longleftrightarrow \MP_{\Gamma}+\Psi$.
\end{thm}

\begin{proof}
Clearly $\Gamma_\Psi$ absorbs $\Gamma$, since any forcing notion
that resurrects $\Psi$ must be in $\Gamma_\Psi$.
\corref{cor_absorbs_MP_Gamma} then gives the forward implication.
The converse follows directly from \thmref{thm_MP_Psi_preserve}.
\end{proof}

This theorem can be applied to obtain useful equiconsistency
results.

The next lemma provides some basic facts.
\begin{lem}
Let $\phi$, $\phi_1$, and $\phi_2$ be sentences in the language of
\ZFC.
\begin{list}
 {(\arabic{listcounter})}{\usecounter{listcounter}\setlength{\rightmargin}{\leftmargin}}
\item
If $\phi_1 $ implies $ \phi_2$, then $\Gamma_{\phi_1} \subseteq
\Gamma_{\phi_2}$
\item
If $\phi$ is provable in \ZFC, then $\phi$ is $\Gamma$-necessary,
for any class $\Gamma$ of forcing notions.
\item
If $\phi$ is provably $\Gamma$-forceable, then $\phi$ is
$\Gamma$-resurrectible.
\end{list}
\end{lem}

\begin{proof}
(1) and (2) are obvious (in modal logic terms, (2) is just the
necessitation principle), and (3) follows from(2).
\end{proof}

Another basic fact is
\begin{lem}\label{lem_provably_forceable}
If $\Psi$ is provably $\Gamma$-forceable, then
$Con(\ZFC+\MP_\Gamma)$ if and only if $Con(\ZFC+\MP_\Gamma+\Psi)$
\end{lem}
\begin{proof}
Let $M \models \MP_\Gamma+ \ZFC$.  Let $M[G]$ be a
$\Gamma$-forcing extension of $M$ such that $M[G] \models \Psi$.
But then $M[G] \models \MP_\Gamma$ as well:  If $M[G] \models$
``$\phi$ is $\Gamma$-forceably necessary" then $M \models$
``$\phi$ is $\Gamma$-forceably necessary", so $M \models$ ``$\phi$
is $\Gamma$-necessary" and therefore $M[G] \models$ ``$\phi$ is
$\Gamma$-necessary".
\end{proof}

Recall that Suslin trees, which are provably forceable either by
Cohen forcing or $<\omega_1$-closed forcing, negate Suslin's
Hypothesis by their existence.  By the above facts, the sentence
$\neg SH$ is resurrectible in the class of all forcing notions.
(In fact, any resurrectible sentence could be used in the
following theorem.)

\begin{thm}\label{thm con_MP_neg_SH}
The following are equivalent:
\begin{list}
 {(\arabic{listcounter})}{\usecounter{listcounter}\setlength{\rightmargin}{\leftmargin}}
\item
$Con(\ZFC+\MP_{\neg SH})$
\item
$Con(\ZFC+\MP)$
\item
$Con(\ZFC)$.
\end{list}
\end{thm}

\begin{proof}
$(1) \Leftrightarrow (2)$: Via \thmref{lem_resurrect_MP_psi} and
\lemref{lem_provably_forceable}.

$(2) \Leftrightarrow (3)$: Via \cite{ham:max}.
\end{proof}

\chapter{$\mathrm{MP_{CCC}}$ and variations}

\section{The consistency strength of $\mathrm{MP_{CCC}}(\mathbb{R})$}
In this section we show a surprising result.  \cite{ham:max} shows
that the maximality principle with real parameters,
$\MP(\mathbb{R})$, has consistency strength strictly greater than
\ZFC, while $\MP$ and $\MPmod{\CCC}$ are both equiconsistent with
\ZFC.  So it would seem that adding the parameter set $\mathbb{R}$
should increase the consistency strength of $\MPmod{\CCC}$ as it
did for $\MP$, especially since \CCC-forcing will certainly add
new reals that are not in the ground model. However, contrary to
expectation, we have the following theorem:

\begin{thm}\label{thm_MP_ccc_R} The following are equivalent:
\begin{list}
 {(\arabic{listcounter})}{\usecounter{listcounter}\setlength{\rightmargin}{\leftmargin}}
\item
$Con(\ZFC)$
\item
$Con(\ZFC+\MPmod{\CCC}(\mathbb{R}))$
\end{list}
\end{thm}

\begin{proof}
$(2)$ implies $ (1)$: trivial.

$(1)$ implies $(2)$:  It will suffice to show that, given a model
of \ZFC, one can produce a model of
$\ZFC+\MPmod{\CCC}(\mathbb{R})$. I will give two proofs of this.

To begin the first proof, we prove the consistency of a weak
version of this principle.  For any set $X$, using our notation,
$\MPmod{\CCC}(X)$ is the modified maximality principle that says
any formula with parameters taken from the set $X$ which is
$\CCC$-forceably necessary is true.  Let $\mathbb{P}$ be a $\CCC$
forcing notion.  Let us confine ourselves to the model
$V^\mathbb{P}$, and denote by $\mathbb{R}^V$ the set of reals of
$V^\mathbb{P}$ which are not introduced by forcing with
$\mathbb{P}$, that is, all reals of the ground model $V$. Let the
principle $\MPmod{\CCC}(\mathbb{R}^V)$ be the form of
$\MPmod{\CCC}(X)$ interpreted in $V^\mathbb{P}$ with parameter set
$\mathbb{R}^V$.

\begin{lem}\label{lem_MP_ccc_R_one_step}
Suppose $V_\delta \prec V$, and $\cof(\delta) > 2^\omega$. Then
there is a forcing notion $\mathbb{P}$ in $\CCC$ such that
$V^{\mathbb{P}} \models \MPmod{\CCC}(\mathbb{R}^V)$ and
$\mathbb{P} \in V_\delta$.
\end{lem}

\begin{proof}
This is an instance of \lemref{lem_MP_S4_X_one_step}.
\end{proof}

\begin{lem}\label{lem_MP_ccc_R}
If there is a model of \ZFC\ then there is a model of
$\ZFC+\MPmod{\CCC}(\mathbb{R})$.
\end{lem}

\begin{proof}[First proof]
Suppose $V\models \ZFC$.  We will construct a forcing extension
which is a model of $\ZFC+\MPmod{\CCC}(\mathbb{R})$. By
\thmref{thm_Vdelta_submod_V_club_class} we may assume that there
is a class club of cardinals $C$ such that for all $\delta$ in
$C$, $V_\delta\prec V$.  Construct a finite-support
$\omega_1$-iteration $\mathbb{P}=\mathbb{P}_{\omega_1}$, such that
$V^{\mathbb{P}} \models \MPmod{\CCC}(\mathbb{R})$, as follows. Let
$\mathbb{P}_0$ be the trivial notion of forcing.  At stage
$\alpha$, select $\delta_\alpha$ from the club $C$ such that the
rank of $\mathbb{P}_\alpha < \delta_\alpha$ (so that
$\mathbb{P}_\alpha$ is in $V_{\delta_\alpha}$) and
$\cof(\delta_\alpha)
> (2^\omega)^{V^{\mathbb{P}_\alpha}}$. Working in
$V^{\mathbb{P}_\alpha}$, define
$\mathbb{P}_{\alpha+1}=\mathbb{P}_\alpha*\dot{\mathbb{Q}}_\alpha$,
where $\dot{\mathbb{Q}}$ is a $\mathbb{P}_\alpha$-name of a $\CCC$
notion of forcing such that
$V^{\mathbb{P}_\alpha*\dot{\mathbb{Q}}_\alpha} \models
\MPmod{\CCC}(\mathbb{R}^{V^{\mathbb{P}_\alpha}})$. Such a
$\mathbb{Q}_\alpha$ is guaranteed to exist by
\lemref{lem_MP_ccc_R_one_step}, since the conditions
$V_{\delta_\alpha} \prec V$ and $ \cof(\delta_\alpha)> 2^\omega$,
are satisfied by $V^{\mathbb{P}_\alpha}$.  This completes the
construction of $\mathbb{P}=\mathbb{P}_{\omega_1}$.

Let $G$ be $V$-generic over $\mathbb{P}$.  I claim that $V[G]
\models \MPmod{\CCC}(\mathbb{R})$.  To see this, let $V[G]
\models$ ``$\phi(r)$ is $\CCC$-forceably necessary", where $r \in
\mathbb{R}$, the reals as interpreted in $V[G]$.  It will suffice
to show that $V[G] \models$ ``$\phi(r)$ is $\CCC$-necessary". Note
that $r$, as a real, is a subset of $\omega$ in $V[G]$, while
$\omega$ itself is in $V$. Therefore, since $\mathbb{P}$ is an
$\omega_1$ iteration with finite support and $\cof(\omega_1) >
|\omega|$, $r$ must be in some $V[G_\alpha]$, where $\alpha <
\omega_1$, $\mathbb{P} =
\mathbb{P}_\alpha*\dot{\mathbb{P}}^{(\alpha)}_{\omega_1}$ is the
factorization of $\mathbb{P}$ at stage $\alpha$, and $G_\alpha$ is
the projection of $G$ to $\mathbb{P}_\alpha$. (This follows from
\cite{kun:sth}, Chapter VIII, Lemma 5.14.)  But the definition of
$\mathbb{P}$ required that $\mathbb{Q}_\alpha$ force
$\MPmod{\CCC}(\mathbb{R}^{V[G_\alpha]})$, which therefore must
hold at stage $\alpha+1$.  Indeed, refactoring $\mathbb{P} =
\mathbb{P}_{\alpha+1}*\dot{\mathbb{P}}^{(\alpha+1)}_{\omega_1}$
and setting $G_{\alpha+1}$ to be the projection of $G$ to
$\mathbb{P}_{\alpha+1}$, we have that $r \in V[G_\alpha]$ and
$V[G_{\alpha+1}] \models \MPmod{\CCC}(\mathbb{R}^{V[G_\alpha]})$
(as well as $V[G_{\alpha+1}] \models$ ``$\phi(r)$ is
$\CCC$-forceably necessary", since
$V[G]=V[G_{\alpha+1}][G^{(\alpha+1)}_{\omega_1}]$ is a
$\CCC$-forcing extension of $V[G_{\alpha+1}]$). Therefore
$V[G_{\alpha+1}] \models$ ``$\phi(r)$ is $\CCC$-necessary".  Since
$V[G]$ is a $\CCC$-forcing extension of $V[G_{\alpha+1}]$, we have
that $V[G] \models$ ``$\phi(r)$ is $\CCC$-necessary", as required.
\end{proof}

\begin{proof}[Second proof]
This time, we use a bookkeeping function style argument.  Again,
suppose $V\models \ZFC$. By
\thmref{thm_Vdelta_submod_V_club_class} we may assume that there
is in $V$ a class club of cardinals $C$ such that for all $\delta$
in $C$, $V_\delta\prec V$.   Let $\pi:ORD \simeq \omega\times ORD
\times ORD$ be a definable bijective class function $\pi: \alpha
\mapsto \langle n,\beta,\mu\rangle$ such that $\beta \leq \alpha$.
Using $\pi$ as a bookkeeping function, we define a sequence of
iterated forcing notions $\mathbb{P}_\alpha$, simultaneously with
a sequence of cardinals $\delta_\alpha$ by transfinite induction
on $\alpha$ in $ORD$, as follows.  Let $\mathbb{P}_0$ be trivial
forcing.  Given $\mathbb{P}_\alpha$, let $\delta_\alpha$ be the
least cardinal in the club $C$ such that $\mathbb{P}_\alpha$ is in
$V_{\delta_\alpha}$.  Define $\dot{\mathbb{Q}}_\alpha$ in
$V_{\delta_\alpha}^\mathbb{P_\alpha}$ as follows: Let $\pi(\alpha)
= \langle n,\beta,\mu\rangle$.  Consider the statement $\phi(x) =
\phi_n(x)$, the $n^{th}$ statement in the language of \ZFC\
according to some enumeration, with single parameter $x = x_\mu$,
the $\mu^{th}$ name for a real in the model $V^{\mathbb{P}_\beta}$
where $\beta \leq \alpha$.\footnote{Assuming an enumeration of the
names for reals in $V^{\mathbb{P}_\beta}$ is based on a
wellordering of the ground model, which can be forced by class
forcing that adds no new sets.  This induces a wellordering of
names, giving a wellordering of any forcing extension.} If, in
$V_{\delta_\alpha}^\mathbb{P_\alpha}$, $\phi(x)$ is
$\CCC$-forceably necessary, let $\dot{\mathbb{Q}}_\alpha$ be the
$V_{\delta_\alpha}$-least $\mathbb{P}_\alpha$-name of a forcing
notion which performs a forcing that $\phi(x)$ is
$\CCC$-necessary. Otherwise let $\dot{\mathbb{Q}}_\alpha$ be the
$\mathbb{P}_\alpha$-name for trivial forcing.  Now let
$\mathbb{P}_{\alpha+1} = \mathbb{P}_\alpha*
\dot{\mathbb{Q}}_\alpha$.  Finally, take finite support at limits.
This defines the sequence $\mathbb{P}_\alpha$ for all $\alpha$ in
$ORD$.  Note that, for all such $\alpha$, $\mathbb{P}_\alpha$ is
$\CCC$ and is contained in $V_{\delta_\alpha}$.

We wish to truncate this sequence at an appropriate length
$\lambda$ to obtain an iterated forcing notion
$\mathbb{P}_\lambda$ which forces a model of
$\MPmod{\CCC}(\mathbb{R})$.  This will occur if all reals in
$V^{\mathbb{P}_\lambda}$ are introduced at some earlier stage of
the iteration and the cofinality of $\lambda$ is greater than
$\omega$.  To ensure this, we define $\lambda$ to be a closure
point of the function $f:ORD\longrightarrow ORD$ which takes
$\beta$, the stage at which a real parameter is introduced, to the
least stage by which all formulas $\phi$ have been applied to all
parameters in $V^{\mathbb{P}_\beta}$.  We now make use of the
technique of using \textit{nice names}, discussed in detail in
\cite{kun:sth}.  Since we only need to count nice names, of which
there are $|\mathbb{P}_\beta|^\omega$ many in
$V^{\mathbb{P}_\beta}$, this gives $f(\beta) =
\mathrm{sup}_{\mu<|\mathbb{P}_\beta|^\omega} \{\pi(\alpha) =
\langle n, \beta, \mu\rangle\}$.  Now let $\lambda$ be the first
closure point of $f:ORD\longrightarrow ORD$ with cofinality
$\omega_1$. (An ordinal $\alpha$ is a \textbf{closure point} of f
if $f``\alpha \subseteq \alpha$.)

Let $\mathbb{P} = \mathbb{P}_\lambda$, and let $G$ be $V$-generic
over $\mathbb{P}$.  By the usual argument, we can now establish
that $V[G] \models \MPmod{\CCC}(\mathbb{R})$: Suppose $V[G]
\models$ ``$\phi(r)$ is $\CCC$-forceably necessary".  Then there
is $\alpha = \langle n, \mu, \beta \rangle$, where $\phi = \phi_n$
and $\dot{r}$ is the $\mu^{th}$ nice $\mathbb{P}_\beta$-name of a
real, for some $\beta \leq \alpha$ (the name $\dot{r}$ appears
before stage $\delta$ by the reasoning used in
\lemref{lem_MP_S4_X_one_step}).  Since $V[G]$ is a $\CCC$-forcing
extension of $V[G_\alpha]$, where $G_\alpha$ is $V$-generic over
$\mathbb{P}_\alpha$, $V[G_\alpha] \models$ ``$\phi(r)$ is
$\CCC$-forceably necessary", whence by elementarity,
$V_{\delta_\alpha}[G_\alpha] \models$ ``$\phi(r)$ is
$\CCC$-forceably necessary".  But by the construction of
$\mathbb{P}$, if $G_{\alpha+1}$ is $\mathbb{P}_{\alpha+1}$-generic
over $V_{\delta_{\alpha+1}}$ then
$V_{\delta_{\alpha+1}}[G_{\alpha+1}] \models$ ``$\phi(r)$ is
$\CCC$-necessary", so by elementarity $V[G_{\alpha+1}] \models$
``$\phi(r)$ is $\CCC$-necessary" and therefore $V[G] \models$
``$\phi(r)$ is $\CCC$-necessary".
\end{proof}

This concludes the proof of \thmref{thm_MP_ccc_R}.
\end{proof}

\noindent\textit{Remark}.  This raises the possibility that
$\MPmod{\CCC}$ directly implies $\MPmod{\CCC}(\mathbb{R})$.  But
we will see that this is false.

\noindent\textit{Remark}.  The proof just given of
\lemref{lem_MP_ccc_R} makes use of the existence of closure points
of the defined function $f:ORD\longrightarrow ORD$.  In order to
know such closure points exist one needs to apply the Replacement
scheme.  Even the first proof makes use of the Replacement Axiom
Scheme enhanced with the symbol $C$, in order to construct the
iteration $\mathbb{P}_\alpha$.  These arguments take place in the
language of \ZFC\ expanded with the symbol $C$ interpreted as a
class club in an expanded model.  This is why we included, in
\thmref{thm_Vdelta_submod_V_club_class}, all instances of
Replacement and Comprehension that mention the class club $C$.

One might expect that \thmref{thm_MP_ccc_R} can be extended to
parameter sets which are power sets of sets of cardinality greater
than $\omega$, such as $\omega_1$ or $\aleph_{\omega_{17}}$, and
in fact this is the case. Let $\kappa$ be a cardinal.  Singling
out the second proof strategy above, one can state and prove the
generalizations of \lemref{lem_MP_ccc_R} and \thmref{thm_MP_ccc_R}
as:

\begin{lem}\label{lem_MP_ccc_kappa}
Let $\kappa$ be any $\CCC$-absolutely definable cardinal.  If
there is a model of \ZFC\ then there is a model of
$\ZFC+\MPmod{\CCC}(H(\kappa))$\end{lem}

\begin{proof}
As before, the proof consists in finding a $\CCC$-forcing
extension model for $\MPmod{\CCC}(H(\kappa))$.  Let $V$ be a model
of \ZFC.  We use the same bookkeeping class function and again
assume in $V$ a class club of cardinals $C$ such that for all
$\delta$ in $C$, $V_\delta\prec V$.  The definition of a sequence
of iterated forcing notions $\mathbb{P}_\alpha$, simultaneously
with a sequence of cardinals $\delta_\alpha$ again proceeds by
transfinite induction on $\alpha$ in $ORD$.   At stage $\alpha$,
in defining  $\dot{\mathbb{Q}}_\alpha$ in
$V_{\delta_\alpha}^\mathbb{P_\alpha}$, $\alpha$ now codes $
\langle n,\beta,\mu\rangle$ where $\mu$ is now an index for the
name for a subset $x$ of $\kappa$ in the model
$V^{\mathbb{P}_\beta}$ where $\beta \leq \alpha$, and
$\dot{\mathbb{Q}}_\alpha$ forces that $\phi(x)$ is
$\CCC$-necessary if such a forcing notion exists.  Again let
$\mathbb{P}_{\alpha+1} =
\mathbb{P}_\alpha*\dot{\mathbb{Q}}_\alpha$ with finite support at
limits. This again gives $\mathbb{P}_\alpha$ for all $\alpha$ in
$ORD$, with $\mathbb{P}_\alpha$ being $\CCC$ and contained in
$V_{\delta_\alpha}$.

The length $\lambda$ of the iteration giving a
$\mathbb{P}_\lambda$ which satisfies $\MPmod{\CCC}(H(\kappa))$
must now have cofinality $\kappa$, to ensure that all names of
parameters in $V^{\mathbb{P}_\lambda}$ are introduced at earlier
stages.  It must also be a closure point of the function
$f:ORD\longrightarrow ORD$ which takes $\beta$, the stage at which
a real parameter is introduced, to the least stage by which all
formulas $\phi$ have been applied to all parameters in
$V^{\mathbb{P}_\beta}$.  Since we only need to count nice names,
of which there are $(|\mathbb{P}_\beta|^\omega)^{<\kappa} =
|\mathbb{P}_\beta|^{<\kappa}$ many in $V^{\mathbb{P}_\beta}$, this
gives $f(\beta) = \mathrm{sup}_{\mu<|\mathbb{P}_\beta|^{<\kappa}}
\{\alpha = \langle n, \beta, \mu\rangle\}$.

Let $\mathbb{P} = \mathbb{P}_\lambda$, and let $G$ be $V$-generic
over $\mathbb{P}$.  By an analogous argument, $V[G] \models
\MPmod{\CCC}(H(\kappa))$.  Since $\cof(\lambda) > \kappa$, all
parameters in $V[G]$ appear at some previous stage. The rest of
the argument is identical to that for \lemref{lem_MP_ccc_R}.
\end{proof}

\begin{thm}\label{thm_MP_ccc_kappa}  Let $\kappa$ be any $\CCC$-absolutely definable cardinal.  Then
the following are equivalent:
\begin{list}
 {(\arabic{listcounter})}{\usecounter{listcounter}\setlength{\rightmargin}{\leftmargin}}
\item
$Con(\ZFC)$
\item
$Con(\ZFC+\MPmod{\CCC}(H(\kappa)))$
\end{list}
\end{thm}

\begin{proof}
Immediate from \lemref{lem_MP_ccc_kappa}
\end{proof}

\noindent\textit{Remark}.  This raises the possibility that
$\MPmod{\CCC}$ directly implies $\MPmod{\CCC}(H(\kappa))$ for any
$\CCC$-absolutely definable cardinal $\kappa$.  But we will see
that this is false.

Finally, we prove what appears to be an optimal result in this
direction. It is optimal in the sense that
$\MPmod{\CCC}(H(2^\omega))$ is known to be equiconsistent with an
inaccessible cardinal, as proved in Lemma 5.9 in \cite{ham:max},
and the following proof is a straightforward modification of that
lemma by only dropping the weak inaccessibility (i.e. regularity)
of $\delta = 2^\omega$.

\begin{thm}\label{thm_MP_ccc_H_cof_2_omega}
The following are equivalent:
\begin{list}
 {(\arabic{listcounter})}{\usecounter{listcounter}\setlength{\rightmargin}{\leftmargin}}
\item
$Con(\ZFC)$
\item
$Con(\ZFC+\MPmod{\CCC}(H(\cof(2^\omega)))$
\end{list}
\end{thm}

\begin{proof}
That $(2)$ implies $(1)$ is trivial.  For the converse suppose
without loss of generality that $V_\delta \prec V$, since this is
equiconsistent with \ZFC.  Define a bijective bookkeeping function
$\delta \longrightarrow \omega\times\delta\times\delta$ such that
if $\alpha = \langle n,\mu,\beta\rangle$ then $\beta < \alpha$.
Enumerate the formulas having one free variable in the language of
\ZFC\ as $\langle \phi_n|n \in \omega\rangle$.  Define a finite
support $\delta$-iteration of $\CCC$ forcing
$\mathbb{P}=\mathbb{P}_\delta$ as follows. Let $\mathbb{P}_0$ be
trivial forcing.  Given $\mathbb{P}_\alpha$, let $\alpha = \langle
n,\mu,\beta\rangle$.  If $\phi_n(x)$ is $\CCC$-forceably necessary
in $V_\delta^{\mathbb{P}_\alpha}$ with $x$ the $\mu^{th}$
$\mathbb{P}_\beta$-name of a parameter from $H(\cof(2^\omega))$ in
$V_\delta$, then let $\mathbb{Q}_\alpha$ be a $\CCC$ forcing
notion such that $V_\delta
^{\mathbb{P}_\alpha\star\dot{\mathbb{Q}}_\alpha}\models$
``$\phi_n(x)$ is $\CCC$-necessary".  Otherwise let
$\mathbb{Q}_\alpha$ be trivial forcing.

Let $G\subseteq \mathbb{P}_\delta$ be $V$-generic for
$\mathbb{P}_\delta$.  We claim $V[G]\models
\MPmod{\CCC}(H(\cof(\delta))$.  Also, the $\CCC$ iteration
$\mathbb{P}_\delta$ preserves the cofinality of $\delta$.

Now let $V[G]\models$ ``$\phi(x)$ is $\CCC$-forceably necessary",
where $x$ is in $H(\cof(\delta))$. It suffices to show that
$V[G]\models$ ``$\phi(x)$ is $\CCC$-necessary".  Since $x$ is in
$H(\cof(\delta))$, it appears at some stage of the iteration,
$\mathbb{P}_\beta$, as the $\mu^{th}$ parameter from $H(\delta)$
in $V^{\mathbb{P}_\beta}$, and $\phi = \phi_n$ for some $n$.  Let
$\alpha = \langle n,\mu,\beta\rangle$ be less than $\delta$.  Then
$V[G]$ is a $\CCC$-extension of $V[G_\alpha]$, where $G_\alpha$ is
$V$-generic over $\mathbb{P}_\alpha$.  So $V[G_\alpha]\models$
``$\phi(x)$ is $\CCC$-forceably necessary". By elementarity,
$V_\delta[G_\alpha]\models$ ``$\phi(x)$ is $\CCC$-forceably
necessary".  But by the construction of $\mathbb{P}_\delta$,
$V_\delta[G_{\alpha+1}]\models$ ``$\phi(x)$ is $\CCC$-necessary",
where $G_{\alpha+1}$ is $V$-generic over $\mathbb{P}_{\alpha+1}$.
By elementarity, $V[G_{\alpha+1}]\models$ ``$\phi(x)$ is
$\CCC$-necessary".  And since $V[G]$ is a $\CCC$-extension of
$V[G_{\alpha+1}]$, $V[G]\models$ ``$\phi(x)$ is $\CCC$-necessary".

Finally, the $\CCC$ iterations below stage $\delta$ will leave
$2^\omega$ less than $\delta$ (being within $V_\delta$, as an
elementary submodel of $V$) but otherwise arbitrarily large (for
any given $\kappa <\delta$, it is \CCC-forceably necessary that
$2^\omega> \kappa$). So $V[G]\models 2^\omega \geq \delta$. But
forcing of size $\leq\delta$ imposes the restriction $V[G]\models
2^\omega \leq \delta$.  So $V[G]\models 2^\omega = \delta$, and
therefore $V[G]$ models $\MPmod{\CCC}(H(\cof(2^\omega))$.
\end{proof}

\noindent\textit{Remark.} Notice that, by modifying the proof to
invoke \thmref{thm_Vdelta_submod_V_club_class}, we can make
$\cof(\delta)$ as large as we like.  Let's call the desired
cofinality $\kappa$.  Just choose $\delta$ from the class club $C$
to have the desired cofinality, by picking the $\kappa^{th}$
element in $C$.  For example, if $\cof(\delta)$ is chosen to be
$\omega_{23}$, we get a model of
$\ZFC+\MPmod{\CCC}(H(\omega_{23}))+\cof(2^\omega)=\omega_{23}$, an
improvement over \thmref{thm_MP_ccc_kappa}.

\noindent\textit{Remark.} The construction in the proof of this
theorem actually provides a stronger parameter set than
$H(cof(2^\omega))$.  At each stage $\alpha$, all parameters from
$H(\delta)$ are allowed, as seen from $V^{\mathbb{P}_\alpha}$.
Notice that such stages exist for all $\alpha<\delta=2^\omega$, so
if $\eta$ is an ordinal below $\alpha$, it is seen by
$V^{\mathbb{P}_\beta}$ as being in $H(\delta)$, for $\beta
<\alpha$.  Since all ordinals below $\delta=2^\omega$ eventually
achieve this status, the set of ordinals below $2^\omega$, i.e.,
the set $2^\omega$, can be used as a parameter set.  So the actual
result of the theorem is $Con(\ZFC)\Longleftrightarrow
Con(\ZFC+\MPmod{\CCC}(H(\cof(2^\omega)\cup 2^\omega))$.

\noindent\textit{Remark.} A similar result is given as theorem 3.9
in \cite{sta:vaa}, but there is a subtle difference between that
theorem and \thmref{thm_MP_ccc_H_cof_2_omega}. Stavi and
V\"a\"an\"anen discuss maximality principles with various
parameter sets including $H(2^\omega)$ and $H(\cof(2^\omega))$, as
well as various classes of forcing notions, including $\CCC$ and
Cohen forcing.  Their statements of maximality principles for a
class $\Gamma$ apply to statements in the language of \ZFC\ which
are provably $\Gamma$-persistent. (A statement $\phi$ is
$\Gamma$-persistent if when $\phi$ is true, then it is true in all
$\Gamma$-forcing extensions.)  They also seem to have provided the
first proofs of maximality principles by iterating over all
instances.  The lineage of such proofs, tracing citations, goes
from \cite{sta:vaa} through \cite{asp:th} to \cite{ham:max}.

\newcounter{bean}
\section{$\Box \mathrm{MP_{CCC}}$}

One certainly expects the consistency strength of $\Box
\MPmod{\CCC}$ with parameters to be significantly higher than $
\MPmod{\CCC}$. Again we can informally argue this case by making
analogies with \MP\ and $\Box\MP(\mathbb{R})$, the consistency
strength of the latter being at least that of the Axiom of
Determinacy, according to a result of Woodin (see \cite{ham:max}).
We recall the point made in the introduction that an uncountable
parameter is not acceptable for the principle $\MP$ (since the
assertion of its uncountability will be falsified by the principle
since its cardinality can be collapsed to $\omega$). With \CCC\
forcing this problem does not occur, so it seems larger parameter
sets should be acceptable. But when we attempt to include
uncountable parameters by using the parameter set $H(\omega_2)$,
the result is unexpected.

\beginThm\label{thm-1}
Let $\Gamma$ be a class of forcing notions necessarily containing
$\CCC$. Then $\Box \MP_{\Gamma}(H(\omega_2))$ is false.
\end{thm}

To prove this theorem, we will utilize some well-known connections
among $\CCC$ forcing, Martin's Axiom (MA), and Suslin trees.
Recall that an ($\omega_1$-) Suslin tree is an $\omega_1$-tree
which has no uncountable chain or antichain.

\begin{lem}\label{lem MA true}
Let $\Gamma$ be a class of forcing notions necessarily containing
$\CCC$. Then $\MP_{\Gamma}(H(\omega_2))$ implies $MA_{\omega_1}$.
\end{lem}

\begin{proof}
Recall that $MA_{\omega_1}$ is the statement that for any forcing
notion $\mathbb{P}$ which is $\CCC$, and any family $\mathcal{D}$
of sets dense in $\mathbb{P}$, where $|\mathcal{D}| = \omega_1$,
there is a filter $G$ meeting each $D \in \mathcal{D}$.  It is a
fact \cite{kun:sth} that to establish $MA_{\omega_1}$ it suffices
to consider forcing notions of size $ \leq \omega_1$.  Assume
$\mathbb{P}$ is a $\CCC$ forcing notion such that $|\mathbb{P}| <
\omega_1$. Assume also that $\mathcal{D}$ is a family of dense
sets in $\mathbb{P}$ where $|\mathcal{D}| = \omega_1$. Let
$\phi(p,d)$ be the assertion that there is a filter $G$ in the
forcing notion $p$ meeting each $D$ in $d$.
$\phi(\mathbb{P},\mathcal{D})$ is included in the scope of this
principle, since the two parameters can be regarded as a single
parameter which is in $H(\omega_2)$. (It is an ordered pair whose
components are each of size $\omega_1$.)

\begin{claim} $\phi(\mathbb{P},\mathcal{D})$ is $\Gamma$-forceably
necessary.
\end{claim}
\begin{proof}
We will find a $\Gamma$-forcing extension (in fact, a
$\CCC$-forcing extension) in which $\phi(\mathbb{P},\mathcal{D})$
is necessary.  We do this by regarding $\mathbb{P}$ as a forcing
notion where $p$ extends $q$ iff $p \leq q$ in $\mathbb{P}$; it
satisfies the countable chain condition, so it is in $\Gamma$. So
in the resulting model $V^{\mathbb{P}}$, there is a generic filter
$G \subseteq \mathbb{P}$ that will meet all dense sets in
$\mathbb{P}$, a fortiori meeting each $D \in \mathcal{D}$. So
$\phi(\mathbb{P},\mathcal{D})$ is $\Gamma$-forceable, since it is
true in $V^{\mathbb{P}}$. Moreover, working in $V^{\mathbb{P}}$,
let $\chi$ be the statement that $G$ is a filter in $\mathbb{P}$
which meets $D$.  $\chi$ is upwards absolute, so
$\phi(\mathbb{P},\mathcal{D})$ will remain true in any forcing
extension of $V^{\mathbb{P}}$. This proves the claim.
\end{proof}

So by $\MP_{\Gamma}(H(\omega_2))$ and the claim,
$\phi(\mathbb{P},\mathcal{D})$ is true (in $V$), for arbitrary
$\CCC$ forcing notion $\mathbb{P}$ of size $\omega_1$ and for any
family of dense sets $\mathcal{D}$ where $|\mathcal{D}| =
\omega_1$, so $MA_{\omega_1}$ is true.
\end{proof}

\begin{lem}\label{lem MA false}
Let $\Gamma$ be a class of forcing notions containing $\CCC$. Then
it is $\Gamma$-forceable that $MA_{\omega_1}$ is false.
\end{lem}

\begin{proof}
To find a $\Gamma$-forcing extension of $V$ in which
$MA_{\omega_1}$ is false it suffices to find one that contains a
Suslin tree, since $MA_{\omega_1}$ implies that there are none.
But forcing to add a single Cohen real will introduce a Suslin
tree (see \cite{jec:sth}, Theorem 28.12, or \cite{tf:amf}, chapter
III). This forcing notion is $\CCC$, hence in $\Gamma$.
\end{proof}

\begin{proof}[Proof of \thmref{thm-1}]
Assume $V\vDash \Box{\MP_{\Gamma}}(H(\omega_2))$, i.e. that
$\MP_{\Gamma}(H(\omega_2))$ is true and holds in any
$\Gamma$-forcing extension. Thus, by \lemref{lem MA true},
$MA_{\omega_1}$ is also true in any $\Gamma$-forcing extension.
But by \lemref{lem MA false}, one can always find a
$\Gamma$-forcing extension in which $MA_{\omega_1}$ is false, a
contradiction. So $\Box{\MP_{\Gamma}}(H(\omega_2))$ is false.
\end{proof}

\noindent\textit{Remark}.  This proof erases the possibility
mentioned in the comment after \thmref{thm_MP_ccc_kappa} that
$\MPmod{\CCC}$ directly implies $\MPmod{\CCC}(H(\kappa))$ for any
$\CCC$-absolutely definable cardinal $\kappa$.  This is because
$\MP_{\Gamma}(H(\omega_2))$ is false in the extension that adds
one Cohen real.  So ``$\MPmod{\CCC}$ implies
$\MPmod{\CCC}(H(\kappa))$" fails for $\kappa=\omega_1$.

\noindent\textit{Remark}.  Furthermore, this proof erases the
possibility mentioned in the remark after the proof of
\thmref{thm_MP_ccc_R} that $\MPmod{\CCC}$ directly implies
$\MPmod{\CCC}(\mathbb{R})$.  If $\MPmod{\CCC}$ holds, then it is
necessary.  So if $\MPmod{\CCC}$ implied
$\MPmod{\CCC}(\mathbb{R})$ then it would imply
$\Box\MPmod{\CCC}(\mathbb{R})$, which implies $\omega_1$ is
inaccessible to reals, but $\MPmod{\CCC}$ is just equiconsistent
with \ZFC.

\begin{cor}\label{cor-1}
The following are all false:
\begin{list}
 {(\arabic{listcounter})}{\usecounter{listcounter}\setlength{\rightmargin}{\leftmargin}}
\item
$\Box{\MPmod{\CCC}}(H(\omega_2))$
\item
$\Box{\MPmod{\CCC}}(H(2^\omega))$
\item
$\Box{\MPmod{\CARD }}(H(\omega_2))$
\item
$\Box{\MP_{\omega_1}}(H(\omega_2))$
\item
$\Box{\MPmod{\PROPER}}(H(\omega_2))$
\item
$\Box{\MPmod{\SEMIPROPER}}(H(\omega_2))$
\end{list}
\end{cor}

\begin{proof}
\textit{(1)}: See \thmref{thm-1}. \textit{(2)}: Under
$\MPmod{\CCC}$, the parameter set $H(2^\omega)$ includes
$H(\omega_2)$, so $\Box{\MPmod{\CCC}}(H(2^\omega))$ implies
$\Box{\MPmod{\CCC}}(H(\omega_2))$, which is false.
\textit{(3)-(6)}: The forcing classes of these principles all
contain $\CCC$, so \thmref{thm-1} applies.
\end{proof}

A question these results still do not answer is, what is the
consistency strength of $\Box{\MPmod{\CCC}}(\mathbb{R})$?

\section{$\omega_1$ is inaccessible to reals under $\Box \mathrm{MP_{CCC}}(\mathbb{R})$}

Restricting parameters to the reals, we now ask about the
consistency strength of $\Box{\MPmod{\CCC}}(\mathbb{R})$.  Hoping
that this restriction on parameters will not lead to
inconsistency, we can at least show that its consistency is
strictly beyond that of \ZFC.  A cardinal $\delta$ is said to be
\textbf{inaccessible to reals} if $L[r] \models ``\delta$ is
inaccessible" for any $r \in \mathbb{R}$, where $L[r]$ is the
model of \ZFC\ consisting of all sets constructible from $r$.  In
the case of $\delta = \omega_1$, the following standard
characterization of this property is useful (see Problem 32.4 in
\cite{jec:sth}, and Proposition 11.5 in \cite{kan:thi}.)

\begin{lem}\label{lem omega1 inaccess}
$\omega_1$ is inaccessible to reals if and only if $\omega_1 \neq
\omega_1^{L[z]}$ for any real $z$.
\end{lem}

\begin{proof}
$\Longrightarrow$: If $\omega_1 = \omega_1^{L[z]}$ for some real
$z$ then $L[z] \models \omega_1 = \omega^+$, therefore $L[z]
\models ``\omega_1$ is accessible" for the real $z$.

$\Longleftarrow$:  Let $\omega_1 \neq \omega_1^{L[z]}$ for any
real $z$.  Since $\omega_1$ is regular in $V$, it must be regular
in any submodel.  Looking at the submodel $L[z]$, this means it
will suffice to show that $\omega_1$ is a strong limit, i.e., that
$L[z]\models 2^\kappa < \omega_1$ for every $\kappa<\omega_1$.
Since $L[z]$ satisfies $GCH$, it is therefore sufficient to show
that, in $L[z]$, $\omega_1$ is not the successor of any cardinal
$\kappa$. Working in $L[z]$, suppose towards contradiction that
$\omega_1 = \kappa^+$. Then $\kappa < \omega_1$, so $\kappa$ is a
countable ordinal in $V$. So we can code the pair $\langle
\kappa,z \rangle$ with a real $w$. We move to $L[w]$, where $L[z]
\subseteq L[w]$ and, by definition, $L[w] \models |\kappa| =
\omega$.  So $\omega_1 = \omega_1^{L[w]}$ for the real $w$, a
contradiction.
\end{proof}

The next lemma constructs a $\CCC$-absolute sequence, which is
used in a subsequent theorem.

\begin{lem}\label{lem almost cohere}
There is a family of functions
$\{e_\alpha:\alpha\rightarrow\omega\mid\omega\leq\alpha<\omega_1\}$
with the property that, for all $\alpha < \omega_1$,

\begin{list}
 {(\arabic{bean})}{\usecounter{bean}\setlength{\rightmargin}{\leftmargin}}
 \item $e_\alpha$ is 1-1.
 \item For all $\beta$ such that $\omega\leq \beta < \alpha$, $|\{\xi<\beta|e_\alpha(\xi)\neq e_\beta(\xi)\}| <
 \omega$ (any disagreement between $e_\alpha$ and $e_\beta$ is finite).
\end{list}
\end{lem}

\noindent\textit{Remark.}  Such a sequence of functions is called
an \textbf{almost-coherent sequence}.

\begin{proof}
Fix a $\subseteq^*$- descending sequence $\{A_\alpha\}_{\alpha <
\omega_1}$ of subsets of $\omega$ (that is, such that for all
$\alpha < \beta < \omega_1$, $|A_\beta\diagdown A_\alpha| <
\omega$ and $|A_\alpha \diagdown A_\beta| = \omega$). Using this
sequence we will construct $\{e_\alpha\}$ inductively, by imposing
a third condition at each stage:
\begin{list}
 {\it{(\arabic{bean}})}{\usecounter{bean}\setlength{\rightmargin}{\leftmargin}\setcounter{bean}{2}}
 \item $range(e_\alpha)\cap A_\alpha
 = \emptyset$.
 \end{list}

Basis step ($\alpha = \omega$): Let $e_\omega$ be any 1-1 function
from $\omega$ into $\omega$ which avoids $A_\omega$.

Induction step: Let $\delta < \omega_1$. We assume that $e_\alpha$
is defined and satisfies (1)--(3) for all $\omega\leq\alpha <
\delta$.

First, suppose $\delta = \alpha+1$, a successor. Define
$e_{\alpha+1}(\xi)=e_\alpha(\xi)$ for all $\xi \in \alpha$, unless
$e_\alpha(\xi)\in A_{\alpha+1}$.  In the latter case, by (3) and
the induction hypothesis, $e_\alpha(\xi) \notin A_\alpha$, that
is, $e_\alpha(\xi) \in A_{\alpha+1}\diagdown A_\alpha$, a finite
set, so choose all such $e_{\alpha+1}(\xi)$ to have distinct
values in $A_\alpha\diagdown A_{\alpha+1}$, an infinite set.
(1)--(3) will thus be preserved.

Next, suppose $\delta$ is a (countable) limit ordinal: enumerate
the set of all ordinals below $\delta$ as an $\omega$ sequence.
There then will be a monotonically increasing subsequence of such
ordinals, $\{\alpha_n\}_{n<\omega}$, unbounded in $\delta$ and
beginning with $\alpha_0$, an arbitrary ordinal below $\delta$.
Our strategy will be to define a sequence of functions
$\{e_\delta^n:\alpha_n\rightarrow\omega\}_{n < \omega}$ whose
union will be our desired function $e_\delta$.

Claim: We can define such a sequence to have the following
properties, for all $n<\omega$:
\begin{list}
 {(\roman{bean})}{\usecounter{bean}\setlength{\rightmargin}{\leftmargin}}
 \item $e_\delta^{n+1}\upharpoonright \alpha_n = e_\delta^n$.
 (This ensures that $e_\delta=\cup_{i<\omega}e_\delta^i$ is
 a function.)
 \item $e_\delta^{n}=^*e_{\alpha_{n}}$.
 \item $range(e_\delta^{n})\cap A_\delta = \emptyset$.
 \item $e_\delta^n$ is 1-1.
 \item $domain(e_\delta^n) = \alpha_n$.
\end{list}

(We write $f=^*g$, of two functions when they disagree in only
finitely many places.)

From the claim we can establish the induction step of our original
argument, and the lemma will be proved:

(1) Suppose $e_\delta$ is not 1-1.  Then by cofinality of
$\{\alpha_n\}$ in $\delta$ and (v), there must be $n<\omega$ such
that $e_\delta^n$ is not 1-1, contradiction.

(2) For all $n<\omega$, $e_\delta\upharpoonright
\alpha_n=_{(i)}e_\delta^n =^*_{(ii)} e_{\alpha_n}$, which is what
(2) says in the case of $\alpha_n < \delta$.  Moreover, for any
ordinal $\beta<\delta$, there is $n<\omega$ such that
$\beta<\alpha_n<\delta$, since $\{\alpha_n\}$ has been defined to
be cofinal in $\delta$.  By the induction hypothesis, the
disagreement of $e_\beta$ with $e_{\alpha_n}$ is finite, and since
this is also the state of affairs between $e_{\alpha_n}$ and
$e_\delta$, so also must it be between $e_\beta$ and $e_\delta$.

(3) $range(e_\delta)\cap A_\delta =
range(\cup_{n<\omega}e_\delta^n)\cap A_\delta =
\cup_{n<\omega}(range(e_\delta^n)\cap A_\delta) = \emptyset$.

Proof of claim (by induction over $n<\omega$):

Basis step:  define $e_\delta^0:\alpha_0\rightarrow\omega$ by
emulating the successor case in the inductive definition of
$e_\delta$: For $\xi \in \alpha_0$, let
$e_\delta^0(\xi)=e_{\alpha_0}(\xi)$, unless $e_{\alpha_0}(\xi) \in
A_\delta$.  By (3), $e_{\alpha_0}(\xi) \notin A_{\alpha_0}$, so
$e_{\alpha_0}(\xi) \in A_\delta \diagdown A_{\alpha_0}$, a finite
set.  So let $e_\delta^0(\xi)$ take an unused value from
$A_{\alpha_0} \diagdown A_\delta$ (an infinite set) in this case.
Thus (ii)--(v) are satisfied.

Induction step: For all $\xi \in \alpha_n$, define
$e_\delta^{n+1}(\xi) = e_\delta^n(\xi)$ (this gives (i)).  For
$\xi \in [\alpha_n,\alpha_{n+1})$, proceed as in the basis step to
establish (ii)--(v) : Let
$e_\delta^{n+1}(\xi)=e_{\alpha_{n+1}}(\xi)$, unless
$e_{\alpha_{n+1}}(\xi) \in A_\delta$, and so forth.
\end{proof}

Note that the essential properties of the almost-coherent sequence
constructed are $\CCC$-absolute: the cardinal $\omega_1$ is
preserved, and the definition is $\Delta_0$.

We now content ourselves to quote a theorem from \cite{tf:amf},
after some introductory definitions.  Suppose $a:[\omega_1]^2
\rightarrow \omega$ is such that if one defines $a_\beta:\omega_1
\rightarrow \omega$ by $a_\beta(\alpha)=a(\alpha,\beta)$, then the
family $\{a_\beta\mid \beta<\omega_1\}$ is almost coherent. Set
$T(a) =
\{a_\beta\upharpoonright\alpha:\alpha\leq\beta<\omega_1\}$, the
set of all initial segments of the functions $a_\beta$. Then
$T(a)$ is an $\omega_1$-tree ordered by inclusion: the chains have
length $\omega$ since initial segments of injective functions from
$\omega_1$ to $\omega$ can only be extended for a countable chain.
And each node has finite branching, so the cardinality of each
level is countable.  Note that the almost coherent sequence of
\lemref{lem almost cohere} can be represented by a function
$e:[\omega_1]^2 \rightarrow \omega$.

Let $\mathcal{C}_\omega$ be the forcing notion consisting of
finite partial functions from $\omega$ to $\omega$. (This forcing
notion is $\CCC$). The union of a resulting
$\mathcal{C}_\omega$-generic filter is then a generic function
from $\omega$ to $\omega$, which we identify with the
corresponding ``Cohen" subset of $\omega$, or Cohen real.  If
$c:\omega\rightarrow\omega$ and $e:[\omega_1]^2 \rightarrow
\omega$, define $e_c = ce:[\omega_1]^2 \rightarrow \omega$.

\begin{thm}[Todorcevic]\label{thm cohere suslin}
If $c$ is a $\mathcal{C}_\omega$-generic function from $\omega$ to
$\omega$ and $e:[\omega_1]^2 \rightarrow \omega$ is almost
coherent, then $T(e_c)$ is a Suslin tree.
\end{thm}

The next theorem is the main result of this section.
\begin{thm}
$\Box \MPmod{\CCC}(\mathbb{R})$ implies that $\omega_1$ is
inaccessible to reals.
\end{thm}
\begin{proof}
Assume $\Box \MPmod{\CCC}(\mathbb{R})$.  By \lemref{lem omega1
inaccess} it is sufficient to show $\omega_1 \neq \omega_1^{L[z]}$
for any real $z$.  So suppose, towards contradiction, that
$\omega_1 = \omega_1^{L[z]}$ for some $z \in \mathbb{R}$. By
\lemref{lem almost cohere} $L[z]$ has an almost coherent family of
functions $\{c_\alpha:\alpha \rightarrow \omega|\alpha <
\omega_1\}$, which is definable in a way that is absolute to
forcing extensions. And since $\omega_1^{L[z]} = \omega_1$, the
family, in $V$, still maps all initial segments of $\omega_1$ to
$\omega$. And $L[z]$ has the correct size of the family.   Let $e
= \{e_\alpha|\alpha < \omega_1\}$ be the $L[z]$-least such family.
The point here is that this defines $e$ in any forcing extension.
So, in the $\CCC$-forcing extension that adds the Cohen real $c$,
the tree $T = T(ce)$, as provided by \thmref{thm cohere suslin},
is definable from $c$ and $z$. Moreover, the definition is
absolute to any extension.

Now consider the statement (with real parameters) $\phi(z,c)$ =
``The tree constructed in $L[z]$ with the real $c$, $T(ce)$, has
an $\omega_1$-branch". This formula has $z$ and $c$ as parameters
since the tree $T(ce)$ is definable from them. In $V[c]$,
$\phi(z,c)$ is false, since $T(ce)$ is Suslin by \thmref{thm
cohere suslin}. But by $\CCC$-forcing (using the tree itself) an
$\omega_1$-branch appears, so $\phi(z,c)$ is \CCC-necessary in
this extension. By the principle $\MPmod{\CCC}(\mathbb{R})$,
$\phi(z,c)$ must be true in $V[c]$, a contradiction.
\end{proof}

\chapter{$\mathrm{MP_{COHEN}}$ and variations}
Let $\delta$ be an uncountable ordinal.  Recall that
$Add(\omega,\delta)$ is the forcing notion, consisting of finite
partial functions from $\omega\times\delta=\delta$ to $2$, that
adds $\delta$ new reals. We define the class of forcing notions
\COHEN\ to be $\{Add(\omega,\theta)\mid \theta$ a cardinal$\}$. In
this chapter we establish equiconsistency of forms of the
principle $\MPmod{\COHEN}$ with \ZFC.  We will make use of the
fact that, if $\theta < \delta$ are cardinals, then forcing with
$Add(\omega,\theta)$ followed by forcing with $Add(\omega,\delta)$
is the same as forcing with $Add(\omega,\delta)$ alone.  In fact,
these forcing notions are absolute, since their conditions are all
finite, so all the names are in the ground model and this two-step
iteration is just product forcing. So
$Add(\omega,\delta)*Add(\omega,\alpha) \sim
Add(\omega,\alpha)*Add(\omega,\delta) \sim Add(\omega,\delta)$,
where $\sim$ denotes forcing-equivalence.

\begin{thm}\label{thm_model_MP_Cohen}
If there is a model of \ZFC\ $+$ $V_\delta\prec V$ then there is a
model of \ZFC\ $+$ $\MPmod{\COHEN}$.
\end{thm}

\begin{proof}
We find such a model by forcing with $Add(\omega,\delta)$.  Let
$G$ be $V$-generic over $Add(\omega,\delta)$.  Note that $G$ has
size $\delta$. We will show that $V[G]\models \MPmod{\COHEN}$.
Suppose $V[G]\models$ ``$\phi$ is \COHEN-forceably necessary" for
an arbitrarily chosen sentence $\phi$ in the language of \ZFC. Now
notice that $\phi$ is \COHEN-forceably necessary in $V$ as well as
$V[G]$---$\phi$ is forced necessary by
$Add(\omega,\delta)*Add(\omega,\alpha)\sim
Add(\omega,\delta+\alpha)$ for some $\alpha$.  Let $\alpha$ be the
least $\alpha$ such that $V[G]\models Add(\omega,\alpha)$ forces
that ``$\phi$ is \COHEN-necessary". We have $\alpha<\delta$,
whence $|\delta+\alpha|=\delta$, because $\alpha$ is definable and
hence in $V_\delta$ by elementarity. But
$Add(\omega,\delta)*Add(\omega,\alpha) \sim Add(\omega,\delta)$
forces ``$\phi$ is \COHEN-necessary". So $V[G]\models$ ``$\phi$ is
\COHEN-necessary".
\end{proof}

\begin{thm}\label{thm_model_MP_Cohen_R}
If there is a model of \ZFC\ $+$ $V_\delta\prec V$ $+$
$\cof(\delta)
> \omega$ then there is a model of \ZFC\ $+$ $\MPmod{\COHEN}(\mathbb{R})$.
\end{thm}

\begin{proof}
We use the same model as in \thmref{thm_model_MP_Cohen}, by
forcing with $Add(\omega,\delta)$. Let $G$ be $V$-generic over
$Add(\omega,\delta)$.  We will show that $V[G]\models
\MPmod{\COHEN}(\mathbb{R})$.  Suppose $V[G]\models$ ``$\phi(r)$ is
\COHEN-forceably necessary and $r \in \mathbb{R}$". Since
$\cof(\delta)>\omega$, $r$ is added to $V[G]$ at some ``stage",
i.e., by $Add(\omega,\theta)$ for $\theta < \delta$ (if not, and
unboundedly many $\theta<\delta$ are needed to decide $r$, then
$\delta$ has cofinality $\omega$, contradiction). So $r\in
V_\delta[r]\subseteq V[G\upharpoonright\theta]$ for $\theta <
\delta$, where $G\upharpoonright\theta= G\cap Add(\omega,\theta)$.
Since $V[G]\models$ ``$\phi(r)$ is \COHEN-forceably necessary and
$r \in \mathbb{R}$", let $\alpha$ be the least $\alpha$ such that
$V[G]\models Add(\omega,\alpha)$ forces that ``$\phi(r)$ is
\COHEN-necessary". We have  $\alpha<\delta$, because it is
definable from $r$, $\theta$, and $G\upharpoonright\theta$ and
hence in $V_\delta$ by elementarity. According to $V$,
$Add(\omega,\delta)*Add(\omega,\alpha)\Vdash \cong
Add(\omega,\delta)$ forces ``$\phi(r)$ is \COHEN-necessary".  So
$V[G]\models$ ``$\phi(r)$ is \COHEN-necessary".
\end{proof}

\noindent\textit{Remark}.  The use of $V_\delta \prec V$ in this
proof was not to circumvent the undefinability of truth, as in
earlier proofs, but rather to obtain the forcing of the
\COHEN-necessity of $\phi(x)$ by $Add(\omega,\theta)$ where
$\theta < \delta$.  This raises the possibility that these
arguments share a common thread involving reflection which is
deeper than what has been presented here.

\begin{thm}\label{thm_model_box_MP_Cohen_R}
If there is a model of \ZFC\ $+$ $V_\delta\prec V$ $+$
$\cof(\delta)
> \omega$ then there is a model of \ZFC\ $+$ $\Box
\MPmod{\COHEN}(\mathbb{R})$.
\end{thm}

\begin{proof}
Take the model $V[G]$ from \thmref{thm_model_MP_Cohen_R}, with $G$
$V$-generic over\linebreak $Add(\omega,\delta)$, where
$V_\delta\prec V$. In $V[G]$, force with $Add(\omega,\beta)$ for
some $\beta$. Let $H\subseteq Add(\omega,\beta)$ be $V[G]$-generic
(so new reals are added to $V[G]$).  Consider $V[G][H]$, and fix
some real $x$ in $V[G][H]$.  Then $x$ is decided by a countable
subset of $Add(\omega,\beta)$; call it $\mathbb{P}_0$. So $x$ is
in $V[G][H_0]$, where $H_0=H\cap \mathbb{P}_0\subseteq
\mathbb{P}_0$. Since $\mathbb{P}_0$ adds the real $x$ and is
countable, $\mathbb{P}_0 \cong Add(\omega,1)$. This gives
$V[G][H_0] = V[G*H_0]$, where $G*H_0$ is $V$-generic over
$Add(\omega,\delta)*Add(\omega,1)\cong Add(\omega,\delta)$. So
$V[G*H_0]$ satisfies the conditions of
\thmref{thm_model_MP_Cohen_R}, so $V[G*H_0]\models
\MPmod{\COHEN}(\mathbb{R}^{V[G*H_0]}))$.  But $x$ is in
$\mathbb{R}^{V[G*H_0]}$, so adding it did not invalidate
$\MPmod{\COHEN}(\mathbb{R})$.  So $V[G]\models \Box
\MPmod{\COHEN}(\mathbb{R})$.
\end{proof}

\begin{cor}\label{cor con MP Cohen R}
The following are equivalent:
\begin{list}
{(\arabic{listcounter})}{\usecounter{listcounter}\setlength{\rightmargin}{\leftmargin}}
\item
$Con(\ZFC)$
\item
$Con(\ZFC$ $+$ $\MPmod{\COHEN})$
\item
$Con(\ZFC$ $+$ $\MPmod{\COHEN}(\mathbb{R}))$
\item
$Con(\ZFC$ $+$ $\Box\MPmod{\COHEN}(\mathbb{R}))$.
\end{list}
\end{cor}

\begin{proof}
$(4)\Longrightarrow (3)\Longrightarrow (2)\Longrightarrow (1)$:
Obvious.

$(1)\Longrightarrow (4)$: By \lemref{lem_Vdelta_submod_V_big_cof},
we know that it is equiconsistent with \ZFC\ that the ground model
$V$ satisfies $V_\delta\prec V$ for some cardinal $\delta$, with
$\cof(\delta) > \omega$.  And by
\thmref{thm_model_box_MP_Cohen_R}, we know, given such a model of
\ZFC $+$ $V_\delta\prec V$ $+$ $\cof(\delta) > \omega$, that there
is a forcing extension model of \ZFC\ $+$
$\Box\MPmod{\COHEN}(\mathbb{R})$.
\end{proof}

These equiconsistency results immediately raise questions as to
whether any of them are, in fact, equivalences.\vspace {12pt}

\noindent\textbf{Question 3.} Is $\MPmod{\COHEN}$ equivalent to
$\MPmod{\COHEN}(\mathbb{R})$?\vspace {12pt}

\noindent\textbf{Question 4.} Is $\MPmod{\COHEN}(\mathbb{R})$
equivalent to $\Box\MPmod{\COHEN}(\mathbb{R})$?\vspace {12pt}

\begin{thm}\label{thm model MP Cohen cof}
If there is a model of \ZFC\ $+$ $V_\delta\prec V$, then there is
a \COHEN-forcing extension which models  \ZFC\ $+$
$\MPmod{\COHEN}(H(\cof(\delta)))$.
\end{thm}

\begin{proof}
The proof is almost identical to \thmref{thm_model_MP_Cohen_R}.
Again, we force with $Add(\omega,\delta)$.  The single difference
is that the parameter of the formula $\phi$ is taken from
$H(\cof(\delta))$. Since it still has size less than
$\cof(\delta)$, the parameter appears from forcing with
$Add(\omega,\theta)$ for some $\theta < \delta$, and the proof
proceeds in the same way.
\end{proof}

\noindent\textit{Remark}.  Since the final model is $V[G]$ where
$G$ is $V$-generic over $Add(\omega,\delta)$, it is clear that
$2^\omega = \delta$ in $V[G]$.  As long as $\delta$ is larger than
the continuum in the ground model, the model created is a model of
$\MPmod{\COHEN}(H(\cof(2^\omega)))$.

\begin{cor}\label{cor con MP Cohen 2^omega}
The following are equivalent:
\begin{list}
{(\arabic{listcounter})}{\usecounter{listcounter}\setlength{\rightmargin}{\leftmargin}}
\item
$Con(\ZFC)$
\item
$Con(\ZFC$ $+$ $\MPmod{\COHEN}(H(\cof(2^\omega))))$
\item
$Con(\ZFC$ $+$
$\MPmod{\COHEN}(H(\omega_{17}))+2^\omega=\omega_{17})$
\item
$Con(\ZFC$ $+$ $\Box \MPmod{\COHEN}(H(\cof(2^\omega))))$.
\item
$Con(\ZFC$ $+$ $\Box
\MPmod{\COHEN}(H(\omega_{17}))+2^\omega=\omega_{17})$.
\end{list}
\end{cor}

\begin{proof}
\lemref{lem_Vdelta_submod_V_big_cof} can be used in the proof of
\thmref{thm_model_MP_Cohen_R} to allow the choice of $\delta$ to
have cofinality larger than any particular definable cardinal.
And by modifying the proof of \thmref{thm_model_box_MP_Cohen_R} in
the obvious way, we get the necessary versions of these
corollaries.
\end{proof}

\begin{lem}\label{lem_MP cohen_2_omega 2_kappa}
If $\MPmod{\COHEN}(2^\omega)$ then $2^\omega = 2^\kappa$ for all
infinite $\kappa < 2^\omega$.
\end{lem}

\begin{proof}
Let $V\models\MPmod{\COHEN}(2^\omega)$.  Suppose
$\omega<\kappa<2^\omega$.  Let $G$ be $V$-generic over
$\mathbb{P}=Add(\omega,\lambda)$ where $\lambda = (2^\kappa)^V$.
We compute the size of $2^\kappa$ in the extension $V[G]$, i.e.,
how many subsets of $\kappa$ have been added by forcing with
$\mathbb{P}$. We can bound $(2^\kappa)^{V[G]}$ above by counting
the nice $\mathbb{P}$-names for subsets of $\kappa$. $\mathbb{P}$
is in \CCC, and $|\mathbb{P}|=\lambda=(2^\kappa)^V$, so there are
$\lambda^\omega$ antichains in $\mathbb{P}$.  So there are
$(\lambda^\omega)^\kappa =
(2^\kappa)^{\omega\cdot\kappa}=2^\kappa$ nice $\mathbb{P}$-names
in $V$.  So the number of subsets of $\kappa$ does not increase
when forcing with $\mathbb{P}$, i.e.,
$(2^\kappa)^{V[G]}=(2^\kappa)^{V}$.  But
$\mathbb{P}=Add(\omega,\lambda)$ has changed $2^\omega$ to
$\lambda = (2^\kappa)^V = (2^\kappa)^{V[G]}$.  So $V[G]\models
2^\omega = 2^\kappa$.  Moreover, further \COHEN-forcing cannot
destroy this statement:  Since $\omega < \kappa$, $2^\omega \leq
2^\kappa$ is absolute.  And $2^\omega < 2^\kappa$ is impossible in
any further \COHEN-forcing extension, by a similar nice names
argument: Let $\mathbb{Q}=Add(\omega, \mu)$, where $V[G]\models
\mu > 2^\omega$, and let $H$ be $V[G]$-generic over $\mathbb{Q}$.
Working in $V[G]$, the number of nice $\mathbb{Q}$-names of
subsets of $\kappa$ is $(\mu^\omega)^\kappa =\mu^\kappa=\mu$
(since $\omega<\kappa<\mu$). So $V[G][H]\models 2^\kappa = \mu =
2^\omega$.  This shows that $2^\omega = 2^\kappa$ is
\COHEN-forceably necessary.  Defining $\phi(\kappa)$ to be
``$2^\kappa = 2^\omega$", $\kappa$ is a parameter allowed by
$\MPmod{\COHEN}(2^\omega)$.  Hence, by applying
$\MPmod{\COHEN}(2^\omega)$, $2^\kappa = 2^\omega$ is true.
\end{proof}

\begin{lem}\label{lem_MP_cohen_implies Vdelta_submod_V_delta_wi}
If there is a model of \ZFC\ $+$ $\MPmod{\COHEN}(H(2^\omega))$,
then it has an inner model of \ZFC\ $+$ $V_\delta\prec V$ $+$
``$\delta$ is inaccessible".
\end{lem}

\begin{proof}
We first show $\MPmod{\COHEN}(H(2^\omega))$ implies that
$2^\omega$ is weakly inaccessible.  By \lemref{lem_MP
cohen_2_omega 2_kappa} $2^\omega$ is regular (for all
$\kappa<2^\omega$, $\cof(2^\omega)=\cof(2^\kappa)>\kappa$).  And
$2^\omega$ cannot be a successor cardinal: for any cardinal
$\kappa$, $Add(\omega,\kappa^{++})$ forces that it is
\COHEN-necessary that $2^\omega > \kappa^+$, so that must be true.
So $2^\omega$ is a regular limit cardinal.  Since it is weakly
inaccessible, $L$ models that $2^\omega$ is (strongly)
inaccessible.

Next we show that, if $\delta = 2^\omega$, $L_\delta \prec L$.
Suppose $L\models \exists y \psi(a,y)$ for the formula $\exists y
\psi(a,y)$ with parameter $a$. Let $\alpha$ be the least cardinal
such that $\exists y \in L_\alpha$ such that $\psi(a,y)$. Consider
$\phi(a)=$ ``the least $\alpha$ such that there is a $y$ in
$L_\alpha$ with $\psi(a,y)^L$ is less than $2^\omega$".  This is
expressed using the parameter $a$ in $L_\delta \subseteq
H(\delta)$.  Since this is \COHEN-forceably necessary, it is true.
So there is a $y$ in $L_\alpha \subseteq L_\delta$ such that
$\psi(a,y)$.  So by the Tarski-Vaught criterion, $L_\delta \prec
L$, and $L$ is the desired inner model.
\end{proof}

\begin{lem}\label{lem_ZFC_Vdelta_submod_V_delta _wi_MPcohen}
If there is a model of \ZFC\ $+$ $V_\delta\prec V$ $+$ ``$\delta$
is inaccessible" then there is a model of \ZFC\ $+$
$\MPmod{\COHEN}(H(2^\omega))$.
\end{lem}

\begin{proof}
Suppose $V\models\ZFC$ $+$ $V_\delta\prec V$ $+$ ``$\delta$ is
inaccessible".  Let $G$ be $V$-generic over
$\mathbb{P}=Add(\omega,\delta)$.  We claim that $V[G]\models
\MPmod{\COHEN}(H(\delta))$.  Suppose $V[G]\models$ ``$\phi(x)$ is
\COHEN-forceably necessary and $x$ is in $H(\delta)$".  It will
suffice to show $V[G]\models$ ``$\phi(x)$ is \COHEN-necessary".
Since $x$ has size less than $\delta$, which is regular, it must
be added by $Add(\omega,\theta)$ for some $\theta<\delta$.

Just as in the proof of \thmref{thm_model_MP_Cohen_R}, let
$\alpha$ be the least $\alpha$ such that $V[G]\models
``Add(\omega,\alpha)$ forces that $\phi(x)$ is \COHEN-necessary".
We have $\alpha<\delta$, because it is definable from $x$,
$\theta$, and $G\upharpoonright\theta$ and hence in $V_\delta$ by
elementarity. According to $V$,
$Add(\omega,\delta)*Add(\omega,\alpha)\cong
Add(\omega,\delta)\Vdash$ ``$\phi(x)$ is \COHEN-necessary". So
$V[G]\models$ ``$\phi(x)$ is \COHEN-necessary".
\end{proof}

\begin{thm}
The following are equivalent:
\begin{list}
{(\arabic{listcounter})}{\usecounter{listcounter}\setlength{\rightmargin}{\leftmargin}}
\item
$Con(\ZFC$ $+$ $\MPmod{\COHEN}(H(2^\omega)))$
\item
$Con(\ZFC$ $+$ $V_\delta \prec V$ $+$ ``$\delta$ is an
inaccessible cardinal"$)$
\item
$Con(\ZFC$ $+$ ``$ORD$ is Mahlo"$)$.
\end{list}
\end{thm}

\begin{proof}
The equivalence of $(1)$ with $(2)$ is from \lemref{lem_MP
cohen_2_omega 2_kappa} and \lemref{lem_MP_cohen_implies
Vdelta_submod_V_delta_wi}.  The proof of the equivalence of $(2)$
with $(3)$ is in \cite{ham:max}.
\end{proof}

So far a strong analogy seems to exists between the results
concerning consistency strengths of maximality principles for the
forcing class \COHEN\ and the class \CCC. We pursue this
correspondence further into the boxed (necessary) versions of
these principles.

\begin{thm}\label{thm_box_MPcohenH2^omega}
$\Box\MPmod{\COHEN}(H(2^\omega))$ is false.
\end{thm}

\begin{proof}
By \lemref{lem_MP cohen_2_omega 2_kappa},
$\MPmod{\COHEN}(H(2^\omega))$ implies that $2^\omega$ is weakly
inaccessible.  So $\Box\MPmod{\COHEN}(H(2^\omega))$ implies that
$2^\omega$ is weakly inaccessible in every \COHEN-forcing
extension.  But there are \COHEN-forcing extensions that add any
desired cardinality of reals, including successor cardinalities,
which are therefore not weakly inaccessible.
\end{proof}

\chapter{$\mathrm{MP_{COLL}}$ and variations}
Let $\delta$ be an uncountable ordinal.  Recall that
$Col(\omega,\theta)$ is the forcing notion, consisting of finite
partial injective functions from $\theta$ to $\omega$, that
collapses the cardinal $\delta$ to $\omega$.  We will use the
notation $Col(\omega,<\delta)$ to represent the L\'evy collapse of
$\delta$ to $\omega_1$, for limit cardinal $\delta$. We define the
class of forcing notions \COLL\ to be $\{Col(\omega,\theta)\mid
\theta$ a cardinal$\}\cup \{Col(\omega,<\theta)\mid \theta$ a
limit cardinal$\}$, and we include, for technical reasons, all
forcing notions that are forcing-equivalent to some
$Col(\omega,\theta)$ or $Col(\omega,<\theta)$. Two forcing notions
are forcing-equivalent if they produce the same forcing
extensions, i.e., they have isomorphic regular open Boolean
Algebras. This does not alter the meaning of $\MPmod{\COLL}$ or
any variation thereof, since saying ``$M$ is a \COLL-forcing
extension" means the same thing in either interpretation of \COLL.
One sees it is natural to include the L\'evy collapse with forcing
notions that collapse any $\theta$ to $\omega$ by recognizing that
$Col(\omega,\theta)$ is equivalent to $Col(\omega, <\theta^+)$.
That is, $Col(\omega,\theta)$ collapses the entire interval
$(\omega,\theta]$, while $Col(\omega,<\theta)$ collapses the
interval $(\omega,\theta)$.  So there are natural closure
properties for this class.

In this chapter we establish equiconsistency of forms of the
principle $\MPmod{\COLL}$ with \ZFC.  We will make use of the fact
that, if $\theta < \delta$ are cardinals, with $\delta$ a limit
cardinal, then forcing with $Col(\omega,<\delta)$ followed by
forcing with $Col(\omega,\theta)$ is the same as forcing with
$Col(\omega,<\delta)$ alone.  This chapter retraces most of the
same kinds of arguments as found in the previous chapter.

Recall the definition of the absorption relation between two
classes of forcing notions: $\Gamma_2$ \textbf{absorbs} $\Gamma_1$
if it is true and $\Gamma_1$-necessary that $\Gamma_2$ is
contained in $\Gamma_1$ and for any $\mathbb{P}$ in $\Gamma_1$
there exists $\mathbb{Q}$ such that $V^\mathbb{P}\models$
``$\dot{\mathbb{Q}}$ is in $\Gamma_2$" and
``$\mathbb{P}*\dot{\mathbb{Q}}$ is in $\Gamma_2$". There is a
well-known fact that can be cast in this language.
\begin{lem}\label{lem_COLL_absorbs_all_forcing}
The class \COLL\ absorbs all forcing.  That is, if $\mathbb{P}$ is
any forcing notion, then there is $\mathbb{Q}$ in \COLL\ (
according to $V^{\mathbb{P}}$) such that
$\mathbb{P}*\dot{\mathbb{Q}}$ is in \COLL.
\end{lem}

To prove this, one makes use of a lemma that Solovay used in his
famous Lebesgue measurability result.  It can be found as
Proposition 10.20 in \cite{kan:thi}.

\begin{lem}[Solovay]\label{lem_P_dense_embed_coll}
Let $\mathbb{P}$ be a separative forcing notion of size $\leq
\alpha$ which collapses $\alpha$ to $\omega$.  Then there is a
dense subset of $Col(\omega,\alpha)$ that densely embeds into
$\mathbb{P}$.
\end{lem}

\begin{proof}[Proof of \lemref{lem_COLL_absorbs_all_forcing}]
Let $\mathbb{P}$ be any forcing notion.  Pick cardinal $\alpha$ to
be at least $|\mathbb{P}|$, and let $\mathbb{Q} =
Col(\omega,\alpha)$. Then $|\mathbb{P}*\dot{\mathbb{Q}}|= \alpha$,
and $\mathbb{P}*\dot{\mathbb{Q}}$ collapses $\alpha$ to $\omega$.
So by \lemref{lem_P_dense_embed_coll} there is a dense embedding
from a dense subset of $\mathbb{Q}=Col(\omega,\alpha)$ into
$\mathbb{P}*\dot{\mathbb{Q}}$.  So, as forcing notions,
$\mathbb{P}*\dot{\mathbb{Q}}$ is equivalent to
$Col(\omega,\alpha)$, which is in \COLL.  So
$\mathbb{P}*\dot{\mathbb{Q}}$ itself is in \COLL.
\end{proof}

\begin{cor}
$\MPmod{\COLL}$ implies $\MP$.
\end{cor}

\begin{proof}
By \lemref{lem_COLL_absorbs_all_forcing} and
\corref{cor_absorbs_MP_Gamma}.
\end{proof}

\begin{cor}
$\MPmod{\COLL}(\mathbb{R})$ implies $\MP(\mathbb{R})$.
\end{cor}

\begin{proof}
By \lemref{lem_COLL_absorbs_all_forcing} and
\corref{cor_absorbs_MP_Gamma_X}.
\end{proof}

\begin{thm}\label{thm_model_MP_Coll}
If there is a model of \ZFC\ $+$ $V_\delta\prec V$ then there is a
model of \ZFC\ $+$ $\MPmod{\COLL}$.
\end{thm}

\begin{proof}
We find such a model by forcing with $Col(\omega,<\delta)$.  Let
$G$ be $V$-generic over $Col(\omega,<\delta)$.  Note that $G$ has
size $\delta$. We will show that $V[G]\models \MPmod{\COLL}$.
Suppose $V[G]\models$ ``$\phi$ is \COLL-forceably necessary" for
an arbitrarily chosen sentence $\phi$ in the language of \ZFC.  So
$V[G]\models$ ``there is some member of \COLL\ that forces that
$\phi$ is \COLL-necessary".  So some finite  condition in
$Col(\omega,<\delta)$ decides that some $Col(\omega,\lambda)$
forces that $\phi$ is \COLL-necessary.  So there is $\theta
<\delta$ such that $V[G\upharpoonright \theta]\models$ ``$\phi$ is
\COLL-forceably necessary".  Now notice that $\phi$ is
\COLL-forceably necessary in $V$ as well as $V[G]$---$\phi$ is
forced necessary by $Col(\omega,<\delta)*Col(\omega,\alpha)$ for
some $\alpha$.  Let $\alpha$ be the least $\alpha$ such that
$V[G\upharpoonright \theta]\models Col(\omega,\alpha)$ forces that
``$\phi$ is \COLL-necessary". We have $\alpha<\delta$ because
$\alpha$ is definable from $G\upharpoonright \theta$ and hence in
$V_\delta[G\upharpoonright \theta]$ (and hence in $V_\delta[G]$)
by elementarity. But $Col(\omega,<\delta)*Col(\omega,\alpha) \cong
Col(\omega,<\delta)$ forces ``$\phi$ is \COLL-necessary". So
$V[G]\models$ ``$\phi$ is \COLL-necessary".
\end{proof}

\begin{cor}\label{cor con MP Coll}
The following are equivalent:
\begin{list}
{(\arabic{listcounter})}{\usecounter{listcounter}\setlength{\rightmargin}{\leftmargin}}
\item
$Con(\ZFC)$
\item
$Con(\ZFC$ $+$ $\MPmod{\COLL})$
\end{list}
\end{cor}

\begin{thm}\label{thm_model_MP_Coll_R}
If there is a model of \ZFC\ $+$ $V_\delta\prec V$ $+$ ``$\delta$
is inaccessible"  then there is a model of \ZFC\ $+$
$\MPmod{\COLL}(\mathbb{R})$.
\end{thm}

\begin{proof}
Let $V$ satisfy \ZFC\ $+$ $V_\delta\prec V$ $+$ ``$\delta$ is
inaccessible".  Now force with $Col(\omega,<\delta)$. Let $G$ be
$V$-generic over $Col(\omega,<\delta)$.  We will show that
$V[G]\models \MPmod{\COLL}(\mathbb{R})$.  Suppose $V[G]\models$
``$\phi(r)$ is \COLL-forceably necessary and $r \in \mathbb{R}$".
Since $\delta$ is inaccessible, $Col(\omega,<\delta)$ has the
$\delta$-cc, so $r$ is added to $V[G]$ ``at some earlier stage",
i.e., by $Col(\omega,\theta)$ for $\theta < \delta$. So $r\in
V_\delta[r]\subseteq V_\delta[G\upharpoonright\theta]$ for $\theta
< \delta$, where $G\upharpoonright\theta= G\cap
Col(\omega,\theta)$. Since $|G\upharpoonright\theta|<\delta$,
forcing with $G\upharpoonright\theta$ is ``small" forcing, and so,
by \lemref{lem_Vdelta_submod_V_forcing_ext},
$V_\delta[G\upharpoonright\theta]\prec V[G\upharpoonright\theta]$.
Working in $V[G\upharpoonright\theta]$, since
$V[G\upharpoonright\theta]\models$ ``$\phi(r)$ is \COLL-forceably
necessary and $r \in \mathbb{R}$", let $\alpha$ be the least
$\alpha$ such that $V[G\upharpoonright\theta]\models
Col(\omega,\alpha)$ forces that ``$\phi(r)$ is \COLL-necessary".
We have $\alpha<\delta$, because it is definable from $r$,
$\theta$, and $G\upharpoonright\theta$ and hence in
$V_\delta[G\upharpoonright\theta]$ by elementarity. According to
$V$, $Col(\omega,<\delta)*Col(\omega,\alpha) \cong
Col(\omega,<\delta)$ forces ``$\phi(r)$ is \COLL-necessary".  So
$V[G]\models$ ``$\phi(r)$ is \COLL-necessary".
\end{proof}

\noindent\textit{Remark}.  This is another alternate use of
$V_\delta \prec V$, not to circumvent the undefinability of truth,
 but to obtain the forcing of the \COLL-necessity of
$\phi(x)$ by $Col(\omega,\theta)$ where $\theta < \delta$.

\begin{thm}\label{thm_model_MP_Coll_R_L_omega_1}
If there is a model of \ZFC\ $+$ $\MPmod{\COLL}(\mathbb{R})$ then
it has an inner model of \ZFC\ $+$ $V_\delta\prec V$ $+$
``$\delta$ is inaccessible".  Specifically, $\omega_1$ is
inaccessible in $L$ and $L_{\omega_1}\prec L$.
\end{thm}

\begin{proof}
This proof just mimics the proof in \cite{ham:max} of the
analogous result for $\MP(\mathbb{R})$.  First, we see that
$\MPmod{\COLL}(\mathbb{R})$ implies that $\omega_1$ is
inaccessible to reals:  For any real $x$, let $\phi(x)$ be ``the
$\omega_1$ of $L[x]$ is countable", which implies $\omega_1^{L[x]}
< \omega_1$, which implies $\delta=\omega_1$ is inaccessible in
$L$. But this is clearly \COLL-forceably necessary, by
$Col(\omega,\omega_1)$, so it is true.  Further, we get
$L_{\omega_1}\prec L$, by the Tarski-Vaught argument.  Suppose
$L\models\exists y \psi(a,y)$ for $a$ in $L_{\omega_1}$.  Let
$\alpha$ be the least such that there is a $y$ in $L_\alpha$ such
that $\psi(a,y)$.  Consider formula $\phi$ which is ``the least
$\alpha$ such that there is a $y$ in $L_\alpha$ such that
$\psi(a,y)$ is countable".  This is \COLL-forceably necessary,
using $Col(\omega, \alpha)$, hence true by
$\MPmod{\COLL}(\mathbb{R})$.
\end{proof}

\begin{cor}\label{cor con MP Coll R}
The following are equivalent:
\begin{list}
{(\arabic{listcounter})}{\usecounter{listcounter}\setlength{\rightmargin}{\leftmargin}}
\item
$Con(\ZFC$ $+$ $\MPmod{\COLL}(\mathbb{R}))$
\item
$Con(\ZFC$ $+$ $V_\delta\prec V$ $+$ ``$\delta$ is
inaccessible"$)$
\item
$Con(\ZFC$ $+$ ``$ORD$ is Mahlo"$)$.
\end{list}
\end{cor}

\begin{proof}
$(2)\Longleftrightarrow (3)$: Proven in the analogous result for
$\MP(\mathbb{R})$ in \cite{ham:max}.

$(2)\Longrightarrow (1)$: By \thmref{thm_model_MP_Coll_R}.

$(1)\Longrightarrow (2)$: By
\thmref{thm_model_MP_Coll_R_L_omega_1}.
\end{proof}

\chapter{Large cardinals and indestructibility}
\section{${\mathrm{MP}_{\Gamma(\kappa)}}$}\label{sec lc indestruct}
Certain classes of forcing notions are defined from a particular
cardinal $\kappa$.  Examples are $\kappa-cc$, $\kappa$-closed,
$<\kappa$-directed closed, and so on.  In this chapter we'll
explore maximality principles based on such classes, striving to
obtain equiconsistency results by applying the same tools used
earlier.  Such classes will be collectively denoted by
$\Gamma(\kappa)$, to show the role of $\kappa$ in the definition
of the class.  The related maximality principles will then be
written $\MP_{\Gamma(\kappa)}$.

Let the formula $\sigma(x)$ express some cardinal property. Then,
if there is a class $\Gamma$ of forcing notions under which
$\sigma(\kappa)$ is $\Gamma$-necessary, $\kappa$ is said to be
\textbf{indestructible} under forcing by $\Gamma$, or
\textbf{$\Gamma$-indestructible}.  In modal notation, this is
expressed by $\Box_\Gamma\sigma(\kappa)$.  The creation via
forcing of a model of \ZFC\ witnessing this indestructibility is
called a \textbf{preparation}.  A typical indestructibility result
says that it is forceable that $\sigma(\kappa)$ is
$\Gamma$-necessary (in modal notation, this is expressed as
$\Diamond\Box_{\Gamma}\sigma(\kappa)$.  Note that $\Diamond$ has
no subscript, indicating that the preparation itself need not be
in the class $\Gamma$.

The literature contains a growing body of indestructibility
results, beginning with Laver's preparation of a model of the
indestructibility of a supercompact cardinal $\kappa$ under
$<\kappa$-directed-closed forcing.  Such a $\kappa$ has come to be
called \textbf{Laver-indestructible}.  Due to its familiarity, it
may be easiest to consider $\Gamma(\kappa)$ as being the class of
$<\kappa$-directed closed forcing notions (abbreviated $<\kappa
-dc$) and $\sigma(\kappa)$ as saying that $\kappa$ is
supercompact. Laver's well-known result will then be expressed as
$\Diamond\Box_{\Gamma(\kappa)}\sigma(\kappa)$.  However, various
other interpretations will also be possible.  In fact, we first
prove an equiconsistency result for $\MP_\Gamma$ in schematic form
that will have, as instances, corollaries corresponding to the
various indestructibility results known so far, as well as those
yet to be discovered.  In each of these cases the class $\Gamma$
of forcing notions inevitably depends on the parameter $\kappa$,
so this will be built into the scheme. So let $\Gamma(\kappa)$
denote a class of forcing notions that is defined in terms of the
parameter $\kappa$, a cardinal. If $\kappa$ is clear from the
context of a statement, we may write $\Gamma$ in place of
$\Gamma(\kappa)$, leaving $\kappa$ as an implicit parameter.

To prove this theorem, we first give a special case of
\lemref{lem_Vdelta_submod_V}.

\begin{lem}\label{lem_Vdelta_submod_V_indestruct_kappa}
Let $\Gamma$ be a class of forcing notions, parametrized by
$\kappa$.  Let $\sigma$ be any definable unary predicate.  If
there is a model $M$ of $\ZFC+ \square_{\Gamma}\sigma(\kappa)$
then there is a model of
$\ZFC+\square_{\Gamma}\sigma(\kappa)+V_\delta\prec V$.
\end{lem}

\begin{proof}
Use \lemref{lem_Vdelta_submod_V}, where theory $T$ is
$\ZFC+\square_{\Gamma}\sigma(\kappa)$.
\end{proof}

\noindent\textit{Remark.}  Again we are working, not merely in the
language or model of \ZFC, but in an expansion which has the
constant $\kappa$.  So if it is consistent that $\kappa$ is
indestructible by forcing with $\Gamma$, then it is also
consistent that this holds and $V_\delta\prec V$ as well.  One can
assume that $\kappa < \delta$.  This follows from the proof of
\lemref{lem_Vdelta_submod_V_big_cof}.  In fact, if $\kappa \geq
\delta$, one could find a suitable $\kappa'<\delta$ by
elementarity.

\noindent\textit{Remark.} We can use the name $\kappa$ as a
constant in the language of both of the above equiconsistent
theories to denote the same object in both. This is because, due
to \lemref{lem_Vdelta_submod_V_in_elem_ext}, a model of
$V_\delta\prec V$ can be found by taking an elementary extension
of any model of the theory on the left to obtain a model of the
theory on the right. This is analogous with large cardinal results
in which a forcing extension can provide new desired properties
(such as indestructibility) for $\kappa$ as interpreted by both
the ground model and the extension.  In the language of modal
logic, $\kappa$ becomes a rigid identifier in this context of
different models connected by elementarity.  This idea was
described earlier, but only in the context of ground models as
elementary substructures of forcing extensions.

\noindent\textit{Remark.} Note that if $\Gamma$ is S4 then
$\square_{\Gamma}\sigma(\kappa)$ directly implies
$\sigma(\kappa)$, since the class $\Gamma$ includes trivial
forcing.

The next theorem draws on the argument in the proofs of Lemma 2.6
and Theorem 5.1 in \cite{ham:max}.

\beginThm\label{thm_con_MP_Gamma_indestruct_kappa}
Let $\Gamma$ be a class of forcing notions, parametrized by
$\kappa$, which is $\Gamma$-necessarily closed under iterations of
length $\omega$ using appropriate support, and which
$\Gamma$-necessarily contains the trivial forcing. Further, let
$\sigma$ be a unary predicate.  Then

$Con(\ZFC+\square_{\Gamma}\sigma(\kappa))$ if and only if
$Con(\ZFC+\square_{\Gamma}\sigma(\kappa)+\MP_{\Gamma})$.
\end{thm}

\noindent\textit{Remark}.  Thus, if we make $\sigma(\kappa)$
indestructible by forcing in $\Gamma$, we have, at the same time,
a model of $\MP_\Gamma$.

\begin{proof}
The implication to the left is trivial.  In the other direction,
suppose $V\models \ZFC+\square_{\Gamma}\sigma(\kappa)$. By
\lemref{lem_Vdelta_submod_V_indestruct_kappa} we can also assume
$V\models V_\delta\prec V$. Let $\{\phi_n\}_{n \in \omega}$
enumerate all sentences in the language of \ZFC.  Now recursively
define an $\omega$-stage forcing iteration $\mathbb{P}$ such that
at each stage $n \in \omega$, $\mathbb{P}_n \in V_\delta$.
Specifically, define $\mathbb{P}_0 = \{\emptyset \}$, the notion
of trivial forcing. Next suppose $\mathbb{P}_n$ has been defined,
and in the model $V_\delta^{\mathbb{P}_n}$  define a
$\mathbb{P}_n$-name for a forcing notion, $\mathbb{Q}_n$, as
follows. If $\phi_n$ is $\Gamma$-forceably necessary then choose
$\dot{\mathbb{Q}}_n$ to be the $\mathbb{P}_n$-name of a notion of
forcing in $\Gamma^{V_\delta^{\mathbb{P}_n}}$ that forces
$\square_\Gamma\phi_n$ \footnote{This is exactly the point at
which we need to be working in $V_\delta$ rather than $V$. Without
a truth predicate---needed for a definition of the forcing
relation---there is no formula applicable to uniformly express
that $\phi_n$ is forceably necessary over $V$, for all $n$.}.
Otherwise let $\dot{\mathbb{Q}}_n$ be the name for trivial
forcing.  Let $\mathbb{P}_{n+1} =
\mathbb{P}_n*\dot{\mathbb{Q}}_n$.  Finally, let $\mathbb{P}$ be
the $\omega$-iteration of $\{\mathbb{P}_n\}_{n \in \omega}$.  Use
appropriate support, as required by the closure conditions of the
family $\Gamma$.

Let $G \subset \mathbb{P}$ be $V$-generic. I now claim that $V[G]
\models \MP_\Gamma$.  To see this, let $\phi$ be a sentence in the
language of \ZFC\ which is $\Gamma$-forceably necessary in $V[G]$.
I will show that $V[G] \models \square_\Gamma\phi$. First, $\phi$
must be $\phi_n$ for some $n \in \omega$.  Factor
$\mathbb{P}=\mathbb{P}_n*\mathbb{P}_{TAIL}$. This gives $V[G] =
V[G_n][G_{TAIL}]$, taking respective generic filters. Since $V[G]$
is a $\Gamma$-forcing extension of $V[G_n]$, $\phi$ is
$\Gamma$-forceably necessary in $V[G_n]$.  And since
$V_\delta\prec V$ and $\mathbb{P}_n \in V_\delta$,
$V_\delta[G_n]\prec V[G_n]$ by
\lemref{lem_Vdelta_submod_V_forcing_ext}.  By elementarity, $\phi$
is $\Gamma$-forceably necessary in $V_\delta[G_n]$.  But the
iterated forcing $\mathbb{P}$ was defined so that
$\mathbb{P}_{n+1}=\mathbb{P}_n*\dot{\mathbb{Q}}_n$, where
$\mathbb{Q}_n$ forces $\square_\Gamma\phi$ in
$V_\delta^{\mathbb{P}_n}$. So take $V$-generic $G_{n+1} \subset
\mathbb{P}_{n+1}$.  We have $V_\delta[G_{n+1}]\models
\square_\Gamma\phi$.  So by elementarity again, $V[G_{n+1}]
\models \square_\Gamma\phi$.  But $V[G]$ is a $\Gamma$-forcing
extension of $V[G_{n+1}]$, so $V[G] \models \square_\Gamma\phi$.

Finally, since $\mathbb{P}_\omega$ is in $\Gamma$,
$\Box_\Gamma\sigma(\kappa)$ still holds.
\end{proof}

One doesn't need to limit applications of
\thmref{thm_con_MP_Gamma_indestruct_kappa} to large cardinals.

\begin{thm}\label{thm kappa card}
If there is a model of \ZFC\ in which $\kappa$ is any uncountable
cardinal, and $\Gamma(\kappa)$ is a class of forcing notions
preserving the cardinality of $\kappa$, defined using $\kappa$,
which is $\Gamma(\kappa)$-necessarily closed under iteration of
length $\omega$ with appropriate support and which contains
trivial forcing then there is a model of
$\ZFC+\MP_{\Gamma(\kappa)}$.
\end{thm}
\begin{proof}
Let $\sigma(\kappa)$ = ``$\kappa$ is an uncountable cardinal".  So
for $\kappa$, $\Box_{\Gamma(\kappa)}\sigma(\kappa)$.  Now apply
\thmref{thm_con_MP_Gamma_indestruct_kappa}.
\end{proof}

Let $<\kappa-dc$ denote the class of $<\kappa$-directed closed
notions of forcing.

\begin{cor}
The following are equivalent:
\begin{list}
{(\arabic{listcounter})}{\usecounter{listcounter}\setlength{\rightmargin}{\leftmargin}}
\item
$Con(\ZFC$ $+$ ``$\kappa$ is an uncountable cardinal")
\item
$Con(\ZFC$ $+$ $\MP_{<\kappa-dc}$ $+$ ``$\kappa$ is an uncountable
cardinal")
\item
$Con(\ZFC$ $+$ $\MP_{\kappa-closed}$ $+$ ``$\kappa$ is an
uncountable cardinal")
\item
$Con(\ZFC$ $+$ $\MP_{\kappa-cc}$ $+$ ``$\kappa$ is an uncountable
cardinal").
\end{list}
\end{cor}

From existing indestructibility results arise various
equiconsistency theorems, by using the following corollary of
\thmref{thm_con_MP_Gamma_indestruct_kappa}.

\begin{cor}\label{cor1}
Let $\Gamma$ be a class of forcing notions, parametrized by
$\kappa$, which is $\Gamma$-necessarily closed under iterations of
length $\omega$ using appropriate support, and which contains the
trivial forcing.  If there is a model of \ZFC\ in which $\kappa$
is a cardinal and $\sigma(\kappa)$ is a property of $\kappa$ that
can provably be made indestructible under $\Gamma$-forcing, then
there is a model of
$\ZFC+\square_\Gamma\sigma(\kappa)+\MP_{\Gamma}$.
\end{cor}

\begin{proof}
In the forward direction, let $M_0\models \ZFC\ $+$
\sigma(\kappa)$. By the hypothesis, there is also $M_1\models
\ZFC\ +\square_\Gamma\sigma(\kappa)$ holds. But then, by
\thmref{thm_con_MP_Gamma_indestruct_kappa}, there is $M\models
\ZFC+\square_\Gamma\sigma(\kappa)+\MP_{\Gamma}$. In the reverse
direction, if we are given $M\models
\ZFC+\square_\Gamma\sigma(\kappa)+\MP_{\Gamma}$, then $M \models
\ZFC+\sigma(\kappa))$ since $\square_\Gamma\sigma(\kappa)$
directly implies $\sigma(\kappa)$.
\end{proof}

\begin{cor}\label{cor2}
$Con(\ZFC+\ ``$There is a supercompact cardinal $\kappa")$ if and
only if $Con(\ZFC$ $+$ ``There is a supercompact cardinal
$\kappa$, indestructible by
$<\kappa$-directed-closed forcing, such that 
$\MP_{<\kappa-dc}$ holds").
\end{cor}

\begin{proof}
By Laver's indestructibility result in \cite{lav:ind},
\corref{cor1}, and the fact that the class of $< \kappa$-directed
closed forcing notions is closed under $\omega$ iterations using
countable support.
\end{proof}

A further result on indestructibility, in \cite{gs:strong}, gives
a preparation that renders a strong cardinal indestructible under
a class of forcing notions, those that are $\kappa^+$-weakly
closed satisfying the Prikry condition.  Herein I will use a
weakening of this result, that strong cardinals can be made
indestructible under $\leq\kappa$-strategically closed forcing. We
will call such a $\kappa$ \textbf{Gitik-Shelah-indestructible}.

\begin{cor}\label{cor4}

$Con(\ZFC+\ ``$There is a strong cardinal $\kappa")$\linebreak
$\Longleftrightarrow Con(\ZFC+\ ``$There is a strong cardinal
$\kappa$, Gitik-Shelah-indestructible, such that 
$\MP_{\leq\kappa-strat-cl}$ holds").
\end{cor}

\begin{proof}
By the indestructibility result of \cite{gs:strong}, \corref{cor1}
and the fact that the class of $\leq\kappa$-strategically closed
forcing notions is closed under $\omega$ iterations using
countable support.
\end{proof}

\noindent\textit{Remark.}  The class of $\leq\kappa$-strategically
closed forcing notions includes, for example, all the
$\kappa^+$-closed forcing notions.\vspace {12 pt}

Separate from the cases just handled is the following result.
\begin{thm}\label{thm sc indest add one}
$Con(\ZFC+\ ``$There is a strongly compact cardinal
$\kappa")$\linebreak $\Longleftrightarrow Con(\ZFC+\ ``$There is a
strongly compact cardinal $\kappa$, indestructible by
$Add(\kappa,1)$ forcing, such that 
$\MP_{Add(\kappa,1)}$ holds").
\end{thm}

\begin{proof}
By \cite{ham:lot}, a strongly compact cardinal can be made
$Add(\kappa,1)$-indestructible.  In this case, $\Gamma =
\{Add(\kappa,1)\}$, which has just one element. So if anything is
forceably necessary by $\Gamma$, then forcing once with it will
make it happen.  So any forcing extension by $Add(\kappa,1)$ is a
model of $\MP_{Add(\kappa,1)}$.
\end{proof}

\section{${\mathrm{MP}_{\Gamma(\kappa)}}$ with Parameters}

Corresponding to the results in Section~\ref{sec lc indestruct} on
modifications of the principle \MP\ which pertains to statements
without parameters are results analogous to the equiconsistency
result regarding $\MP(H(\omega_1))$--which takes real
parameters--in \cite{ham:max}.  As before, let $\kappa$ be a large
cardinal.  Let $\Gamma(\kappa)$ be a class of forcing notions
definable from $\kappa$ under which the largeness of $\kappa$ can
be made indestructible in some $\Gamma(\kappa)$-forcing extension.
In addition, we now specify a parameter set $S(\kappa)$ to be used
in expressing the maximality principle
$\mathrm{\MP}_{\Gamma(\kappa)}(S(\kappa))$.  If we express the
largeness of $\kappa$ with the formula $\sigma(\kappa)$, the
generalized theories for which equiconsistency is sought are:

\begin{list}
{(\arabic{listcounter})}{\usecounter{listcounter}\setlength{\rightmargin}{\leftmargin}}
\item
$\mathrm{\ZFC}$ $+$ $\Box_{\Gamma(\kappa)}\sigma(\kappa)$ $+$
$\mathrm{\MP}_{\Gamma(\kappa)}(S(\kappa))$
\item
$\mathrm{\ZFC}$ $+$ $V_\delta \prec V$ $+$ ``$\delta$ is
inaccessible" $+$ $\kappa < \delta$ $+$ $\sigma(\kappa)$
\item
$\mathrm{\ZFC}$ $+$ ``$ORD$ is Mahlo" $+$ $\sigma(\kappa)$.
\end{list}

\noindent\textit{Remark.} The principle ``$ORD$ is Mahlo" is known
as the L\'evy Scheme.  It says that every club class in $ORD$
contains a regular cardinal.

Theory (1) asserts that the largeness of $\kappa$ is
indestructible under $\Gamma(\kappa)$-forcing, and that any
statement with parameter in $S(\kappa)$ which is
$\Gamma(\kappa)$-forceably necessary is
$\Gamma(\kappa)$-necessary.  Theories (2) and (3) are identical to
those in \cite{ham:max} with the addition of the largeness of
$\kappa$.  Indeed, the proof that (2) and (3) are equiconsistent
follows that in \cite{ham:max}; the forward implication is direct:
let $C \subseteq ORD$ be any definable club; $C\cap \delta$ is
unbounded in $\delta$ so $\delta$, which is regular, is in $C$.
The converse is proved by constructing a model through a
compactness argument using L\'evy Reflection, which can be done by
including $\kappa$ in an expanded ground model and suitably
modifying the argument of
\lemref{lem_Vdelta_submod_V_in_elem_ext}.

An unsolved problem for the moment is to prove that
$Con(1)\Longrightarrow Con(2)$ for any specific kind of large
cardinal.  The analogous proof in \cite{ham:max} uses $L[x]$, for
real $x$, as an model of (2) within a model of (1).  But here we
need to preserve the largeness of $\kappa$ when we move to an
inner model which in fact becomes the core model.  Results for
strong cardinals have been developed in this model, so it may work
here.  On the other hand, as of yet, supercompact cardinals have
no core model, so for them the problem will probably be more
difficult.

At any rate, the best we can do here is prove the direction
$Con(2)\Longrightarrow Con(1)$ for specific cases, which at least
gives an upper bound on the consistency strength of
$\MP_{\Gamma(\kappa)}(S(\kappa))$, namely, that of a cardinal
satisfying property $\sigma$ together with the L\'evy scheme.
Rather than prove a general version of this, we will prove it for
specific cases, addressing their specific issues.  The first
result will bound the consistency strength of the principle
$\MP_{<\kappa-dc}(H(\kappa^+))$, where $\kappa$ is a supercompact
cardinal and $<\kappa-dc$ is the class of $<\kappa$-directed
closed forcing notions.  We will consider the relative consistency
of the following theories:

\begin{list}
{(\arabic{listcounter})}{\usecounter{listcounter}\setlength{\rightmargin}{\leftmargin}}
\item
$\mathrm{\ZFC}$ $+$ ``there is a supercompact $\kappa$,
indestructible by $<\kappa$-directed closed forcing" $+$
$\mathrm{\MP}_{<\kappa-dc}(H(\kappa^+))$
\item
$\mathrm{\ZFC}$ $+$ ``$\kappa$ is supercompact" $+$ $V_\delta
\prec V$ $+$ $\kappa < \delta$ $+$ ``$\delta$ is inaccessible"
\item
$\mathrm{\ZFC}$ $+$ ``$\kappa$ is supercompact" $+$ ``$ORD$ is
Mahlo".
\end{list}

We will prove the implication $Con(2)\Longrightarrow Con(1)$ by
constructing a model of theory (1) via a $\delta$ iteration,
beginning with a model of (2) as ground model. But first we need
the following standard lemma to assure us that this iteration does
not collapse $\delta$.

\begin{lem}\label{lem delta system}
Let $\delta$ be inaccessible, $\kappa < \delta$, and $\mathbb{P} =
\mathbb{P}_\delta$ a $\delta$-stage, $<\kappa$-support iterated
notion of forcing such that for all $\alpha < \delta$,
$\mathbb{P}_\alpha \Vdash ``\dot{\mathbb{Q}}_\alpha$ is in
$V_\delta$" where
$\mathbb{P}_{\alpha+1}=\mathbb{P}_\alpha*\dot{\mathbb{Q}}_\alpha$.
Then $\mathbb{P}$ is $\delta$-cc.
\end{lem}

\begin{proof}
For all $\alpha<\delta$, $\mathbb{P}_\alpha$ is in $V_\delta$ and
therefore is $\delta$-cc.  This is true by induction: for
successor stages each $\dot{\mathbb{Q}}_\alpha$ is in $V_\delta$,
and at limit stages, $\mathbb{P}_\alpha$ is in $V_\delta$ because
$\delta$ is inaccessible.  It remains to show that
$\mathbb{P}_\delta$ itself is $\delta$-cc.  Suppose, towards
contradiction, that $A=\{p^\beta|\beta<\delta\}$ is an antichain
in $\mathbb{P}_\delta$. The iteration is with $<\kappa$-support.
Thus, since $\kappa < \delta$ and $\delta$ is inaccessible,
$|\delta^{<\kappa}|<\delta$. So by passing to a subset $B
\subseteq A$ where $|B| = \delta$, we can assume that
$\{\supp(p^\beta)|\beta <\delta\}$ forms a $\Delta$-system. Let
$r$ be the root of this $\Delta$-system.  Since $\delta$ is
inaccessible, we can fix $\zeta < \delta$ with $r \subset \zeta$.

For any $p$ and $p'$ in $B$, we have $p$ is incompatible with $p'$
as $B$ is an antichain.  However, it is possible to choose such
$p$ and $p'$ so that $p\upharpoonright_\zeta$ is compatible with
$p'\upharpoonright_\zeta$, since for all $\zeta < \delta$,
$\mathbb{P}_\zeta$ has the $\delta$-cc.  We have $supp(p)\cap
supp(p')=r\subset\zeta$.  But by Lemma 5.11(f) of \cite{kun:sth},
Chapter VIII, this implies that $p$ is incompatible with $p'$ if
and only if $p\upharpoonright_\zeta$ is incompatible with
$p'\upharpoonright_\zeta$, a contradiction.
\end{proof}

\begin{thm}\label{thm sc indest MP <kdc H}
If there is a model of \ZFC\ $+$ ``$\kappa$ is a supercompact
cardinal" $+$ $V_\delta \prec V$ $+$ ``$\delta$ is inaccessible"
$+$ $\kappa<\delta$, then there is a model of \ZFC\ $+$ ``$\kappa$
is a supercompact cardinal indestructible by $<\kappa$-directed
closed forcing" $+$ $\mathrm{\MP}_{<\kappa -dc}(H(\kappa^+))$.
\end{thm}
\begin{proof}
Suppose $V\models \ZFC$ $+$ $V_\delta\prec V$ $+$ ``$\kappa$ is
supercompact" $+$ $\kappa<\delta$ $+$ ``$\delta$ is inaccessible".
In $V$, by the well-known result of Laver \cite{lav:ind}, there is
a $<\kappa$-directed closed notion of forcing $\tilde{\mathbb{P}}$
that produces a model in which the supercompactness of $\kappa$ is
indestructible by further $<\kappa$-directed closed forcing. Since
$\tilde{\mathbb{P}}$ has size $\kappa<\delta$,
$\tilde{\mathbb{P}}$ is $\delta-cc$, so it preserves $\delta$ and
its inaccessibility, together with the statement $\kappa <\delta$.

Let $G$ be $V$-generic over $\tilde{\mathbb{P}}$.  We now work in
$V[G]$.  By \lemref{lem_Vdelta_submod_V_forcing_ext}, we have that
$V_\delta[G]\prec V[G]$.  We will construct a $<\kappa$-support,
$\delta$-iteration $\mathbb{P} = \mathbb{P}_\delta$ of
$<\kappa$-directed closed forcing notions as follows. Let $\Gamma$
be the class of $<\kappa$-directed closed forcing notions.  At
stage $\alpha$ let $\alpha$ encode the triple
$\langle\beta,\gamma,n\rangle$, where $\beta < \alpha$, $\gamma <
\delta$, and $n \in \omega$. Consider $\phi(z)$, where
$\phi=\phi_n$, $z$ is the $\gamma^{th}$ element of
$H(\delta)^{V^{\mathbb{P}_\beta}}$ for $\beta < \alpha$.  If
$\phi(z)$ is $\Gamma$-forceably necessary over $V_\delta[G]$,
then, in $V_\delta[G]^{\mathbb{P}_\alpha}$, let
$\dot{\mathbb{Q}}_\alpha$ be the $\mathbb{P}_\alpha$-name of a
forcing notion forcing that $\phi(z)$ is necessary, otherwise let
$\dot{\mathbb{Q}}_\alpha$ be trivial forcing.  This defines
$\mathbb{P} = \mathbb{P}_\delta$.  Notice that since every
$\dot{\mathbb{Q}}_\alpha$ is in $V_\delta[G]^{\mathbb{P}_\alpha}$,
$|\dot{\mathbb{Q}}_\alpha| < \delta$, so $\dot{\mathbb{Q}}_\alpha$
is $\delta-cc$.  $\mathbb{P}$ is an iteration of these with
$<\kappa$ support, so by \lemref{lem delta system}, it is also
$\delta-cc$. So $\delta$ is not collapsed in the extension, and
its cofinality is preserved.

Let $H$ be $V[G]$ generic over $\mathbb{P}$.  Since $\mathbb{P}$
is a $\delta$ iteration with $<\kappa$ support of
$<\kappa$-directed closed forcing notions, $\mathbb{P}$ is itself
$<\kappa$-directed closed. And the preceding Laver Preparation
$\tilde{\mathbb{P}}$ ensures that $\kappa$ is still indestructibly
supercompact in $V[G][H]$.  The final claim is that
$V[G][H]\models \mathrm{\MP}_{\Gamma}(H(\delta))$. To show this,
suppose $V[G][H]$ satisfies that $\phi(z)$ is $\Gamma$-forceably
necessary, where $z$ is in $H(\delta)$. Since $|z| <
\cof(\delta)=\delta$, $z$ will have been introduced as a
$\mathbb{P}_\beta$-name at some stage $\beta < \delta$, say, as
the $\gamma^{th}$ element of $V^{\mathbb{P}_\beta}$, and
$\phi=\phi_n$ for some $n$ in $\omega$.  So there is a stage
$\alpha<\delta$ where $\alpha=\langle \beta,\gamma,n\rangle$.
$V[G][H_\alpha]$ also satisfies that $\phi(z)$ is
$\Gamma$-forceably necessary, where $H_\alpha$ is taken to be
$\mathbb{P}_\alpha$ generic over $V[G]$ and
$V[G][H]=V[G][H_\alpha][H_{TAIL}]$ via the factoring
$\mathbb{P}=\mathbb{P}_\alpha*\dot{\mathbb{P}}_{TAIL}$. And by
elementarity, $V_\delta[G][H_\alpha]$ satisfies that $\phi(z)$ is
$\Gamma$-forceably necessary.  But by the definition of
$\mathbb{P}$,
$\mathbb{P}_{\alpha+1}=\mathbb{P}_\alpha\ast\dot{\mathbb{Q}}_\alpha$,
where $V_\delta[G][H_\alpha]\models$ ``$\dot{\mathbb{Q}}_\alpha$
forces $\phi(z)$ to be $\Gamma$-necessary".  This implies that by
refactoring
$\mathbb{P}=\mathbb{P}_{\alpha+1}*\dot{\mathbb{P}}_{TAIL2}$ and
taking $H_{\alpha+1}$ to be $\mathbb{P}_{\alpha+1}$ generic over
$V[G]$, $V_\delta[G][H_{\alpha+1}]\models$ ``$\phi(z)$ is
$\Gamma$-necessary".  Elementarity then gives
$V[G][H_{\alpha+1}]\models$ ``$\phi(z)$ is $\Gamma$-necessary",
from which we finally have $V[G][H]\models$ ``$\phi(z)$ is
$\Gamma$-necessary".  This gives $V[G][H]\models
\mathrm{\MP}_{\Gamma}(H(\delta))$.  Notice that $\mathbb{P}$ has
$<\kappa$-directed closed factors that preserve cardinals
$<\kappa$ but eventually every ordinal between $\kappa$ and
$\delta$ appears as a parameter in
$H(\delta)^{V^{\mathbb{P}_\beta}}$ at some stage $\alpha$ where
$\beta < \alpha<\delta$. So $\delta$ becomes $\kappa^+$ in $V[G]$.
This finally gives $V[G][H]\models
\mathrm{\MP}_{\Gamma}(H(\kappa^+))$
\end{proof}

\begin{cor}\label{cor sc indest le ord mahlo}
If there is a model of a supercompact cardinal where the L\'evy
Scheme holds, then there is a model of a Laver-indestructible
supercompact cardinal where
$\mathrm{\MP}_{<\kappa-dc}(H(\kappa^+))$ holds.
\end{cor}

\begin{cor}\label{cor sc indest lt mahlo}
If there is a supercompact cardinal with a Mahlo cardinal above
it, then there is a model of a Laver-indestructible supercompact
cardinal where $\mathrm{\MP}_{<\kappa-dc}(H(\kappa^+))$ holds.
\end{cor}

\begin{con}
The consistency strength of the theory \ZFC\ $+$ ``there is a
supercompact $\kappa$, indestructible by $<\kappa$-directed closed
forcing" $+$ \linebreak $\mathrm{\MP}_{<\kappa-dc}(H(\kappa^+))$
is the same as that of the theory \ZFC\ $+$ ``$\kappa$ is
supercompact" $+$ $V_\delta \prec V$ $+$ $\kappa < \delta$ $+$
``$\delta$ is inaccessible".
\end{con}

Looking next at the Gitik-Shelah result of \cite{gs:strong} for
strong cardinals, we see that the previous argument adapts to
handle this case.  Again, we will consider the class
$\Gamma_{\leq\kappa-strat-cl}$ of $\leq\kappa$-strategically
closed forcing notions, a strengthening of the Prikry condition
Gitik and Shelah impose. An important difference between this
class and the class of $<\kappa$-directed closed forcing notions
is that the most general parameter set of
$\MP_{\leq\kappa-strat-cl}(X)$
 is $H(\kappa^{++})$, the sets
of hereditary size less than $\kappa^{++}$, by the usual
consideration: one can regard a parameter from this set as encoded
by a subset of $\kappa^{+}$, any larger parameter set will have
members of size greater than $\kappa^{+}$.  Suppose , toward
contradiction, that $|x| > \kappa^+$ for some parameter $x$ in a
parameter set $X$.  Such a parameter $x$ can have its cardinality
collapsed by $\leq\kappa$-strategically closed forcing.  So the
statement $|x| \leq \kappa^+$ is
$\Gamma_{\leq\kappa-strat-cl}$-forceably necessary. By applying
$\mathrm{\MP}_{\leq\kappa-strat-cl}(X)$, this falsifies the
statement that $|x|
> \kappa^+$.

\begin{thm}
If there is a model of \ZFC\ $+$ ``$\kappa$ is a strong cardinal"
$+$ $V_\delta \prec V$ $+$ ``$\delta$ is inaccessible" $+$
$\kappa<\delta$, then there is a model of \ZFC\ $+$ ``$\kappa$ is
a strong cardinal indestructible by $\leq\kappa$-strategically
closed forcing" $+$
$\mathrm{\MP}_{\leq\kappa-strat-cl}(H(\kappa^{++}))$.
\end{thm}
\begin{proof}
Here we emulate the previous proof, performing the necessary
forcing preparation $\tilde{\mathbb{P}}$ in a model of theory (2)
(\ZFC\ $+$ ``$\kappa$ is a strong cardinal" $+$ $V_\delta \prec V$
$+$ ``$\delta$ is inaccessible" $+$ $\kappa<\delta$).  By
elementarity, the construction can be assumed to take place in
$V_\delta$.  This means $\tilde{\mathbb{P}}$ has cardinality below
$\delta$, and is therefore $\delta$-cc and the resulting forcing
extension $V[G]$, where $G$ is $M$-generic over
$\tilde{\mathbb{P}}$, preserves theory (2).  Of course,
indestructibility of the strongness of $\kappa$ also holds in the
extension.

Let $G$ be $V$-generic over $\tilde{\mathbb{P}}$.  We now work in
$V[G]$.  By \lemref{lem_Vdelta_submod_V_forcing_ext}, we have that
$V_\delta[G]\prec V[G]$.  We will construct a
$\leq\kappa$-support, $\delta$-iteration $\mathbb{P} =
\mathbb{P}_\delta$ of $\leq\kappa$-strategically closed forcing
notions as follows. At stage $\alpha$ let $\alpha$ encode the
triple $\langle\beta,\gamma,n\rangle$, where $\beta < \alpha$,
$\gamma < \delta$, and $n \in \omega$. Consider $\phi(z)$, where
$\phi=\phi_n$, $z$ is the $\gamma^{th}$ element of
$H(\delta)^{V^{\mathbb{P}_\beta}}$ for $\beta < \alpha$. If
$\phi(z)$ is $\Gamma_{\leq\kappa-strat-cl}$-forceably necessary
over $V_\delta[G]$, then, in $V_\delta[G]^{\mathbb{P}_\alpha}$,
let $\dot{\mathbb{Q}}_\alpha$ be the $\mathbb{P}_\alpha$-name of a
$\leq\kappa$-strategically closed forcing notion forcing that
$\phi(z)$ is necessary, otherwise let $\dot{\mathbb{Q}}_\alpha$ be
trivial forcing.  This defines $\mathbb{P} = \mathbb{P}_\delta$.
Notice that since every $\dot{\mathbb{Q}}_\alpha$ is in
$V_\delta[G]^{\mathbb{P}_\alpha}$, $|\dot{\mathbb{Q}}_\alpha| <
\delta$, so $\dot{\mathbb{Q}}_\alpha$ is $\delta-cc$. $\mathbb{P}$
is an iteration of these with $\leq\kappa$ support, so by
\lemref{lem delta system}, using $\kappa^+$ in place of $\kappa$,
it is also $\delta-cc$. So $\delta$ is not collapsed in the
extension, and its cofinality is preserved.

Let $H$ be $V[G]$ generic over $\mathbb{P}$. An argument similar
to that of the previous theorem shows that $V[G][H]\models
\mathrm{\MP}_{\leq\kappa-strat-cl}(H(\kappa^{++}))$.
\end{proof}

\noindent\textit{Remark.}   Notice that $\mathbb{P}$, in this
proof, contains $\leq\kappa$-strategically closed forcing factors
that preserve cardinals $<\kappa^+$ but eventually collapse all
cardinals between $\kappa^+$ and $\delta$.  So $\delta$ becomes
$\kappa^{++}$ in $V[G]$.

\begin{cor}\label{}
If there is a model of a strong cardinal where the L\'evy Scheme
holds, then there is a model of a Gitik-Shelah-indestructible
strong cardinal $\kappa$ where
$\mathrm{\MP}_{\leq\kappa-strat-cl}(H(\kappa^{++}))$ holds.
\end{cor}

\begin{cor}\label{cor str indest lt mahlo}
If there is a strong cardinal with Mahlo cardinal above it, then
there is a model of a Gitik-Shelah-indestructible strong cardinal
$\kappa$ where $\mathrm{\MP}_
{\leq\kappa-strat-cl}(H(\kappa^{++}))$ holds.
\end{cor}

\begin{con} The theory \ZFC\ $+$ ``$\kappa$ is a strong cardinal,
indestructible by $\leq\kappa$-strategically closed forcing" $+$
$\mathrm{\MP}_{\leq\kappa-strat-cl}(H(\kappa^{++}))$ is the same
consistency strength as the theory \ZFC\ $+$ ``$\kappa$ is strong"
$+$ $V_\delta \prec V$ $+$ $\kappa < \delta$ $+$ ``$\delta$ is
inaccessible".
\end{con}

Other results could follow in this vein.  Hamkins' own results on
establishing indestructibility of large cardinals by the Lottery
preparation in \cite{ham:lot} should be be usable here. By
performing all such preparations in $V_\delta$, they have the
$\delta-cc$ property. All that is required is closure under
appropriate support $\delta$-iterations of the class of forcing
notions for which indestructibility is asserted.  The resulting
forcing notion will be $\delta$-cc, preserving $\delta$, and a
generic extension over this forcing notion will satisfy the
corresponding maximality principle.

\backmatter

\end{document}